\documentclass[10pt]{article}
\usepackage{CJKutf8}    
\usepackage[backend=biber, style=numeric-comp, sorting=none, maxnames=2, minnames=1,uniquename=false, giveninits=false]{biblatex}

\usepackage{titlesec}   
\usepackage{booktabs}   
\usepackage{threeparttable} 
\usepackage{tabularx}   
\usepackage{multirow}   
\usepackage{graphicx}   
\usepackage{caption}    
\usepackage{subcaption} 
\usepackage{float}      
\usepackage{placeins}   
\usepackage{geometry}   
\usepackage{times}      
\usepackage{hyperref}   
\usepackage{enumitem}   

\usepackage[most]{tcolorbox}
\tcbuselibrary{breakable, skins}

\usepackage{amsmath,amsfonts,amsthm}
\usepackage{bm}
\usepackage{mathrsfs}

\usepackage{algorithm}
\usepackage{algorithmicx}
\usepackage{algpseudocode}

\usepackage{cleveref}   

\newcommand{\Title}[1]{
            \par\vspace{2em}\noindent
            \begin{minipage}{\textwidth}
                \centering\LARGE\textbf{#1}
            \end{minipage}
            \par\vspace{2em}}
\newcommand{\Authorlist}[1]{
            \vspace{1em}\noindent
            \begin{minipage}{\textwidth}
                \large\textbf{#1}
            \end{minipage}
            \vspace{1em}}
\newcommand{\Author}[2]{{#1}\textsuperscript{#2}}
\newcommand{\Affiliation}[1]{
            \par\noindent
            \begin{minipage}{\textwidth}
                \large\textit{#1}
            \end{minipage}
            \par\vspace{0.5em}}
\newcommand{\Coauthor}[1]{
            \par\vspace{1em}\noindent
            \begin{minipage}{\textwidth}
                \large{†The authors contribute equally to this work. #1} 
            \end{minipage}
            \par\vspace{0.25em}}
\newcommand{\Corresponding}[1]{
            \par\vspace{1em}\noindent
            \begin{minipage}{\textwidth}
                \large{*Corresponding author. #1} 
            \end{minipage}
            \par\vspace{0.25em}}
\newcommand{\Email}[1]{
            \par\vspace{0.25em}\noindent
            \begin{minipage}{\textwidth}
                \large{E-mail address: #1} 
            \end{minipage}
            \par\vspace{1em}}
\newcommand{\Abstract}[1]{
            \par\vspace{1em}\noindent
            \begin{minipage}{\textwidth}
                \setlength{\leftskip}{0.5cm}
                \setlength{\rightskip}{0.5cm}
                \textbf{Abstract:} #1
            \end{minipage}
            \par\vspace{1em}}
\newcommand{\Keywords}[1]{
            \par\vspace{0.5em}\noindent
            \begin{minipage}{\textwidth}
                \setlength{\leftskip}{0.5cm}
                \setlength{\rightskip}{0.5cm}
                \textbf{Keywords:} #1
            \end{minipage}
            \par\vspace{0.5em}}

\titleformat*{\section}{\large\bfseries}
\titlespacing*{\section}{0pt}{6ex plus 3ex minus.2ex}{3ex plus.1ex}
\titleformat*{\subsection}{\bfseries\itshape}
\titlespacing*{\subsection}{0pt}{4.5ex plus 2ex minus.1ex}{3ex plus.1ex}
\titleformat*{\subsubsection}{\bfseries\itshape}
\titlespacing*{\subsubsection}{0pt}{4.5ex plus 2ex minus.1ex}{3ex plus.1ex}
\hypersetup{
    colorlinks=true,
    linkcolor=blue,    
    citecolor=blue,    
    urlcolor=blue,      
    hypertexnames=false
}

\Crefname{figure}{Fig.}{Figs.}
\Crefname{table}{Tab.}{Ts.}
\Crefname{section}{Sec.}{Secs.}

\title{}
\author{}
\date{}

\begin{document}

\Title{A Systematic Analysis of Automatic Differentiation versus Discretization-based Constraints for Physics-Informed PDE Solvers}
\Authorlist{\Author{Xing Guo}{a, b, †}, \Author{Hongwei Tang}{b, †}, \Author{Zewei Meng}{b}, \Author{Yidong Zhang}{a, b}, \Author{Shaoqiu Xiao}{a,*}, \Author{Feng Liu}{b,*}}
\Affiliation{a. School of Systems Science and Engineering, Sun Yat-sen University, Guangzhou, China}
\Affiliation{b. National Key Laboratory of Aerospace Physics in Fluids, Mianyang, China}
\Coauthor{}
\Corresponding{}
\Email{\href{guox89@mail2.sysu.edu.cn}{guox89@mail2.sysu.edu.cn} (X. Guo), 
        \href{thw1021@nuaa.edu.cn}{thw1021@nuaa.edu.cn} (H. Tang), 
        \href{mengzw94@163.com}{mengzw94@163.com} (Z. Meng), 
        \href{zhangyd8@mail2.sysu.edu.cn}{zhangyd8@mail2.sysu.edu.cn} (Y. Zhang), 
        \href{xiaoshq8@mail.sysu.edu.cn}{xiaoshq8@mail.sysu.edu.cn} (S. Xiao),
        \href{liufengmaple@foxmail.com}{liufengmaple@foxmail.com} (F. Liu)}
\Abstract{Physics-informed neural networks (PINNs) represent a growing frontier in using artificial intelligence to solve partial differential equations (PDEs). Automatic differentiation (AD) plays a central role in this paradigm, which is mesh-free and replaces traditional iterative solvers with gradient-based optimization in continuous space. However, the inherent limitations of AD, particularly in handling higher-order derivatives and discontinuous solutions, pose significant challenges for complex problems. This has motivated a growing number of researchers to explore discretization-based constraints as an alternative path. Yet, the respective applicability of these two paradigms remains largely unexplored. In this work, we conduct systematic experiments across a wide spectrum of problems, from simple linear Poisson to high-Mach hypersonic flows with strong discontinuities. Through a rigorous decomposition of approximation, optimization, and truncation errors, we systematically elucidate the fundamental trade-offs and error-governing mechanisms of both paradigms, as well as two representative network architectures: multi-layer perceptron (MLP) and graph neural network (GNN). Our results reveal a consistent trend: as nonlinearity strengthens, the accuracy advantage of discretization-based constraints becomes increasingly pronounced, with smaller optimization errors compensating for the truncation errors. Moreover, the more complex the nonlinearity and boundary conditions, the greater the advantage of GNN over MLP. These insights offer a robust practical guideline for configuring neural PDE solvers in demanding engineering applications. Our source data and code are available at \href{https://github.com/guoxing0809/neuropde_analysis}{https://github.com/guoxing0809/neuropde\_analysis}.}

\section{Introduction}\label{sec:introduction}
Physics-informed neural networks (PINNs) have become an established methodology since Raissi’s seminal work \parencite{2019_jcp_vanillapinn}, evolving from earlier attempts to integrate artificial neural networks with physical constraints for solving differential equations. The PINN framework advances this paradigm by preserving physical interpretability while incorporating deep learning’s representational capacity, demonstrating superior effectiveness for ill-posed \parencite{2025_illcondition} and inverse problems \parencite{2020_cmame_conservativepinn} where purely physics-driven or data-driven methods face limitations, and exhibiting enhanced robustness to noisy observational data \parencite{2020_hiddenfluid}. In fluid dynamics, PINN variants have branched into multiple directions, reflecting different design choices in data sampling, network architecture and the formulation of physical constraints \parencite{2022_jsc_pinnreview, 2025_air_pinnreview}.

Unlike direct numerical simulation (DNS), which exhibits well-defined convergence orders and rigorous error bounds, PINNs rely on a nonlinear optimization framework with mean squared error (MSE) as the loss function. Consequently, PINNs are inherently subject to the so-called spectral bias \parencite{2019_icml_spectralbias}---a tendency to prioritize learning low-frequency components of the solution while struggling to capture high-frequency features. This pathology, which often hinders the training of PINNs in more complex cases, has been systematically analyzed by \parencite{2022_jcp_pinnfail} through the lens of neural tangent kernel theory. Several strategies have been adopted by researchers to mitigate such limitations. 
Dynamic loss weighting \parencite{2021_nn_dualdimer, 2026_nd_csdpinn} is a typical example. In addition, a number of techniques have been developed to improve the optimization convergence of PINNs, including variable-scaling techniques for stiff and high-frequency problems \parencite{2025_jcp_vspinn}, gradient-enhanced residual formulation \parencite{2022_cmame_gpinn}, and curvature-aware optimization method \parencite{2026_arxiv_curvatureaware}. More recently, an operator learning framework has been proposed to learn the Cholesky factors of regularized Newton updates, unrolling the iteration into a neural architecture to emulate partial differential equations (PDEs) solvers \parencite{2026_jcp_operatorlearning}. 
Although these advances have improved PINN training in certain settings, their broader practical applicability remains limited. Dynamic loss weighting strategies, for instance, are highly problem-dependent and often lack generalizability across different problem settings. Furthermore, higher-order optimizers, despite reported accuracy gains in certain cases, often suffer from convergence instability \parencite{2024_nips_pinnacle}. Methods that incorporate higher-order parameter-space information into neural networks may reduce the number of outer-loop iterations, but this comes at a computational cost that can be several times that of standard iterative procedures, and their reliability on more complex physical problems remains unclear. Such efforts, however, may diverge from the original motivation for applying machine learning to PDEs. Rather than pursuing marginal improvements over classical solvers on problems where the latter already perform well, neural network methods should focus on more demanding problems where traditional methods fall short \parencite{2025_nmi_betterbenchamarks}.

The alternative physics-informed paradigm embeds the PDE residual through discretization-based constraints, which can take differential, integral, variational, or mixed forms. A notable strength of this paradigm lies in its ability to naturally incorporate elemental conservation laws into the learning framework \parencite{2020_cmame_conservativepinn}. While the conventional finite difference (FD) form is often constrained to regular grids \parencite{2025_jsc_automaticdiff}, the integral finite volume (FV) form can be straightforwardly extended to computational domains with unstructured meshes \parencite{2021_cmame_discretizationnet, 2022_icml_physicsaware, 2024_pof_fvmpignn, 2025_jcp_fvmpignn}. Variational formulations offer another route within the discretization-based constraints. For instance, hp-VPINNs introduced domain decomposition and high-order polynomial test functions to enforce PDE residuals in an integral sense \parencite{2021_cmame_hpvipinn}. Along similar lines, finite element (FE) frameworks have been employed to integrate weak-form discretization with neural network architecture \parencite{2022_cmame_graphneuralgalerkin, 2023_jsc_meshinformed}. Hybrid strategies with mixed forms include coupling automatic differentiation (AD) with auxiliary discrete operators \parencite{2024_nd_mixeddiff}, embedding Rankine-Hugoniot jump conditions as additional loss constraints for discontinuous shock problems \parencite{2024_jsc_discontinuous}, combining spectral discretization for fractional derivatives with AD for integer-order derivatives \parencite{2025_nd_spectralfpinn}, and leveraging prior field information from pre-trained PINNs to enforce conservation laws on node topology networks \parencite{2026_jcp_pginns}. In addition to constraint formulation, the network architecture itself also plays a critical role in determining solver performance. Multi-layer perceptrons (MLPs) are widely used due to their simplicity and efficiency, but they lack the coupling of neighbourhood information, which limits their ability to capture local interactions in complex geometries. Convolutional neural networks (CNNs) have been employed to address this limitation by leveraging local stencils on structured grids \parencite{2021_cmame_discretizationnet, 2026_arxiv_bowshock}, yet they offer limited flexibility for enforcing hard constraints, especially on unstructured meshes. As a result, increasing attention has turned to graph neural networks (GNNs), which offer greater flexibility by operating directly on graph structures \parencite{2022_cmame_graphneuralgalerkin, 2024_pof_fvmpignn, 2025_jcp_fvmpignn}. This flexibility has been demonstrated in hybrid frameworks that couple GNNs with differentiable PDE solvers \parencite{2020_icml_differentiablepdes}, as well as in mesh-based simulation pipelines that learn dynamics directly on adaptive meshes \parencite{2021_iclr_meshbasedgnn}. A practical limitation, however, is that the cost of message passing scales with the number of mesh edges, which can become significant for large meshes. Thus, it is crucial to determine when the accuracy gains of GNNs outweigh their computational overhead, and whether such trade-offs are justified for a given problem.

Much of the existing work has primarily emphasized accuracy, or has sought to reproduce the efficiency of high-order iterative schemes within neural network training, yet their performance on unseen problems remains highly uncertain. Such post-hoc validation is often neither feasible nor reliable for engineering practice. More significantly, according to a report from \parencite{2024_nmi_weakbaselines}, nearly 80\% of machine learning approaches for fluid-related PDEs benchmark their proposed models against weak baselines. This striking statistic reveals persistent overoptimism in machine learning applications for this domain. Accuracy, though important, is not the sole criterion for practical utility. What is equally, if not more, needed in engineering practice are guidelines grounded in physical priors, such as actual efficiency (not the one-sided efficiency that ignores additional overhead), robustness, stability, and the applicability domain of each method \parencite{2025_nmi_betterbenchamarks}. These considerations, however, remain severely underexplored in the current literature. 
For instance, \parencite{2025_jsc_automaticdiff} compares AD with FD on simple linear and weakly nonlinear problems, concluding that AD outperforms FD in training. This conclusion is drawn from a relatively narrow set of test cases and may not generalize well for several reasons: weak nonlinearity favors AD by design, higher-order FD schemes were not tested, and lower training loss does not guarantee fidelity in high-frequency regions. Convergence of neural network training loss alone, without a posteriori verification, is an insufficient metric for evaluating solver performance in engineering practice \parencite{2021_nips_pinnfailure}. Recent work formalises this optimisation pathology under the term the term "hallucination": optimisers can attain low training loss while generating physically inadmissible solutions, meaning that additional safeguarding mechanisms become necessary particularly for strongly-nonlinear problems \parencite{2025_arxiv_hallucinations}.
Nevertheless, broadly‑applicable evaluation criteria decoupled from specific solver formulations and training pipelines are still lacking within the community, which motivates further investigation into reliable assessment strategies for neural PDE solvers.
To this end, we take a step toward filling this gap by systematically and fairly comparing AD and discretization-based formulations. Our goal is not to declare a winner, but to provide practical insights into when and why each paradigm is preferable, and to offer physically informed guidelines for engineering applications.
In this work, our major contributions are generalized as follows:
\begin{itemize}
    \item We conduct a systematic comparative study of AD and discretization-based formulations across a wide spectrum of PDEs with increasing nonlinearity and boundary complexity, revealing that the accuracy advantage of discretization-based methods becomes increasingly pronounced as nonlinearity intensifies, while the approximation and optimization errors of AD escalate more rapidly in strongly nonlinear regimes.
    \item We analyze errors through a rigorous decomposition framework that distinguishes approximation, discretization, and optimization errors, using the numerical solution obtained under identical discretization settings as a reference for truncation error. This enables a clear separation of the neural solver's approximation error from the underlying grid discretization error.
    \item We provide a comprehensive and fair evaluation of different paradigms and network architectures across diverse scenarios, covering accuracy, efficiency, and robustness, and derive consistent empirical guidelines for selecting appropriate neural solver configurations in engineering practice. Our analysis further reveals key practical bottlenecks for achieving robust training under general hyperparameter configurations without extensive case-specific tuning or a posteriori reference solutions.
\end{itemize}
\FloatBarrier

\section{Method}\label{sec:method}
Consider a steady state PDE with general formulation illustrated as
\begin{equation}\label{eq:general_pde}
    G(\mathbf{x}, \mathbf{u}, \nabla \mathbf{u}, \nabla^{2} \mathbf{u}, \cdots) = f(\mathbf{x}), \quad \mathbf{x} \in \Omega,
\end{equation}
\noindent where $\Omega \subset \mathbb{R}^d$ is the spatial domain, $f: \Omega \to \mathbb{R}$ denotes the source term, and $\mathbf{u} = u(\mathbf{x}): \Omega \to \mathbb{R}^n$ is the unknown vector-valued solution with spatial gradient $\nabla \mathbf{u}$. The governing equation $G$ is subject to boundary conditions $\mathcal{B}(\mathbf{x}, \mathbf{u}, \nabla \mathbf{u}, \cdots) = h(\mathbf{x})$ for $\mathbf{x} \in \partial\Omega$. Let $\hat{\mathbf{u}}$ be a trial solution approximated by a neural network $N(\mathbf{x}; \boldsymbol{\theta})$ with trainable parameters $\boldsymbol{\theta}$. Then we reformulate the problem as an optimization task:
\begin{equation}\label{eq:pde_constraints}
\min_{\boldsymbol{\theta}} \left| G(\mathbf{x}, \hat{\mathbf{u}}, \nabla \hat{\mathbf{u}}, \nabla^{2} \hat{\mathbf{u}}, \cdots) - f(\mathbf{x}) \right|^2 \quad \text{s.t.} \quad \mathcal{B}(\mathbf{x}, \hat{\mathbf{u}}, \nabla \hat{\mathbf{u}}, \cdots) = h(\mathbf{x}).
\end{equation}
In PINN training, this constrained problem is reduced to an unconstrained form in two ways. The first is to embed hard constraints directly into the trial function. For example, Dirichlet boundaries can be formulated as
\begin{equation}\label{eq:hard_bc_formulation}
    \hat{\mathbf{u}} = A(\mathbf{x}) + \Psi\big( \mathbf{x}, N(\mathbf{x};\boldsymbol{\theta}) \big), 
\end{equation}
where $A$ satisfies $\mathcal{B}$ exactly and $\Psi$ nullifies boundary contributions \parencite{1998_itnn_annode}. However, such hard constraint construction is not always feasible for problems with complex boundary conditions, though some works introduce auxiliary extra learnable fields for formal hard Neumann-type enforcement \parencite{2022_nips_hardbc}, yet they only shift soft-constraint penalties into domain-wide consistency terms, still susceptible to derivative-induced solution contamination. In such cases, the second approach implements soft constraints through additional penalty terms:
\begin{equation}\label{eq:pde_optimization}
\min_{\boldsymbol{\theta}} (| G - f |^2 + | \mathcal{B} - h |^2).
\end{equation}

\begin{figure}[!ht]
    \centering
    \includegraphics[width=1\linewidth]{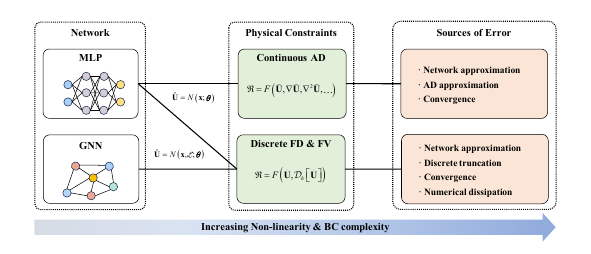}
    \caption{An overview of our analysis framework, covering network architectures, constraint formulations, error sources, and increasing problem complexity.}
    \label{fig:analysis_framework}
\end{figure}

This multi-objective optimization is typically solved by minimizing a weighted sum of the MSE losses associated with each learning objective over a set of collocation points, with effectiveness influenced by many factors---including the intrinsic properties of the PDE itself, network architectures, optimization methods, and training heuristics. A fully exhaustive enumeration of all factors is infeasible, as controlled experiments for causal attribution require careful isolation of variables. Hence, we focus on three key aspects that are fundamental to the neural solver's performance: two representative network architectures (MLP vs. GNN), two constraint formulations (AD vs. discretization-based), and the effects of problem nonlinearity and boundary complexity. Our analysis framework is illustrated in \cref{fig:analysis_framework}, and the details are elaborated in the following subsections.

\subsection{Network Architectures}\label{sec:network_architecture}
MLP, also known as the fully-connected neural network, applies successive linear transformations followed by nonlinear activations to achieve functional approximation, defined by
\begin{equation}\label{eq:mlp_network}
    \mathbf{h}^{(k + 1)} =  \sigma \big(\mathbf{W}^{\top}\mathbf{h}^{(k)} + \mathbf{b} \big),
\end{equation}
\noindent where $\mathbf{h}^{(k+1)}$ is the layer output, $\mathbf{W}$ and $\mathbf{b}$ are the learnable weight matrix and bias, and $\sigma$ is the activation.

The pointwise nature of MLP, where each input is processed independently without explicit neighbourhood coupling, makes it naturally compatible with AD, which satisfies the requirements for continuous function approximation. However, this learning approach can only extract features from the local coordinate information of each point to map the desired solution, lacking the receptive field of architectures such as CNNs or GNNs, which explicitly aggregate information from neighbouring points or graph nodes. Yet, the extent to which such additional neighbourhood information can improve performance, and when it justifies sacrificing the efficiency and convenience of mesh-free pointwise evaluation, remains unclear, as there is currently no established reference or guideline. This is precisely the question we aim to address.

We choose GNN as the representative architecture for incorporating neighbourhood information, rather than CNN. This choice is motivated not only by CNN's inherent limitation to structured grids, but also by the fact that its input representation and parameter structure differ substantially from those of MLP. This makes it difficult to maintain consistent data storage formats and hard constraint enforcement across architectures, which would introduce confounding factors that undermine the controlled-variable design of our comparison. The key difference between GNN and MLP lies in the message-passing aggregation mechanism, given by \cref{eq:gnn_aggr}, which enables the exchange of information between neighbouring nodes. This becomes increasingly relevant as the PDE nonlinearity strengthens, since nonlinear solutions tend to exhibit stronger local interactions that cannot be captured by pointwise mappings alone. However, stacking too many aggregation layers introduces known limitations such as over-smoothing and over-squashing \parencite{2022_arxiv_graphbottlenecks, 2026_tmlr_gnnsurvey}, which discourages arbitrarily deep architectures.
\begin{equation}\label{eq:gnn_aggr}
    \mathbf{h}_{i}^{(k + 1)} =  \phi \big(\mathbf{h}_{i}^{(k)}, \psi \underset{j \in N_{i}}{\bigoplus}\big(\mathbf{h}_{i}^{(k)}, \mathbf{h}_{j}^{(k)}, \mathbf{e}_{ij} \big) \big),
\end{equation}
\noindent where $\bigoplus$ denotes a permutation-invariant aggregation operator (e.g., sum, mean, or attention); $\phi$ and $\psi$ are learnable functions (e.g., MLPs).

\subsection{Continuous and discrete constraints}\label{sec:physical_constraints}
Theoretically, feedforward networks with sufficient width can approximate continuous functions defined over compact domains to arbitrary accuracy \parencite{1989_nn_approximationtors}. This property underpins AD-based PINNs, where derivatives are obtained analytically through chain-rule differentiation of the network ansatz. Nevertheless, the theorem only asserts the existence of suitable network parameters but offers no guarantee that gradient-based training via backpropagation, which does not strictly follow a deterministic optimization path as numerical iteration, will practically converge toward the desired solutions. There accordingly exists a gap between theoretical expressive capacity and real-world performance. Recognizing this practical approximation ceiling and its dominant influencing factors constitutes a delineation of practical operating bounds for engineering-oriented applications.

Discretization-based constraint formulations for neural networks share the same forward pass procedure as their AD-based counterparts. While they introduce considerable truncation errors originating from Taylor-series expansions within the discretized space, such errors are deterministic and stable, and are directly imposed on the predicted solution itself. By contrast, for AD-based constraints, derivatives are computed through the full network propagation chain originating from the predicted solution. This chain grows longer with increasing network depth. Such embedding of physical constraints tightly couples the entire training process, giving rise to accumulated error effects.
Our comparison of the two constraint-informed paradigms is not merely aimed at judging which approach performs better. More importantly, it seeks to answer what governs the practical approximation ceiling and to what extent performance can be expected under general conditions free of ad-hoc empirical tricks. Thus, we resort to a more rigorous error analysis to achieve this goal.

\subsection{Decomposition of errors}\label{sec:error_sources}
In practical applications, most PDE problems do not possess exact analytical solutions. Evaluating the performance of neural-network PDE solvers solely through standalone error metrics is therefore insufficient. On one hand, reference numerical solutions obtained from lower-fidelity solvers with limited formal accuracy can lead to misleading assessments: a model prediction that better captures the true physical behaviour may exhibit larger errors with respect to such imperfect reference solutions. On the other hand, the achievable approximation upper bound varies substantially across problems of different inherent difficulty. Isolated error-metric figures (e.g., a certain reported $L_{2}$ or MSE error value) provide no broadly accepted reference frame to indicate to what extent a model performs. 
Accordingly, for error computation, we adopt high-fidelity DNS solutions with suppressed numerical errors as substitutes for unavailable analytical solutions. Numerical results obtained under identical discretization settings are further included as baseline anchors. Using these DNS-derived anchors, we can characterize how closely a neural-network model can approach the accuracy level of a given conventional numerical method. Nevertheless, the aforementioned two issues cannot be entirely eliminated. Since error comparisons require evaluation at identical spatial locations, projecting the high-fidelity reference solution onto target discrete points inevitably introduces additional interpolation errors. This relative evaluation paradigm still yields objective and interpretable performance descriptions. It further enables us to delineate the applicability of neural-network approaches for practical engineering scenarios. For problems where higher accuracy can be achieved, neural-network solvers may serve as direct alternatives to conventional numerical solvers. For complex problems where only low-accuracy solutions are attainable, these methods are better suited for preliminary screening, and research efforts should instead be directed toward improving computational efficiency.

To obtain the present analysis, we mainly divide errors into three types: approximation error, optimization error, and truncation error. It should be noted that this decomposition does not follow strict mathematical definitions of error sources, since reliable analytical tools for the decomposition of strongly coupled approximation- and optimization-related errors in neural-network-based PDE solvers are still lacking. Although several works \parencite{2021_nips_pinnfailure,2022_jcp_pinnfail} have analyzed this issue from the perspective of loss landscapes and attributed many difficulties to ill-posed optimization, complete disentanglement between these two error components remains difficult given limited test cases and network parameter configurations.
Therefore, tailored to the focus of this work, we give specialized definitions: approximation error specifically denotes the bias induced by network architectures, with a primary focus on MLP and GNN structures; optimization error specifically denotes the bias introduced by physics-constraint embedding strategies as well as soft- and hard-constraint enforcement schemes for boundary conditions, mainly reflecting influences brought by AD versus discretization-based constraint formulations under MLP networks; truncation error serves as an auxiliary anchor term to assist our analysis of neural-network errors under discretized settings. Since discrete truncation errors dominate the overall error budget for most conventional numerical methods \parencite{2002_springer_cfd}, well-converged numerical solutions can therefore be adopted as reference anchors in our analysis.  For those indispensable approximations introduced for numerical stability or completeness, which are problem-dependent, we collectively categorize them as dissipation-related contributions and will analyse them individually for each physical scenario.
Our comparative analysis strictly adheres to controlled-variable principles. Given the computational overhead of double-precision neural network training and the fact that round-off errors are negligible relative to the three preceding error components, all neural network models employ single-precision arithmetic, in contrast to the uniform double-precision setting used for DNS solutions.
\FloatBarrier

\section{Experiments}\label{sec:experiments}
\subsection{Implementation details}\label{sec:implementation_details}
\begin{figure}[!ht]
    \centering
    \begin{minipage}{0.2\textwidth}
        \centering
        \includegraphics[width=1\textwidth]{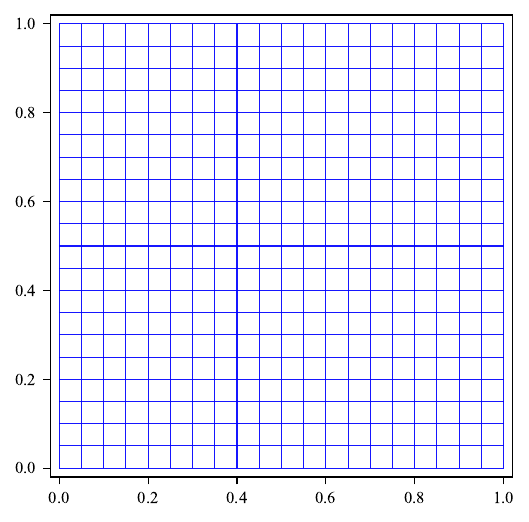}
        \caption*{(a)}
    \end{minipage}
    \begin{minipage}{0.21\textwidth}
        \centering
        \includegraphics[width=1\textwidth]{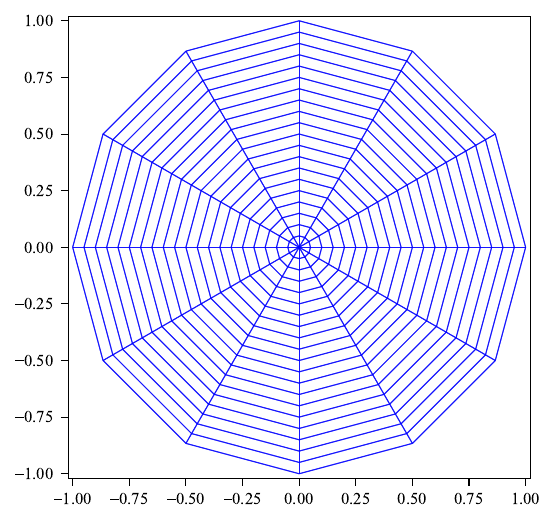}
        \caption*{(b)}
    \end{minipage}
    \begin{minipage}{0.37\textwidth}
        \centering
        \includegraphics[width=1\textwidth]{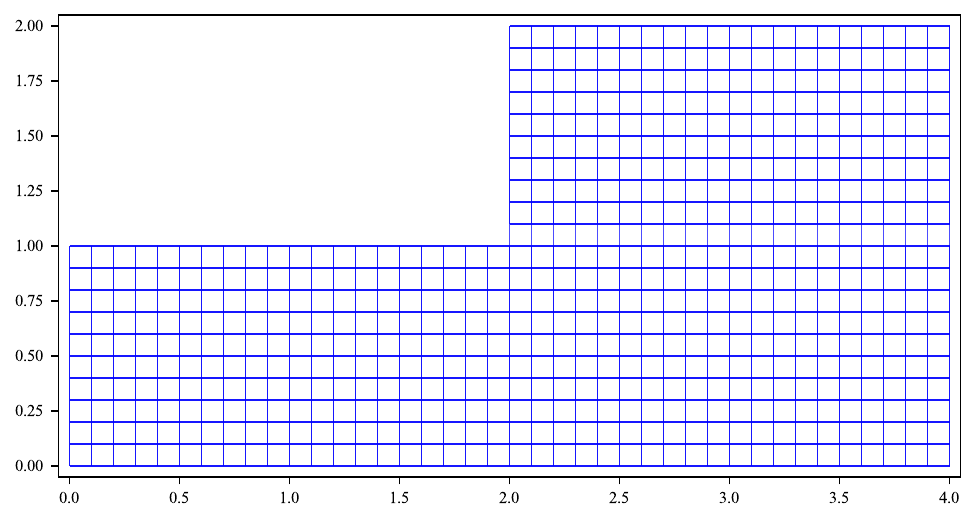}
        \caption*{(c)}
    \end{minipage}
    \begin{minipage}{0.081\textwidth}
        \centering
        \includegraphics[width=1\textwidth]{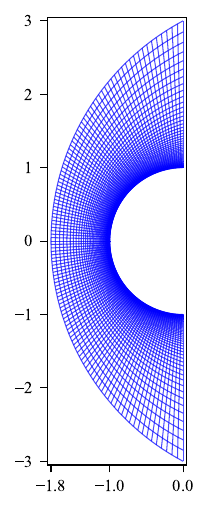}
        \caption*{(d)}
    \end{minipage}
    \begin{minipage}{0.081\textwidth}
        \centering
        \includegraphics[width=1\textwidth]{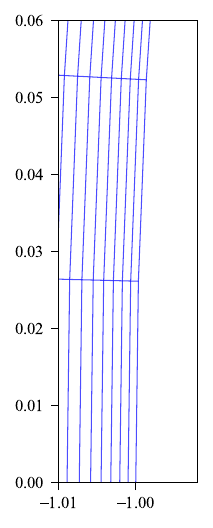}
        \caption*{(e)}
    \end{minipage}
    \caption{Mesh configurations: (a) isotropic uniform Cartesian box mesh; (b) uniform polar mesh; (c) uniform Cartesian mesh for the non-convex backward-facing step domain; (d) anisotropic hypersonic mesh with boundary-layer refinement; (e) zoom-in view of the boundary-layer region shown in (d).}
    \label{fig:mesh}
\end{figure}

To isolate the effects of spatial discretization and nonlinearity, we choose two-dimensional steady-state problems for all test cases, so as to avoid introducing additional confounding factors. The selected cases are designed with progressively increasing nonlinearity and boundary condition complexity, covering a hierarchical spectrum of physical problems: canonical Poisson-type equations with analytical solutions, including linear sine-forced, nonlinear polynomial-forced, and stiff exponential-forced (Liouville) problems, referred to hereafter as Sine, Polynomial, and Liouville cases respectively. We further extend the test suite to coupled incompressible Navier-Stokes (NS) flows including lid-driven cavity and backward-facing step problems, and finally to hypersonic inviscid flow featuring bow shock discontinuities over a blunt cylinder. 
Mesh configurations differ across test cases according to geometrical requirements:
the Sine, Polynomial, and lid-driven cavity cases adopt an isotropic uniform Cartesian box mesh (\Cref{fig:mesh}(a)); the Liouville case uses a uniform polar mesh (\Cref{fig:mesh}(b));
for the backward-facing step problem, a uniform Cartesian mesh is constructed over the non-convex domain (\Cref{fig:mesh}(c)); the hypersonic inviscid flow employs an anisotropic mesh with boundary-layer refinement (\Cref{fig:mesh}(d)).

For DNS solutions, standard FD discretization is adopted for the Poisson-family test cases.
The Polynomial case relies on conventional Newton iteration for nonlinear updates \parencite{2002_springer_cfd}. Owing to the stiffness of the Liouville case, standard Newton iteration readily diverges unless an excessively loose tolerance is enforced. The higher-order Halley iteration is therefore deployed instead \parencite{1997_bams_halley}.
Incompressible NS computations are carried out on staggered grids using a fully-coupled FD formulation \parencite{1986_jcp_blockimplicit}. Contrary to the reference implementation, we employ exact analytical Jacobians for implicit iterations rather than symmetric-coupled Gauss-Seidel updates, aiming to avoid unnecessary additional approximation errors.
Global-scale $L_{2}$-norm errors serve as the evaluation metric for the preceding five cases, where neural network models are capable of achieving relatively high accuracy.
For the challenging hypersonic flow, however, neural network training cannot yield high-fidelity results, notably failing to resolve accurate shock positions.
It is only adequate for fast preliminary engineering estimations, so local quantity metrics are adopted for assessment, including stagnation-point and wall-surface observables. The corresponding DNS reference is computed using vertex-centered FV discretization with dual control volumes \parencite{2016_moukalled_fvm}, where the Rusanov scheme \parencite{1962_cmmp_rusanov} is combined with Venkatakrishnan-type MUSCL reconstruction \parencite{1995_jcp_venkata}. Explicit first-order pseudo-time stepping is adopted to produce a coarse-level reference solution, ensuring a fair comparison with the low-fidelity network predictions. All DNS solvers are fully implemented in Python, with Numpy and SciPy for CPU computation and CuPy \parencite{2017_nips_cupy} for GPU acceleration. This consistent coding framework ensures fair efficiency evaluations under identical hardware conditions. Complete implementation details are documented in Appendix~\ref{appendix_a}.

For neural network-based methods, comprehensive benchmark investigations in PINNacle \parencite{2024_nips_pinnacle} have demonstrated that existing PINN variants exhibit no consistent and significant performance advantages over vanilla PINNs. For detailed benchmark comparisons of various PINN architectures, we refer readers to this published work. Accordingly, we adopt the vanilla PINN with a standard MLP structure as our representative baseline for AD-constrained models. We further construct three network configurations with identical hyperparameter settings, including continuous MLP-based PINN, discrete MLP-based PINN, and GNN-based PINN. These models are denoted as PINN\textsubscript{c}, PINN\textsubscript{d}, and PIGNN, respectively, in the following discussions. 
For Poisson cases, discretization-based PINN models adopt identical discretization formulations as the DNS solver. For incompressible NS flows, collocated grids are used to facilitate boundary constraint imposition. Similarly, vertex-centered FV discretization is selected for compressible hypersonic flows to enable robust hard constraint enforcement. Hard boundary constraints are prioritized in this work, as they improve the overall accuracy by several orders of magnitude under identical hyperparameter settings, which is validated by our soft-constraint comparative experiments in Appendix~\ref{appendix_c}. Full hard constraint implementation is applied whenever feasible. For cases on more complex domains that cannot accommodate global hard constraints, including the backward-facing step and hypersonic flow problems, localized hard constraint strategies are adopted. 

Importantly, this work adheres to a strict general-purpose experimental principle across all validation test cases. No case-specific priors are incorporated into the optimization pipeline. Training relies exclusively on the governing PDEs and physically meaningful boundary conditions, with boundary constraints assigned equal weighting alongside the PDE residual terms. Case-specific hand-tuned configurations can only deliver performance gains for isolated test scenarios and cannot constitute evidence of general algorithmic capability. In contrast, generic, problem-agnostic priors are considered legitimate. For example, near-wall clustering, either for mesh generation or collocation-point sampling, is adopted herein as a standard, broadly applicable practice for wall-bounded flow problems. Nevertheless, manually pre-specified knowledge such as exact shock positions corresponding to a particular geometry is strictly excluded from our entire framework. Detailed constraint enforcement schemes and loss function configurations are documented in Appendix~\ref{appendix_b}.

We have tested several activation functions. Tanh is selected for the linear Sine case, while Swish \parencite{2017_arxiv_swish} is adopted for all other nonlinear cases. Regarding graph aggregators, attention-based operators exhibit good robustness across general scenarios; therefore, the graph attention network (GAT) \parencite{2018_arxiv_gat} with five attention heads is used for all PIGNN models. Unless otherwise specified, all training configurations follow the settings summarized in \cref{tab:training_params}.
All neural network models are implemented in PyTorch \parencite{2019_nips_pytorch}.
We employ Glorot initialization \parencite{2010_aistats_glorotinit} for network weight parameters. The Adam optimizer \parencite{2017_arxiv_adam} is utilized with an initial learning rate of $0.01$. Step-wise learning-rate decay is applied, with a decay factor of $0.9$ every $1000$ epochs, and training terminates at a fixed epoch count. To eliminate sampling-related interference and mitigate randomness, PINN\textsubscript{c} adopts exactly the same discrete points as the discretization-based solver. Each experiment is repeated five times with different random initializations, and we report the averaged results.
It should be noted that we do not perform exhaustive hyperparameter searches to pursue maximum achievable model performance for individual cases. Instead, the adopted hyperparameters serve as general purpose baseline settings oriented toward practical engineering applications under efficiency and robustness considerations. This choice is motivated by the intrinsic lack of separate validation datasets for PINN-based PDE solving. Only basic network dimensions and training epoch counts are briefly explored, since fair comparisons under such general settings are prioritized over case-specific optimal performance. Hyperparameter analysis limited to PINN\textsubscript{c} models is provided in Appendix~\ref{appendix_d}.

\begin{table}[!ht]
    \centering
    \caption{Network architecture and training hyperparameters for each test PDE case.}
    \label{tab:training_params}
    \small
    \begin{tabularx}{0.75\textwidth}{@{\hspace{0.5em}}llllll@{\hspace{0.5em}}}
        \toprule
        PDE case & Epochs & Layers & Hidden units & Activation & BC enforcement \\
        \midrule
        Poisson1 (Sine) & 20000 & 2 & 16 & Tanh & Fully-hard \\
        Poisson2 (Polynomial) & 20000 & 2 & 16 & Swish & Fully-hard \\
        Poisson3 (Liouville) & 20000 & 2 & 16 & Swish & Fully-hard \\
        Lid-driven cavity flow & 20000 & 4 & 32 & Swish & Fully-hard \\
        Backward-facing step flow & 20000 & 4 & 32 & Swish & Partially-hard \\
        Hypersonic inviscid flow & 10000 & 4 & 32 & Swish & Partially-hard \\
        \bottomrule
    \end{tabularx}
\end{table}

\subsection{Linear to nonlinear Poisson}\label{sec:poisson}
We consider a family of elliptic Poisson problems governed by
\begin{equation}\label{eq:poisson_formulation}
    -\nabla^{2}u = f(x, y, u), \quad (x, y) \in \Omega.
\end{equation}

Nonlinearity of these test problems is introduced and gradually increased via the source term $f$.
The linear Sine case, $f(x, y) = 2\pi^{2} \sin (\pi x) \sin (\pi y)$, is defined on $\Omega:(x, y) \in [0, 1] \times [0, 1]$, with Dirichlet boundary condition $u=0$ at $x=0,1$ and $y=0,1$, and analytical solution $u = \sin(\pi x) \sin(\pi y)$;
the nonlinear Polynomial case, $f(u) = -4(u + u^{3})$, is defined on $\Omega:(x, y) \in [-\frac{\pi}{6}, \frac{\pi}{6}] \times [-\frac{\pi}{6}, \frac{\pi}{6}]$, with Dirichlet boundary condition $u=\tan(y \pm \frac{\pi}{6})$ on $x=\pm\frac{\pi}{6}$ and $u=\tan(x \pm \frac{\pi}{6})$ on $y=\pm\frac{\pi}{6}$, and analytical solution $u = \tan(x+y)$;
the Liouville case, $f(u) = 2\mathrm{e}^{u}$, is defined on the unit disk $\Omega: x^{2} + y^{2} \le 1$, with Dirichlet boundary condition $u=0$ on $x^{2}+y^{2}=1$, and analytical solution $u=\ln 4 - 2\ln\bigl(1+x^{2}+y^{2}\bigr)$.

\begin{figure}[!ht]
    \centering
    \includegraphics[width=1\linewidth]{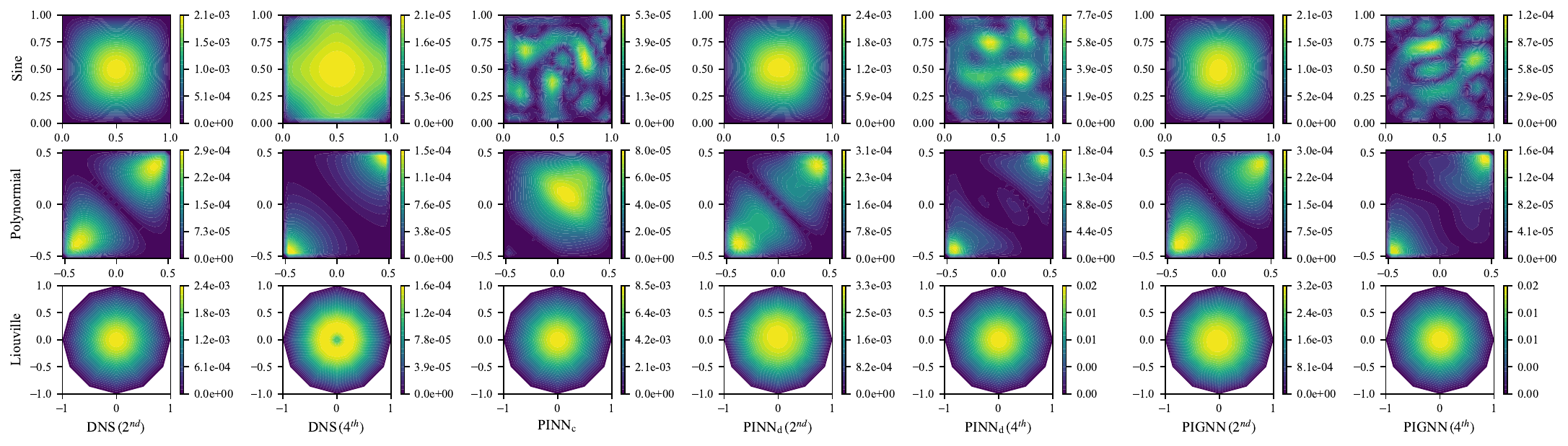}
    \caption{Absolute error maps for the three Poisson-type test cases. Rows from top to bottom correspond to the Sine, Polynomial, and Liouville cases.}
    \label{fig:poisson_error}
\end{figure}

Analytical solutions allow precise error quantification, so we incorporate high-order discretization schemes in our comparisons. To examine the trends of truncation and approximation errors, four mesh resolutions are adopted for each test case. Results for the Sine, Polynomial, and Liouville cases are reported in \cref{tab:sine_l2}, \cref{tab:polynomial_l2}, and \cref{tab:liouville_l2}, respectively.
From the experimental results, a clear trend can be observed: approximation and optimization errors of neural network solvers grow with increasing nonlinearity and mesh refinement. A similar trend holds for high-order discretization schemes, whose benefits gradually diminish and may even become inferior to low-order schemes under stronger nonlinearity and finer resolution. For the linear Sine case at $20\times 20$ resolution, the second-order central difference discretization-constrained neural solver achieves good optimization performance, with errors dominated solely by truncation error, as its error levels closely match those of the second-order DNS solution. When switching to fourth-order discretization, noticeable deviations from the fourth-order DNS solution emerge, which can also be clearly observed in the error maps, as illustrated in \cref{fig:poisson_error}.
Although fourth-order stencils yield larger condition numbers, they still produce positive effects for the lower-resolution Sine and Polynomial cases where optimization ill-conditioning remains mild. For the stiff Liouville case, however, fourth-order discretization exacerbates spectral bias.
From our observations, AD-based PINN\textsubscript{c} exhibit behaviours reminiscent of much higher-order discretization operators. Their advantages are prominent in weakly nonlinear scenarios, where high accuracy can be achieved with relatively few collocation points. This behaviour is also reflected in the statistical stability: the standard deviation and error magnitude of PINN\textsubscript{c} generally lie at the same order of magnitude, whereas discretization-constrained models exhibit more stable training, with standard deviations consistently one order of magnitude lower than their prediction errors at lower resolutions. As nonlinearity intensifies, the inherent instability risk of such high-order-like operators gradually dominates. In contrast, explicit discrete constraints effectively regularize ill-behaved optimization trajectories and yield more robust training and prediction performance under strong nonlinearity. Nevertheless, such benefits are only observed at lower resolutions, and this degradation holds for all considered models. The averaging nature of MSE based-optimization inherently limits performance gains under mesh refinement. We hypothesize that a critical resolution exists, beyond which further mesh refinement brings more drawbacks than benefits, and this critical resolution decreases as nonlinearity becomes stronger.

\begin{table}[!ht]
    \centering
    \caption{Global relative $L_{2}$ errors for the linear Sine test case across different mesh resolutions.}
    \label{tab:sine_l2}
    \small
    \begin{threeparttable}
    \begin{tabularx}{0.85\textwidth}{@{\hspace{0.5em}}lXXXX@{\hspace{0.5em}}}
        \toprule
        \multirow{2}{*}{Method} & \multicolumn{4}{c}{Resolution} \\
        \cmidrule{2-5}
        & $20 \times 20$ & $50 \times 50$ & $100 \times 100$ & $200 \times 200$ \\
        \midrule
        DNS ($2^{nd}$) & $2.06\mathrm{e}{-3}$ & $3.29\mathrm{e}{-4}$ & $8.23\mathrm{e}{-5}$ & $2.06\mathrm{e}{-5}$ \\
        DNS ($4^{th}$) & $2.14\mathrm{e}{-5}$ & $3.24\mathrm{e}{-7}$ & $1.55\mathrm{e}{-8}$ & $8.24\mathrm{e}{-10}$ \\
        PINN\textsubscript{c} & $6.74\mathrm{e}{-5} \pm 6.27\mathrm{e}{-5}$ & $5.74\mathrm{e}{-5} \pm 4.46\mathrm{e}{-5}$ & $5.69\mathrm{e}{-5} \pm 3.45\mathrm{e}{-5}$ & $4.07\mathrm{e}{-5} \pm 1.15\mathrm{e}{-5}$ \\
        PINN\textsubscript{d} ($2^{nd}$) & $2.06\mathrm{e}{-3} \pm 1.13\mathrm{e}{-4}$ & $3.18\mathrm{e}{-4} \pm 9.92\mathrm{e}{-5}$ & $1.17\mathrm{e}{-4} \pm 8.05\mathrm{e}{-5}$ & $7.36\mathrm{e}{-5} \pm 3.49\mathrm{e}{-5}$ \\
        PINN\textsubscript{d} ($4^{th}$) & $5.56\mathrm{e}{-5} \pm 1.60\mathrm{e}{-5}$ & $6.47\mathrm{e}{-5} \pm 2.84\mathrm{e}{-5}$ & $6.76\mathrm{e}{-5} \pm 4.19\mathrm{e}{-5}$ & $2.11\mathrm{e}{-4} \pm 7.67\mathrm{e}{-5}$ \\
        PIGNN ($2^{nd}$) & $2.02\mathrm{e}{-3} \pm 7.38\mathrm{e}{-5}$ & $3.53\mathrm{e}{-4} \pm 6.73\mathrm{e}{-5}$ & $1.33\mathrm{e}{-4} \pm 6.49\mathrm{e}{-5}$ & \underline{$2.85\mathrm{e}{-4} \pm 3.02\mathrm{e}{-5}$} \\
        PIGNN ($4^{th}$) & $5.82\mathrm{e}{-5} \pm 8.16\mathrm{e}{-6}$ & $1.30\mathrm{e}{-4} \pm 3.41\mathrm{e}{-5}$ & $1.29\mathrm{e}{-4} \pm 6.43\mathrm{e}{-5}$ & \underline{$3.26\mathrm{e}{-4} \pm 3.89\mathrm{e}{-5}$} \\
        \bottomrule
    \end{tabularx}
    \begin{tablenotes}
        \footnotesize 
        \item Note: Underlined PIGNN results exhibit optimization instability and adopt early stopping at 10000 epochs.
    \end{tablenotes}
    \end{threeparttable}
\end{table}

\begin{table}[!ht]
    \centering
    \caption{Global relative $L_{2}$ errors for the nonlinear Polynomial test case across different mesh resolutions.}
    \label{tab:polynomial_l2}
    \small
    \begin{tabularx}{0.85\textwidth}{@{\hspace{0.5em}}lXXXX@{\hspace{0.5em}}}
        \toprule
        \multirow{2}{*}{Method} & \multicolumn{4}{c}{Resolution} \\
        \cmidrule{2-5}
        & $20 \times 20$ & $50 \times 50$ & $100 \times 100$ & $200 \times 200$ \\
        \midrule
        DNS ($2^{nd}$) & $1.69\mathrm{e}{-4}$ & $2.79\mathrm{e}{-5}$ & $7.00\mathrm{e}{-6}$ & $1.75\mathrm{e}{-6}$ \\
        DNS ($4^{th}$) & $8.76\mathrm{e}{-5}$ & $5.24\mathrm{e}{-6}$ & $4.74\mathrm{e}{-7}$ & $3.67\mathrm{e}{-8}$ \\
        PINN\textsubscript{c} & $2.33\mathrm{e}{-5} \pm 1.01\mathrm{e}{-5}$ & $1.80\mathrm{e}{-5} \pm 1.21\mathrm{e}{-5}$ & $1.01\mathrm{e}{-5} \pm 4.58\mathrm{e}{-6}$ & $2.88\mathrm{e}{-5} \pm 9.52\mathrm{e}{-6}$ \\
        PINN\textsubscript{d} ($2^{nd}$) & $2.10\mathrm{e}{-4} \pm 6.05\mathrm{e}{-6}$ & $3.59\mathrm{e}{-5} \pm 2.31\mathrm{e}{-6}$ & $2.05\mathrm{e}{-5} \pm 1.17\mathrm{e}{-5}$ & $2.53\mathrm{e}{-5} \pm 1.63\mathrm{e}{-5}$ \\
        PINN\textsubscript{d} ($4^{th}$) & $7.77\mathrm{e}{-5} \pm 3.82\mathrm{e}{-6}$ & $1.58\mathrm{e}{-5} \pm 3.86\mathrm{e}{-6}$ & $4.22\mathrm{e}{-5} \pm 1.21\mathrm{e}{-5}$ & $6.44\mathrm{e}{-5} \pm 6.31\mathrm{e}{-6}$ \\
        PIGNN ($2^{nd}$) & $2.05\mathrm{e}{-4} \pm 6.07\mathrm{e}{-6}$ & $4.15\mathrm{e}{-5} \pm 8.44\mathrm{e}{-6}$ & $1.70\mathrm{e}{-5} \pm 2.90\mathrm{e}{-6}$ & $3.19\mathrm{e}{-5} \pm 3.16\mathrm{e}{-6}$ \\
        PIGNN ($4^{th}$) & $7.24\mathrm{e}{-5} \pm 3.17\mathrm{e}{-6}$ & $1.52\mathrm{e}{-5} \pm 5.87\mathrm{e}{-6}$ & $2.88\mathrm{e}{-5} \pm 1.08\mathrm{e}{-5}$ & $6.92\mathrm{e}{-5} \pm 6.95\mathrm{e}{-6}$ \\
        \bottomrule
    \end{tabularx}
\end{table}

\begin{table}[!ht]
    \centering
    \caption{Global relative $L_{2}$ errors for the stiff Liouville test case across different mesh resolutions.}
    \label{tab:liouville_l2}
    \small
    \begin{threeparttable}
    \begin{tabularx}{0.85\textwidth}{@{\hspace{0.5em}}lXXXX@{\hspace{0.5em}}}
        \toprule
        \multirow{2}{*}{Method} & \multicolumn{4}{c}{Resolution} \\
        \cmidrule{2-5}
        & $N_{r}=20, \, N_{\theta}=12$ & $N_{r}=50, \, N_{\theta}=24$ & $N_{r}=100, \, N_{\theta}=36$ & $N_{r}=200, \, N_{\theta}=48$ \\
        \midrule
        DNS ($2^{nd}$) & $1.75\mathrm{e}{-3}$ & $2.80\mathrm{e}{-4}$ & $7.00\mathrm{e}{-5}$ & $1.75\mathrm{e}{-5}$ \\
        DNS ($4^{th}$) & $1.13\mathrm{e}{-4}$ & $3.00\mathrm{e}{-6}$ & $1.88\mathrm{e}{-7}$ & $1.79\mathrm{e}{-8}$ \\
        PINN\textsubscript{c} & $7.36\mathrm{e}{-3} \pm 1.78\mathrm{e}{-3}$ & $5.65\mathrm{e}{-2} \pm 1.59\mathrm{e}{-2}$ & $4.24\mathrm{e}{-2} \pm 1.83\mathrm{e}{-2}$ & $6.44\mathrm{e}{-2} \pm 8.39\mathrm{e}{-4}$ \\
        PINN\textsubscript{d} ($2^{nd}$) & $2.34\mathrm{e}{-3} \pm 2.52\mathrm{e}{-4}$ & $3.70\mathrm{e}{-2} \pm 1.96\mathrm{e}{-2}$ & $5.70\mathrm{e}{-2} \pm 1.32\mathrm{e}{-2}$ & $5.77\mathrm{e}{-2} \pm 1.51\mathrm{e}{-2}$ \\
        PINN\textsubscript{d} ($4^{th}$) & $1.09\mathrm{e}{-2} \pm 2.13\mathrm{e}{-3}$ & $4.74\mathrm{e}{-2} \pm 1.98\mathrm{e}{-2}$ & $4.04\mathrm{e}{-2} \pm 1.44\mathrm{e}{-2}$ & $5.19\mathrm{e}{-2} \pm 1.64\mathrm{e}{-2}$ \\
        PIGNN ($2^{nd}$) & $5.59\mathrm{e}{-3} \pm 1.41\mathrm{e}{-3}$ & $1.91\mathrm{e}{-2} \pm 2.30\mathrm{e}{-3}$ & \underline{$2.45\mathrm{e}{-2} \pm 2.25\mathrm{e}{-3}$} & \underline{\underline{$2.40\mathrm{e}{-2} \pm 1.62\mathrm{e}{-2}$}} \\
        PIGNN ($4^{th}$) & $1.77\mathrm{e}{-2} \pm 4.25\mathrm{e}{-3}$ & $2.37\mathrm{e}{-2} \pm 3.24\mathrm{e}{-3}$ & \underline{$2.26\mathrm{e}{-2} \pm 3.56\mathrm{e}{-3}$} & \underline{\underline{$2.60\mathrm{e}{-2} \pm 2.46\mathrm{e}{-2}$}} \\
        \bottomrule
    \end{tabularx}
    \begin{tablenotes}
        \footnotesize 
        \item Note: Underlined and double-underlined PIGNN results exhibit optimization instability and adopt early stopping at 10000 and 2000 epochs, respectively.
    \end{tablenotes}
    \end{threeparttable}
\end{table}

\begin{figure}[!ht]
    \centering
    \begin{minipage}{0.33\textwidth}
        \centering
        \includegraphics[width=1\textwidth]{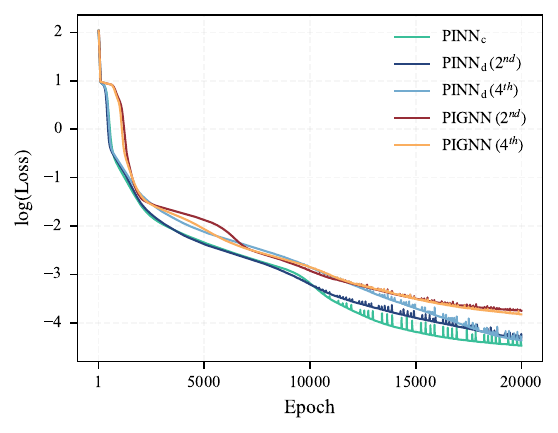}
        \caption*{(a)}
    \end{minipage}
    \begin{minipage}{0.33\textwidth}
        \centering
        \includegraphics[width=1\textwidth]{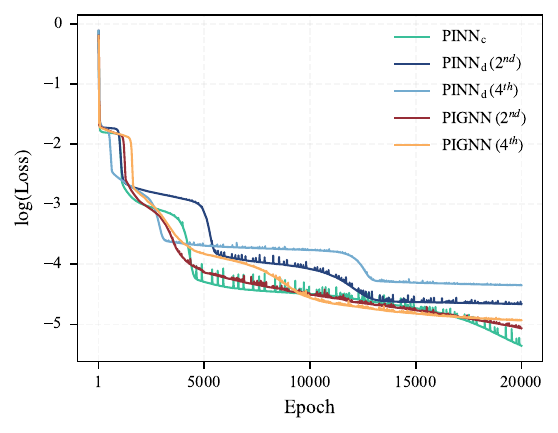}
        \caption*{(b)}
    \end{minipage}
    \begin{minipage}{0.33\textwidth}
        \centering
        \includegraphics[width=1\textwidth]{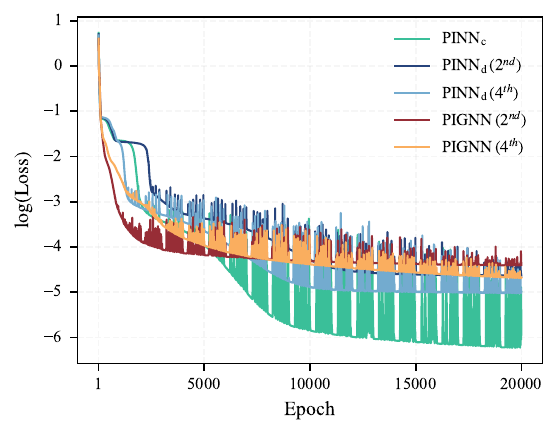}
        \caption*{(c)}
    \end{minipage}
    \begin{minipage}{0.33\textwidth}
        \centering
        \includegraphics[width=1\textwidth]{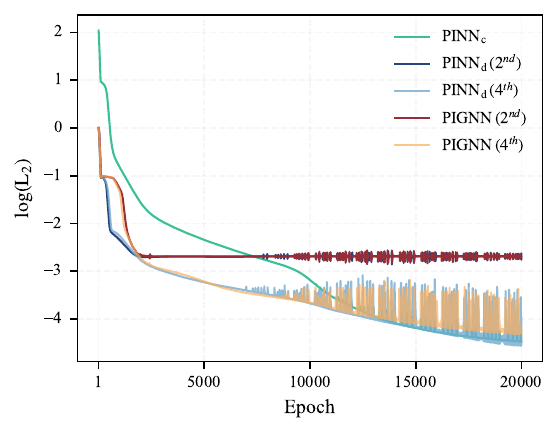}
        \caption*{(d)}
    \end{minipage}
    \begin{minipage}{0.33\textwidth}
        \centering
        \includegraphics[width=1\textwidth]{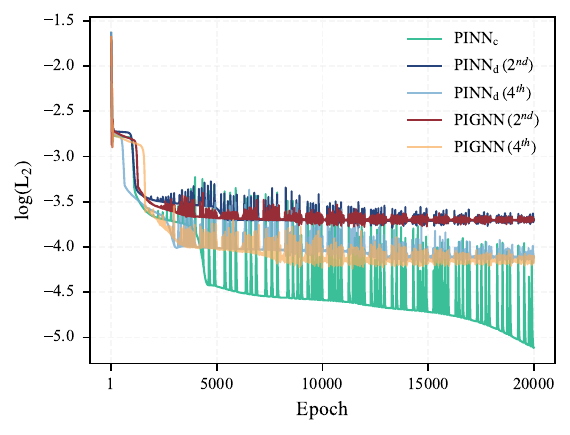}
        \caption*{(e)}
    \end{minipage}
    \begin{minipage}{0.33\textwidth}
        \centering
        \includegraphics[width=1\textwidth]{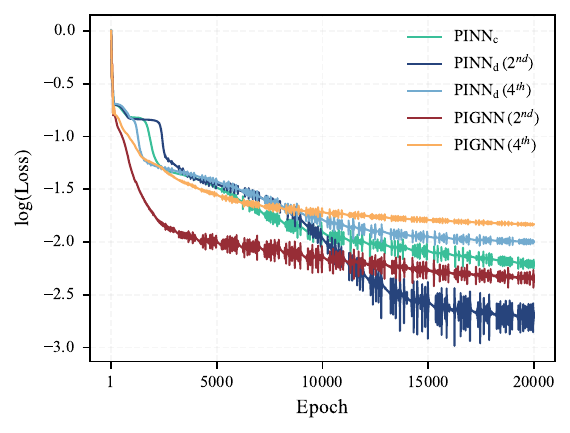}
        \caption*{(f)}
    \end{minipage}
    \caption{Mismatch phenomenon between residual loss and its corresponding $L_{2}$ error for three test cases at the lowest mesh resolution: loss (a) and $L_{2}$ error (d) for the Sine case; loss (b) and $L_{2}$ error (e) for the Polynomial case; loss (c) and $L_{2}$ error (f) for the Liouville case.}
    \label{fig:poisson_convergence}
\end{figure}

Loss convergence is adopted as a performance metric in many existing studies. However, since neural network optimization is inherently non-convex, the averaged loss only reflects a global convergence level and cannot serve as a reliable evaluation metric. Particularly under strong nonlinearity, lower loss values may even amplify errors in high-frequency regions. As illustrated in \cref{fig:poisson_convergence}, mismatches between converged loss and $L_{2}$ error are frequently observed. In the Polynomial case, the second-order scheme with larger error achieves better loss convergence than the fourth-order counterpart. Conversely, in the Liouville case, the fourth-order scheme, which yields larger errors, exhibits superior loss convergence compared with the second-order one. As can be seen from the error maps, the improved loss convergence of the higher-order scheme originates from smoothed regions, while larger errors are shifted toward the central high-frequency region. Different network architectures and constraint-embedding strategies lead to distinct optimization trajectories, and loss values obtained along these divergent trajectories are not directly comparable. 

In terms of accuracy, under full hard-constraint enforcement, the PDE fitting capabilities of MLPs and GNNs are nearly comparable for these simple single-equation test cases. Neither shows a statistically significant advantage, with errors fluctuating within the same order of magnitude; either approach may yield better results given minor parameter perturbations. Nevertheless, MLPs outperform GNNs in terms of computational efficiency and training stability, particularly at high resolutions and for purely soft-constrained models (see Appendix~\ref{appendix_c}). GNNs are more prone to stagnant or even rising loss and training divergence. Such behaviour can be mitigated by reducing the initial learning rate or increasing the learning-rate decay factor. To maintain identical hyper-parameter settings across all models in the present work, early stopping is adopted instead. We conjecture that this phenomenon likely arises from the low problem complexity, where additional neighbourhood information cannot be fully exploited and introduces excessive overhead for GNNs. From an engineering perspective, however, the tendency of GNNs toward loss oscillation or divergence under poor learning performance is not entirely disadvantageous. It serves as a useful diagnostic signal: the loss better correlates with actual model performance. Poor model behaviour is directly reflected in the loss trajectory, which flags potential issues for further adjustment. Consequently, GNN loss can be more trustworthy than MLP loss when no post-hoc reference solution is available.

\subsection{Lid-driven cavity Flow}\label{sec:lid_driven}
We consider the steady-state incompressible NS equations formulated in non-dimensional form:
\begin{equation}\label{eq:ns_formulation}
    \begin{array}{ll}
    \mathbf{u} \cdot \nabla \mathbf{u} + \nabla p - \frac{1}{\mathrm{Re}} \nabla^{2} \mathbf{u} = \mathbf{0} \\[2pt]
    \nabla \cdot \mathbf{u} = 0
    \end{array}, \quad (x, y) \in \Omega,
\end{equation}
\noindent where $\mathbf{u} = [u, v]^{\top}$ denotes the non-dimensional velocity vector,
$p$ is the non-dimensional pressure, and $\mathrm{Re}$ is the Reynolds number.

The lid-driven cavity flow is governed by \cref{eq:ns_formulation}, defined on the domain $\Omega:(x, y) \in [0, 1] \times [0, 1]$, with $\mathrm{Re}=100$.
The boundary conditions are given as
\begin{equation}\label{eq:lid_driven_bc}
    \begin{array}{l}
    u=4x(1-x), \, (x,y)\in \Omega_{\text{top}}; \quad u=0, \, (x,y) \in \partial\Omega \setminus \Omega_{\text{top}} \\[2pt]
    v=0, \quad (x,y) \in \partial\Omega \\[2pt]
    p=0, \quad (x,y)=(0,0)
    \end{array}.
\end{equation}

Multiple mesh resolutions are also tested in this case. Since no analytical solution exists, a high-fidelity DNS solution at $1000 \times 1000$ resolution is taken as the reference for error evaluation. The quantitative results are summarized in \cref{tab:lid_driven_l2}. The optimal performance lies between $100\times 100$ and $150\times 150$. A prominent error reduction is observed when refining from $50\times 50$ to $100\times 100$. Upon further mesh refinement, velocity errors exhibit stagnation, with no meaningful reduction and occasional error increments. Although pressure error still exhibits a decreasing trend under mesh refinement, the corresponding marginal benefit becomes very limited. 
PINN\textsubscript{c} still exhibits noticeable advantages for these coupled systems under weak-nonlinear conditions. However, somewhat unexpectedly, discretization-constrained PINN\textsubscript{d} and PIGNN outperform the numerical DNS solutions at identical lower resolutions. This observation suggests that the overall error of neural network solvers is dominated by truncation errors under such settings. We attribute this accuracy improvement to two primary factors. First lies in the difference of discretization stencils (see Appendix~\ref{appendix_a2} and Appendix~\ref{appendix_b2}). The DNS solver adopts a staggered grid, where velocity components $u$ and $v$ are stored at staggered locations. This arrangement introduces additional approximation errors via boundary interpolation, while neural network solvers enforce full hard boundary conditions with boundary values imposed directly at their physical locations. Second, differences in boundary condition treatment contribute to the discrepancy. To obtain a well-posed system, the DNS solver introduces approximate Neumann-like conditions $\frac{\partial p}{\partial \mathbf{n}} = 0$, which are not physically exact. Further numerical tests indicate that the error contribution from these approximate Neumann-type boundary treatments is considerably smaller than the interpolation errors originating from the staggered-grid layout, as imposing this additional constraint yields only minor changes in the observed error magnitude.

\begin{table}[!ht]
    \centering
    \caption{Global relative $L_{2}$ errors for the lid-driven cavity flow at $\mathrm{Re} = 100$ across different mesh resolutions.}
    \label{tab:lid_driven_l2}
    \small
    \begin{tabularx}{0.85\textwidth}{@{\hspace{0.5em}}llXXXX@{\hspace{0.5em}}}
        \toprule
        \multirow{2}{*}{Method} & \multirow{2}{*}{Metric} & \multicolumn{4}{c}{Resolution} \\
        \cmidrule{3-6}
        & & $50 \times 50$ & $100 \times 100$ & $150 \times 150$ & $200 \times 200$ \\
        \midrule
        \multirow{3}{*}{DNS}
        & $u$ & $3.84\mathrm{e}{-2}$ & $1.53\mathrm{e}{-2}$ & $8.99\mathrm{e}{-3}$ & $6.15\mathrm{e}{-3}$ \\
        & $v$ & $6.79\mathrm{e}{-2}$ & $2.26\mathrm{e}{-2}$ & $1.23\mathrm{e}{-2}$ & $8.11\mathrm{e}{-3}$ \\
        & $p$ & $0.1759$ & $7.32\mathrm{e}{-2}$ & $4.28\mathrm{e}{-2}$ & $2.87\mathrm{e}{-2}$ \\
        \midrule
        \multirow{3}{*}{PINN\textsubscript{c}}
        & $u$ & $9.70\mathrm{e}{-3} \pm 2.95\mathrm{e}{-3}$ & $6.19\mathrm{e}{-3} \pm 9.44\mathrm{e}{-4}$ & $9.79\mathrm{e}{-3} \pm 1.26\mathrm{e}{-3}$ & $6.99\mathrm{e}{-3} \pm 1.54\mathrm{e}{-3}$ \\
        & $v$ & $1.49\mathrm{e}{-2} \pm 4.16\mathrm{e}{-3}$ & $8.79\mathrm{e}{-3} \pm 1.45\mathrm{e}{-3}$ & $1.54\mathrm{e}{-2} \pm 1.93\mathrm{e}{-3}$ & $1.07\mathrm{e}{-2} \pm 2.99\mathrm{e}{-3}$ \\
        & $p$ & $0.1221 \pm 4.28\mathrm{e}{-3}$ & $6.24\mathrm{e}{-2} \pm 1.56\mathrm{e}{-3}$ & $5.04\mathrm{e}{-2} \pm 2.77\mathrm{e}{-3}$ & $4.30\mathrm{e}{-2} \pm 3.62\mathrm{e}{-3}$ \\
        \midrule
        \multirow{3}{*}{PINN\textsubscript{d}}
        & $u$ & $1.22\mathrm{e}{-2} \pm 7.60\mathrm{e}{-4}$ & $8.32\mathrm{e}{-3} \pm 1.15\mathrm{e}{-3}$ & $7.87\mathrm{e}{-3} \pm 2.64\mathrm{e}{-3}$ & $8.62\mathrm{e}{-3} \pm 1.90\mathrm{e}{-3}$ \\
        & $v$ & $2.59\mathrm{e}{-2} \pm 1.26\mathrm{e}{-3}$ & $1.37\mathrm{e}{-2} \pm 3.00\mathrm{e}{-3}$ & $1.19\mathrm{e}{-2} \pm 3.34\mathrm{e}{-3}$ & $1.36\mathrm{e}{-2} \pm 2.66\mathrm{e}{-3}$ \\
        & $p$ & $0.1205 \pm 5.47\mathrm{e}{-3}$ & $6.37\mathrm{e}{-2} \pm 4.56\mathrm{e}{-3}$ & $4.74\mathrm{e}{-2} \pm 3.68\mathrm{e}{-3}$ & $4.31\mathrm{e}{-2} \pm 3.56\mathrm{e}{-3}$ \\
        \midrule
        \multirow{3}{*}{PIGNN}
        & $u$ & $1.14\mathrm{e}{-2} \pm 8.17\mathrm{e}{-4}$ & $7.19\mathrm{e}{-3} \pm 1.38\mathrm{e}{-3}$ & $7.43\mathrm{e}{-3} \pm 2.11\mathrm{e}{-3}$ & $7.21\mathrm{e}{-3} \pm 1.55\mathrm{e}{-3}$ \\
        & $v$ & $2.44\mathrm{e}{-2} \pm 1.01\mathrm{e}{-3}$ & $1.36\mathrm{e}{-2} \pm 3.18\mathrm{e}{-3}$ & $1.21\mathrm{e}{-2} \pm 4.04\mathrm{e}{-3}$ & $1.09\mathrm{e}{-2} \pm 2.36\mathrm{e}{-3}$ \\
        & $p$ & $0.1287 \pm 1.11\mathrm{e}{-2}$ & $6.14\mathrm{e}{-2} \pm 2.20\mathrm{e}{-3}$ & $4.72\mathrm{e}{-2} \pm 3.89\mathrm{e}{-3}$ & $3.67\mathrm{e}{-2} \pm 1.44\mathrm{e}{-3}$ \\
        \bottomrule
    \end{tabularx}
\end{table}

\begin{figure}[!ht]
    \centering
    \includegraphics[width=1\linewidth]{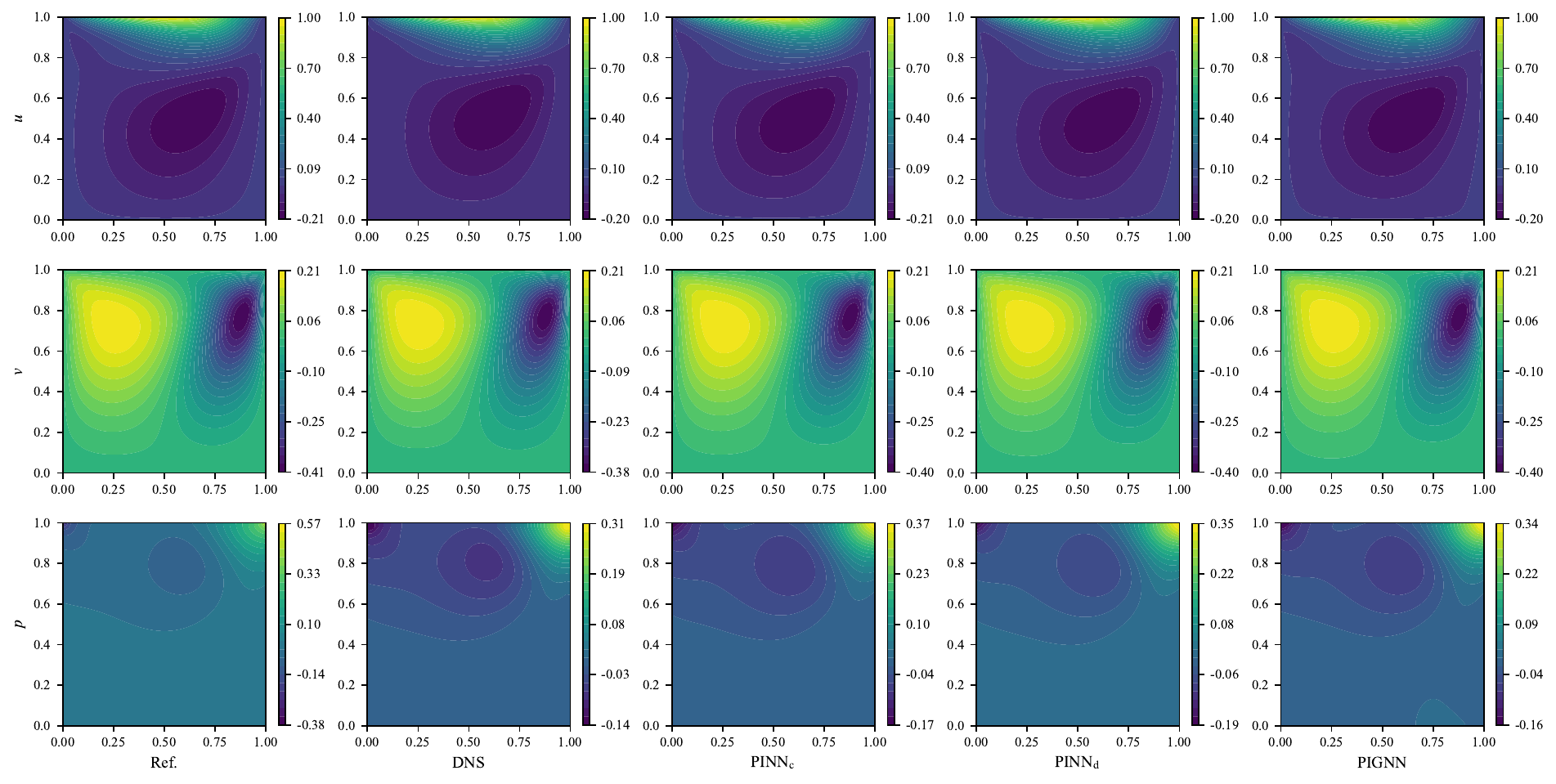}
    \caption{Flow-field comparisons of different methods for lid-driven cavity flow at $\mathrm{Re}=100$. Columns from left to right: reference DNS solution ($1000 \times 1000$), followed by DNS, PINN\textsubscript{c}, PINN\textsubscript{d}, PIGNN at $50 \times 50$ resolution.}
    \label{fig:lid_driven_comparison}
\end{figure}

Qualitative comparison at $50 \times 50$ resolution is illustrated in \cref{fig:lid_driven_comparison}. Neural network solvers faithfully capture fine flow-field textures. From the color-bar value range, it can be observed that neural network solutions are closer to the reference than the DNS result at the same resolution. This validates the necessity of adopting a high-fidelity reference solution, and evaluating neural network models solely against the DNS solution of identical resolution would likely yield misleading conclusions.

\begin{figure}[!ht]
    \centering
    \begin{minipage}{0.33\textwidth}
        \centering
        \includegraphics[width=1\textwidth]{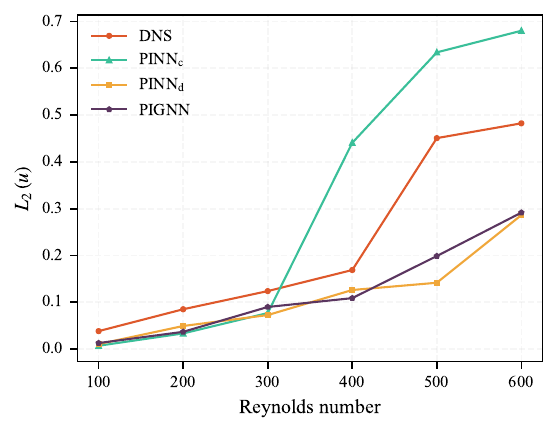}
        \caption*{(a)}
    \end{minipage}
    \begin{minipage}{0.33\textwidth}
        \centering
        \includegraphics[width=1\textwidth]{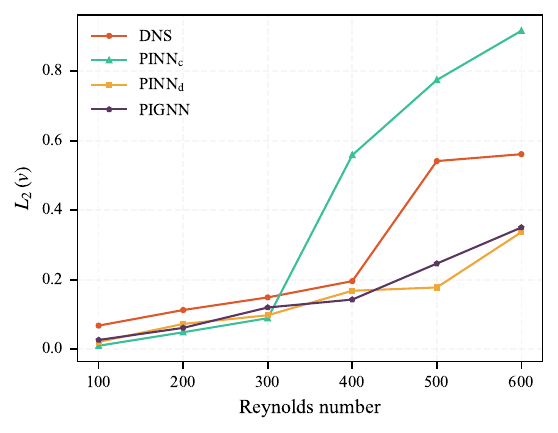}
        \caption*{(b)}
    \end{minipage}
    \begin{minipage}{0.33\textwidth}
        \centering
        \includegraphics[width=1\textwidth]{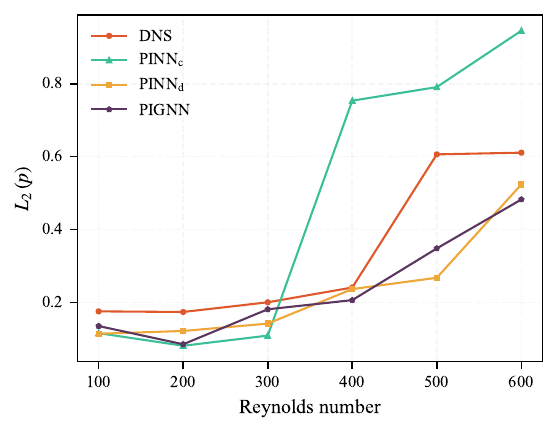}
        \caption*{(c)}
    \end{minipage}
    \caption{Global $L_2$ error trends for lid-driven cavity flow with increasing Reynolds number ($50 \times 50$) against reference DNS solution ($500 \times 500$); neural network models trained for $20000$ epochs except $\mathrm{Re}=500, \, 600$ with $30000$ epochs: (a) horizontal velocity $u$; (b) vertical velocity $v$; (c) pressure $p$.}
    \label{fig:reynolds_curve}
\end{figure}

To further investigate whether the trends observed in the Poisson-type cases hold for these coupled-equation problems under varying nonlinearity, we increase the Reynolds number to strengthen the nonlinear level. A $500 \times 500$ DNS solution is used as the reference, and tests are performed at $50 \times 50$ resolution. 
As illustrated in \cref{fig:reynolds_curve}, PINN\textsubscript{c} maintains a modest advantage over discretization-constrained models for Reynolds numbers up to 300. When Re rises to 400, the error of PINN\textsubscript{c} increases sharply and degrades completely at higher Reynolds numbers, where the predicted outputs lose essential flow-field features. Meanwhile, DNS also exhibits an abrupt deviation from its formerly linear error-growth trend beyond $\mathrm{Re}=400$, which arises from amplified boundary-interpolation errors under intensified flow gradients. 
By contrast, PINN\textsubscript{d} and PIGNN with discrete PDE constraints based on collocated grids exhibit steady and near-linear error growth, though higher Reynolds numbers require longer training convergence. Specifically, 30000 epochs are adopted for the cases at Re = 500 and Re = 600. Additionally, for the failing PINN\textsubscript{c} under higher Reynolds numbers, loss converges normally, offering no visible precursor of solution breakdown. This kind of hidden failure has been previously observed for AD-based training \parencite{2021_nips_pinnfailure}. Such hidden failure modes pose substantial risks for engineering applications, which aligns with our earlier argument that explicit instability signatures are not necessarily detrimental. This observation indicates that neural network training remains robust when discrete PDE constraints are enforced on collocated grids together with global hard constraints. This setup mitigates convex optimization difficulties and achieves lower optimization errors to compensate for deficiencies in truncation errors, avoiding the cumulative degradation observed in AD-based formulations, consistent with the trends seen in the Poisson test cases.

\subsection{Backward-facing step Flow}\label{sec:bward_step}
The backward-facing step flow shares the same governing equations as the lid-driven cavity flow in \cref{eq:ns_formulation}, with $\mathrm{Re}=100$. It is defined on a non-convex computational domain illustrated in \cref{fig:mesh}(c), where $x=0$ denotes the inlet boundary, $x=4$ denotes the outlet boundary, and $\Omega_{\text{wall}}: \big\{(x,y)\,\big|\, y=0 \cup y=2 \cup \big(y=1,\,x\in[0,2]\big) \cup \big(x=2,\,y\in[1,2]\big)\big\}$.
The boundary conditions are given as
\begin{equation}\label{eq:backward_step_bc}
    \begin{array}{ll}
    u=4y(1-y),\; (x,y) \in \Omega_{\text{inlet}}; \quad u=0,\; (x,y) \in \Omega_{\text{wall}} \\[2pt]
    v=0, \quad (x,y)\in \partial \Omega \setminus \Omega_{\text{out}} \\[2pt]
    p=0, \quad (x,y) \in \Omega_{\text{outlet}}
    \end{array}.
\end{equation}

\begin{figure}[!ht]
    \centering
    \includegraphics[width=1\linewidth]{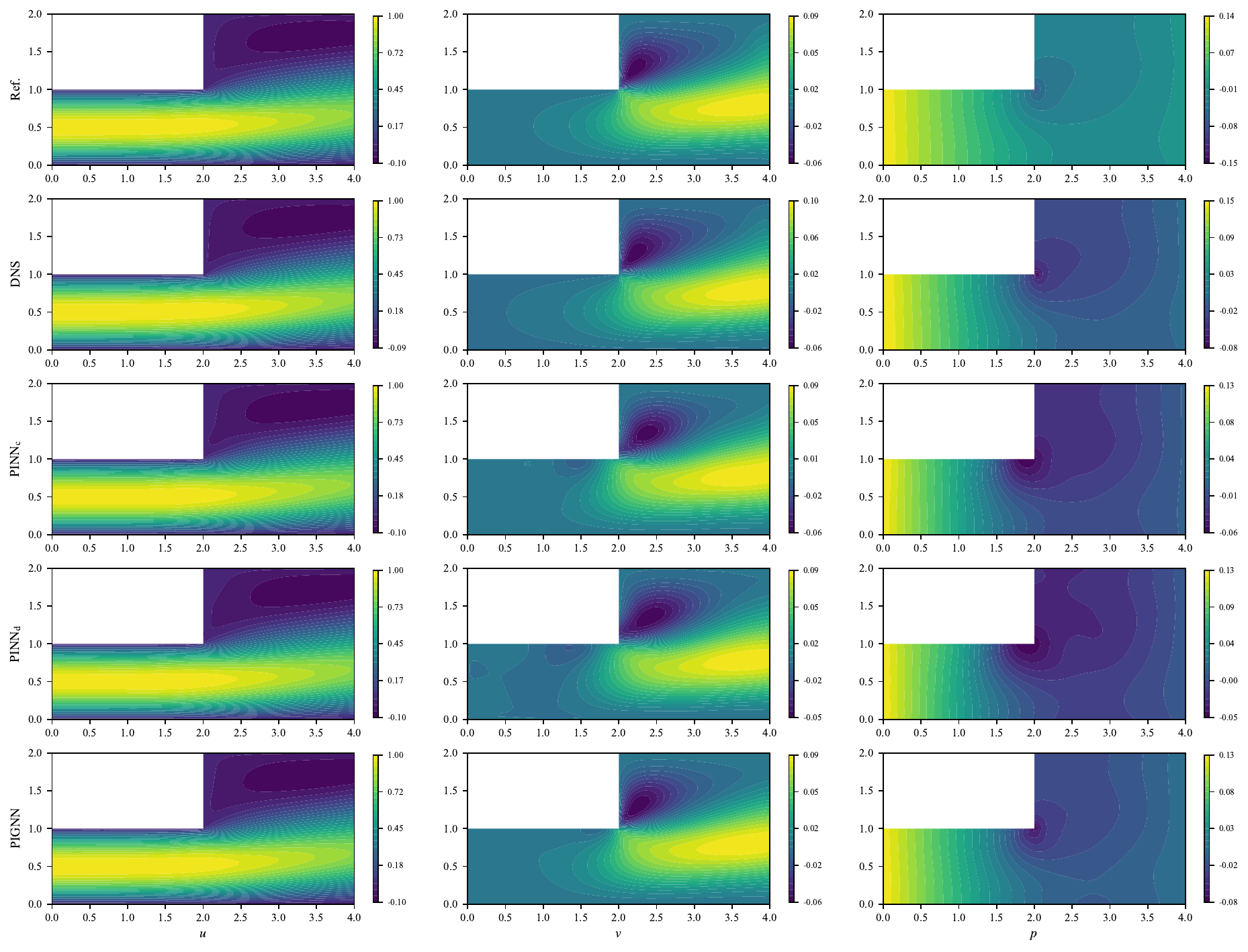}
    \caption{Flow-field comparisons of different methods for backward-facing step flow at  $\mathrm{Re}=100$. Rows from top to bottom: reference DNS solution ($1000 \times 500$), and DNS, PINN\textsubscript{c}, PINN\textsubscript{d}, PIGNN at $100 \times 50$ resolution.}
    \label{fig:bward_step_comparison}
\end{figure}

Unlike the lid-driven cavity, where global hard boundary constraints can be readily constructed, this test case features re-entrant concave corners together with inlet and outlet boundaries. Consequently, the previous global hard-constraint strategy cannot be directly applied here. Only partial hard boundary enforcement is achievable. Specifically, for the horizontal velocity component $u$, hard constraints are enforced on the top and bottom solid walls, while soft constraints are adopted for the inlet and step-side walls. The vertical velocity $v$ is imposed via soft constraints on the step-side walls (for details see Appendix~\ref{appendix_b3}).
As for the discretization-constrained neural models, one may adopt the extrapolation-type outlet boundary conditions used in DNS, namely $\frac{\partial u}{\partial \mathbf{n}} = 0$ and $\frac{\partial v}{\partial \mathbf{n}} = 0$, or skip these outlet conditions entirely. Omitting these conditions yields slightly larger pressure error yet smaller velocity error, whereas activating them reduces pressure error at the cost of slightly elevated velocity error (see \cref{tab:bward_step_outlet_neumann} for quantitative comparison). For this incompressible flow setup, velocity field reconstruction constitutes our primary objective, as pressure is a derived quantity determined by the velocity field. We therefore choose not to apply these outlet extrapolation conditions. The corresponding experimental comparisons are presented in \cref{tab:bward_step_l2}. 
In this case, under identical parameter settings and in the absence of global hard constraints, the performance advantages of MLP-based PINN\textsubscript{c} and PINN\textsubscript{d} relative to DNS vanish, while for the $100 \times 50$ mesh, PIGNN yields better velocity predictions than the DNS at the identical mesh resolution. This can be attributed to differences in discretization schemes and boundary treatments, which we have discussed in the lid-driven cavity section. From the comparison plots in \cref{fig:bward_step_comparison}, we can clearly observe that PIGNN shows good agreement with DNS in capturing intricate local flow-field details around the step corner. By contrast, PINN\textsubscript{c} and PINN\textsubscript{d} suffer from insufficient convergence under the present training budget. This may stem from the singular nature at the step corner, which hinders MLP-based models during multi-objective optimization. Further supplementary experiments indicate that PINN\textsubscript{c} and PINN\textsubscript{d} can achieve improved prediction quality when training is extended to nearly $50000$ epochs. However, such improvements are accompanied by considerable stochasticity, particularly pronounced for the AD-based PINN\textsubscript{c}, and the overall robustness remains inferior to that of PIGNN. In comparison, the GNN leverages additional edge topological information to extract richer features, which helps enhance the network’s approximation capacity and enables more reliable and robust optimization convergence.

\begin{table}[!ht]
    \centering
    \caption{Global relative $L_{2}$ errors for the backward-facing step flow at $\mathrm{Re} = 100$ across different mesh resolutions.}
    \label{tab:bward_step_l2}
    \small
    \begin{threeparttable}
    \begin{tabularx}{0.725\textwidth}{@{\hspace{0.5em}}llXXX@{\hspace{0.5em}}}
        \toprule
        \multirow{2}{*}{Method} & \multirow{2}{*}{Metric} & \multicolumn{3}{c}{Resolution} \\
        \cmidrule{3-5}
        & & $100 \times 50$ & $200 \times 100$ & $300 \times 150$ \\
        \midrule
        \multirow{3}{*}{DNS}
        & $u$ & $1.06\mathrm{e}{-2}$ & $6.64\mathrm{e}{-3}$ & $6.00\mathrm{e}{-3}$ \\
        & $v$ & $9.08\mathrm{e}{-2}$ & $5.30\mathrm{e}{-2}$ & $4.62\mathrm{e}{-2}$ \\
        & $p$ & $0.1077$ & $4.51\mathrm{e}{-2}$ & $5.60\mathrm{e}{-2}$ \\
        \midrule
        \multirow{3}{*}{PINN\textsubscript{c}}
        & $u$ & $1.59\mathrm{e}{-2} \pm 3.71\mathrm{e}{-3}$ & $2.04\mathrm{e}{-2} \pm 5.47\mathrm{e}{-3}$ & $1.53\mathrm{e}{-2} \pm 5.76\mathrm{e}{-3}$ \\
        & $v$ & $0.1315 \pm 9.16\mathrm{e}{-3}$ & $0.1496 \pm 1.93\mathrm{e}{-2}$ & $0.1499 \pm 2.60\mathrm{e}{-2}$ \\
        & $p$ & $0.1217 \pm 1.09\mathrm{e}{-2}$ & $0.1461 \pm 9.07\mathrm{e}{-3}$ & $0.1497 \pm 1.62\mathrm{e}{-2}$ \\
        \midrule
        \multirow{3}{*}{PINN\textsubscript{d}}
        & $u$ & $1.82\mathrm{e}{-2} \pm 6.19\mathrm{e}{-3}$ & $1.87\mathrm{e}{-2} \pm 5.18\mathrm{e}{-3}$ & $1.69\mathrm{e}{-2} \pm 5.89\mathrm{e}{-3}$ \\
        & $v$ & $0.1473 \pm 4.60\mathrm{e}{-2}$ & $0.1497 \pm 3.26\mathrm{e}{-2}$ & $0.1445 \pm 3.33\mathrm{e}{-2}$ \\
        & $p$ & $0.1245 \pm 2.26\mathrm{e}{-2}$ & $0.1545 \pm 9.05\mathrm{e}{-3}$ & $0.1518 \pm 1.30\mathrm{e}{-2}$ \\
        \midrule
        \multirow{3}{*}{PIGNN}
        & $u$ & $6.13\mathrm{e}{-3} \pm 1.28\mathrm{e}{-3}$ & $5.19\mathrm{e}{-3} \pm 1.31\mathrm{e}{-3}$ & $7.74\mathrm{e}{-3} \pm 1.86\mathrm{e}{-3}$ \\
        & $v$ & $8.47\mathrm{e}{-2} \pm 2.32\mathrm{e}{-2}$ & $7.29\mathrm{e}{-2} \pm 5.67\mathrm{e}{-3}$ & $0.1011 \pm 1.66\mathrm{e}{-2}$ \\
        & $p$ & $0.1159 \pm 1.10\mathrm{e}{-2}$ & $0.1427 \pm 1.67\mathrm{e}{-2}$ & $0.1221 \pm 8.29\mathrm{e}{-3}$ \\
        \bottomrule
    \end{tabularx}
    \begin{tablenotes}
    \footnotesize 
    \item Note: (1) For the resolution $w \times h$, $w$ denotes the streamwise grid number of the bottom wall, and $h$ denotes the vertical grid number at the outlet. (2) Cross-resolution error trends are affected by localized interpolation artifacts near the corner singularity, where elevated local errors can occasionally obscure global convergence trends despite qualitative fidelity improvements at finer resolutions.
    \end{tablenotes}
    \end{threeparttable}
\end{table}

\subsection{Hypersonic inviscid cylinder flow}\label{sec:bowshock}
The hypersonic inviscid compressible flow is governed by the Euler equations
\begin{equation}\label{eq:euler_formulation}
    \begin{array}{lll}
    \nabla \cdot (\rho \mathbf{u}) = 0  \\[2pt]
    \nabla \cdot (\rho \mathbf{uu}) + \nabla p = \mathbf{0} \\[2pt]
    \nabla \cdot (\rho E \mathbf{u}) + \nabla \cdot (p \mathbf{u}) = 0
    \end{array}, \quad (x, y) \in \Omega,
\end{equation}
\noindent where $\rho$ is the density ($\mathrm{kg/m^3}$), $\mathbf{u}=[u,v]^\top$ the velocity vector ($\mathrm{m/s}$), $p$ the pressure ($\mathrm{Pa}$), and $E$ the total specific energy, which is given by
\begin{equation}\label{eq:energy_formulation}
    E = \frac{p}{(\gamma - 1) \rho} + \frac{1}{2} \mathbf{u} \cdot \mathbf{u}.
\end{equation}

In our test case, the calorically perfect gas model is adopted, with a constant specific heat ratio $\gamma=1.4$. The equation of state reads $p = \rho R T$, where $R$ denotes the specific gas constant set to $287 \, \mathrm{J}/(\mathrm{kg}\cdot \mathrm{K})$.
Due to the large disparity in variable magnitudes, free-stream sound speed normalization is adopted (see Appendix~\ref{appendix_a3}). Under this normalization, the dimensionless free-stream density is unity, the velocity magnitude equals the Mach number $\mathrm{Ma}$, the free-stream temperature is $1$, and the free-stream pressure reads $\frac{1}{\gamma}$. All flow-field variables presented hereinafter are dimensionless quantities. The corresponding equation of state becomes $p=\frac{\rho T}{\gamma}$, and the sound speed is given by $c = \sqrt{T}$.
Regarding the mesh, progressive refinement is applied near the solid wall (\cref{fig:mesh}(d)). Starting from a near-wall thickness of $10^{-3}$, the mesh is stretched exponentially with a growth factor of $1.08$ (\cref{fig:mesh}(e)). The circumferential direction around the cylinder is discretized into $120$ cells, and the wall-normal direction contains $54$ cells.
The far-field inlet boundary is defined on the outer circular boundary $(x-1.8)^{2} + y^{2} = {3.4}^{2}$, the solid wall corresponds to the unit circle $x^{2}+y^{2}=1.0$, and the outflow boundaries are located on the two lateral sides at $x=0$.

We impose a free-stream boundary condition with $\mathrm{Ma}=8$ at the far-field inlet, where all free-stream variables are strictly enforced via hard constraints. For the Euler equations, the solid-wall slip boundary condition reads $\mathbf{u}\cdot\mathbf{n}=0$. Complete numerical boundary treatment requires certain extrapolation approximations, yet constraint implementations differ across neural network solvers. For the AD-based PINN\textsubscript{c} defined in the continuous space, PDE residuals are directly enforced at outflow boundaries. For the discretization-constrained models, the network can fit boundary relations via discrete coupling terms. Nevertheless, we observe that applying extrapolation-based boundary constraints accelerates convergence considerably for this strongly nonlinear case. Unlike the low Reynolds number incompressible lid-driven cavity and backward-facing-step flows with moderate nonlinearity, where such additional extrapolation approximations introduce slight errors and yield minor side-effects, these treatments produce beneficial effects in the present hypersonic case. For further details on constraint enforcement for different approaches, see Appendix~\ref{appendix_b4}.

\begin{figure}[!ht]
    \centering
    \begin{minipage}{0.47\textwidth}
        \centering
        \includegraphics[width=1\textwidth]{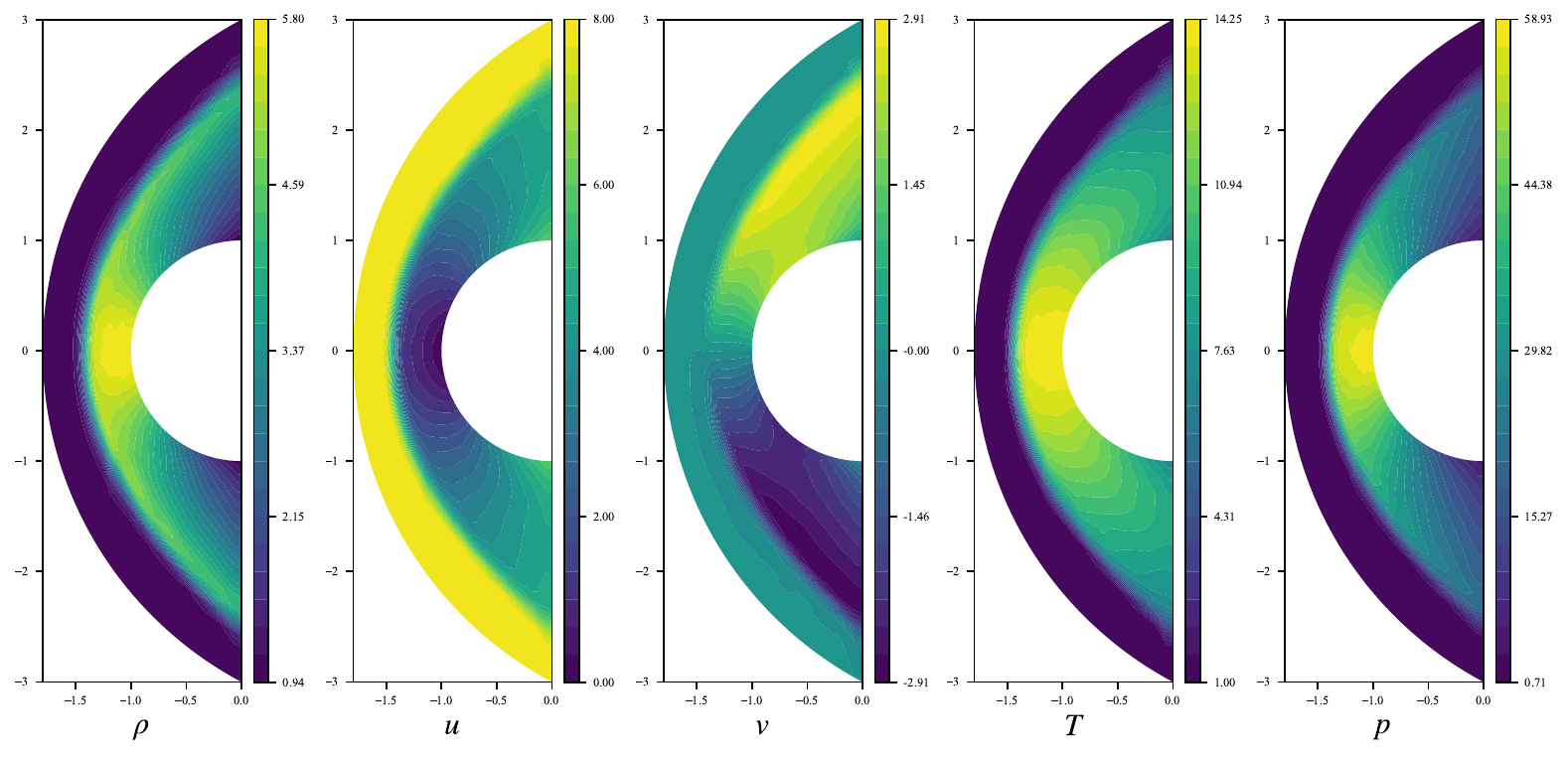}
        \caption*{(a) DNS}
    \end{minipage}
    \hspace{0.04\textwidth}
    \begin{minipage}{0.47\textwidth}
        \centering
        \includegraphics[width=1\textwidth]{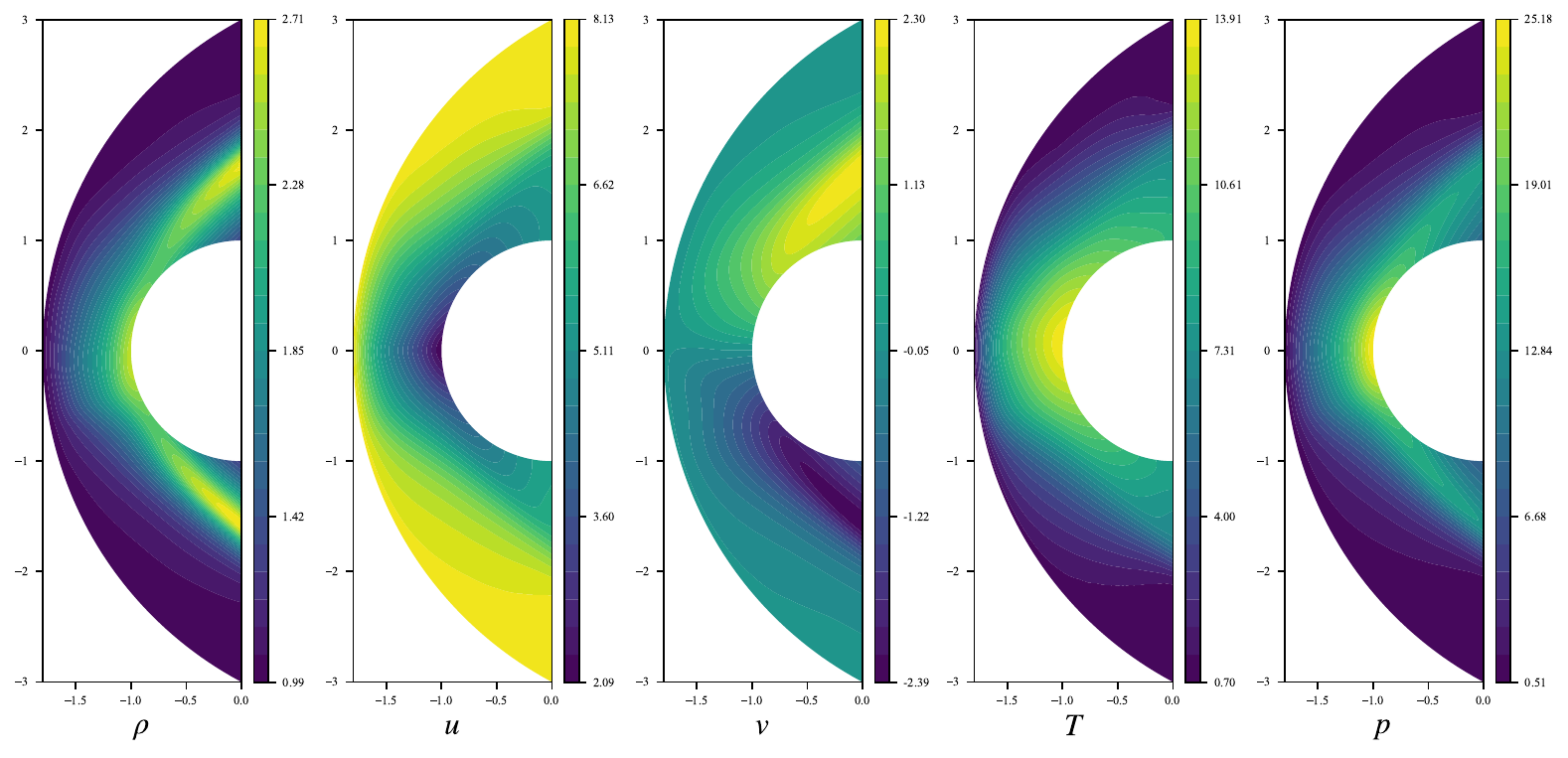}
        \caption*{(b) PINN\textsubscript{c}}
    \end{minipage}
    \begin{minipage}{0.47\textwidth}
        \centering
        \includegraphics[width=1\textwidth]{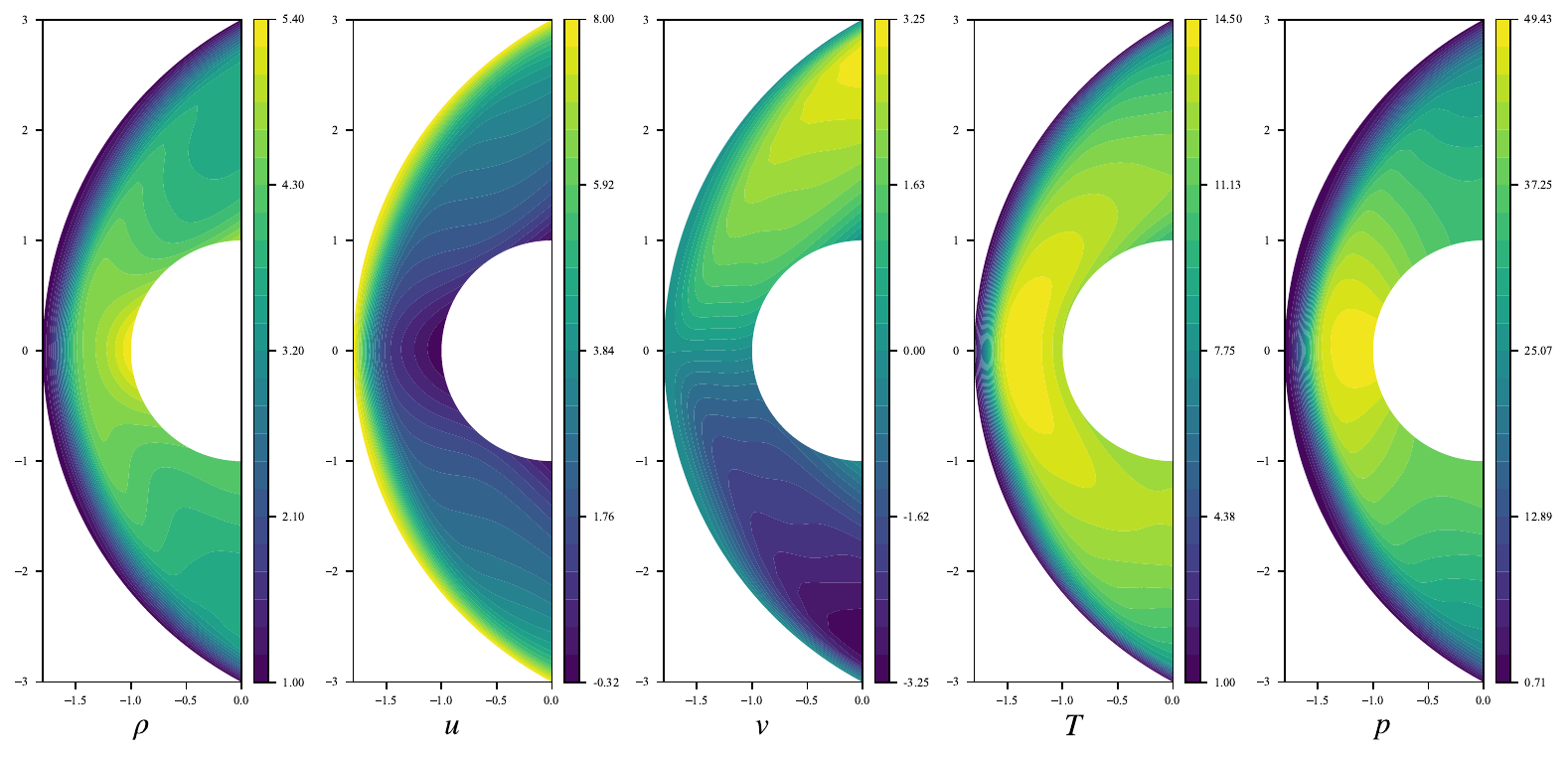}
        \caption*{(c) PINN\textsubscript{d}}
    \end{minipage}
    \hspace{0.04\textwidth}
    \begin{minipage}{0.47\textwidth}
        \centering
        \includegraphics[width=1\textwidth]{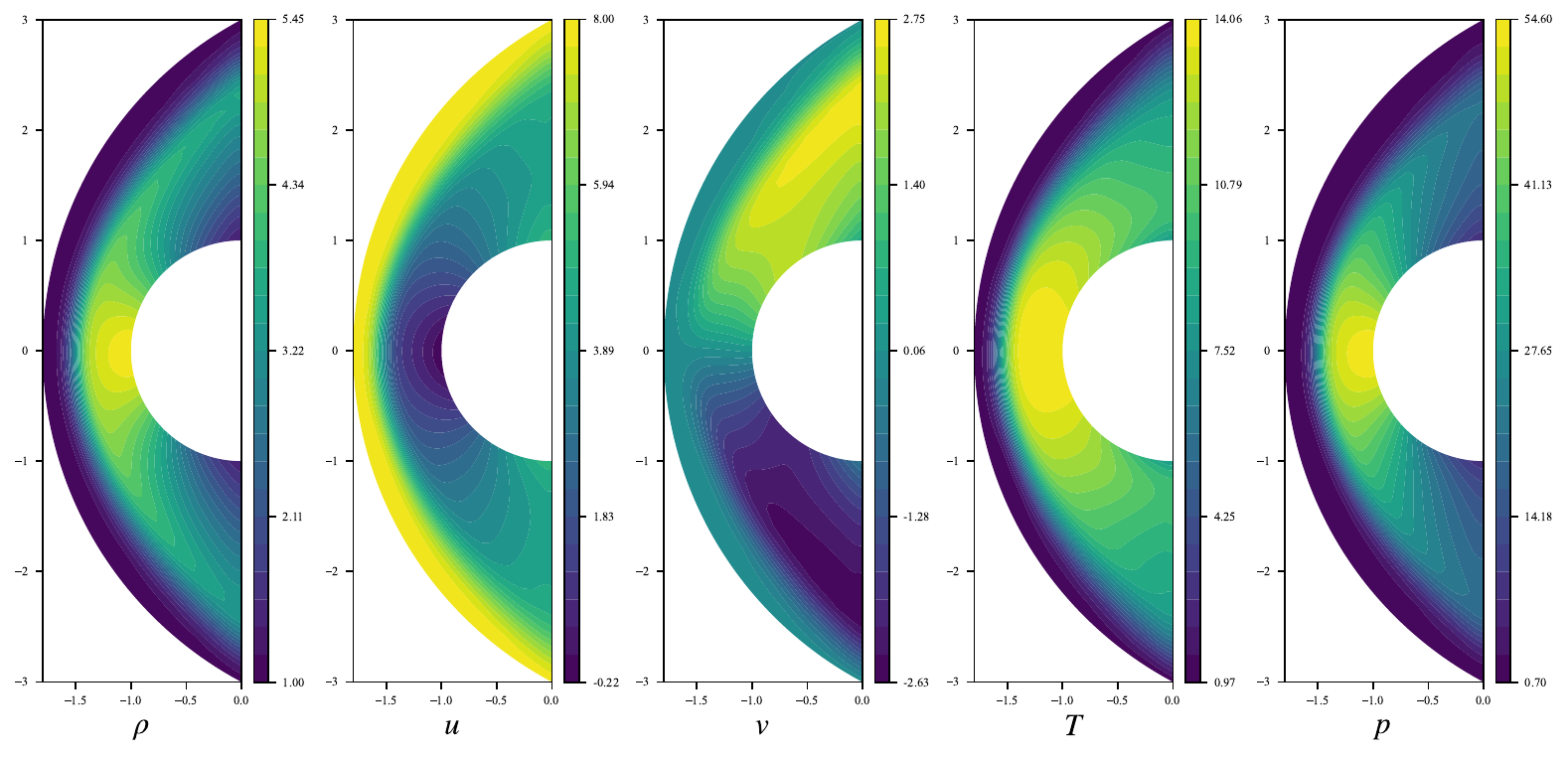}
        \caption*{(d) PIGNN}
    \end{minipage}
    \caption{Flow-field comparisons of different methods for hypersonic inviscid cylinder flow at $\mathrm{Ma}=8$.}
    \label{fig:bowshock_comparison}
\end{figure}

Hypersonic flows feature a wide dynamic range of compressible flow variables, and hard constraints based on distance functions do not behave as stably as in previous test cases when using raw coordinate inputs. We thus normalize input coordinates for stable training. For discretization-constrained neural solvers, a stable first order FV scheme with Rusanov flux is adopted, consistent with DNS. Representative flow-field comparisons are presented in \cref{fig:bowshock_comparison}.
The PINN\textsubscript{c} becomes trapped in a sub-optimal optimization state. Symmetric shock-like structures form at the cylinder shoulders yet cannot propagate toward the front stagnation region. Extending training further degrades predictions, even without anomalous loss history signals. Discretization constrained models capture the general flow-field topology but yield heavily smeared shocks. Specifically, PINN\textsubscript{d} shows multiple regions of unphysical reverse flow in horizontal velocity $u$ and fails to reproduce expansion features at upper and lower shoulders. In contrast, PIGNN performs much better. Despite local mild unphysical reverse flow, it recovers major expansion region flow features.

\begin{figure}[!ht]
    \centering
    \begin{minipage}{0.47\textwidth}
        \centering
        \includegraphics[width=1\textwidth]{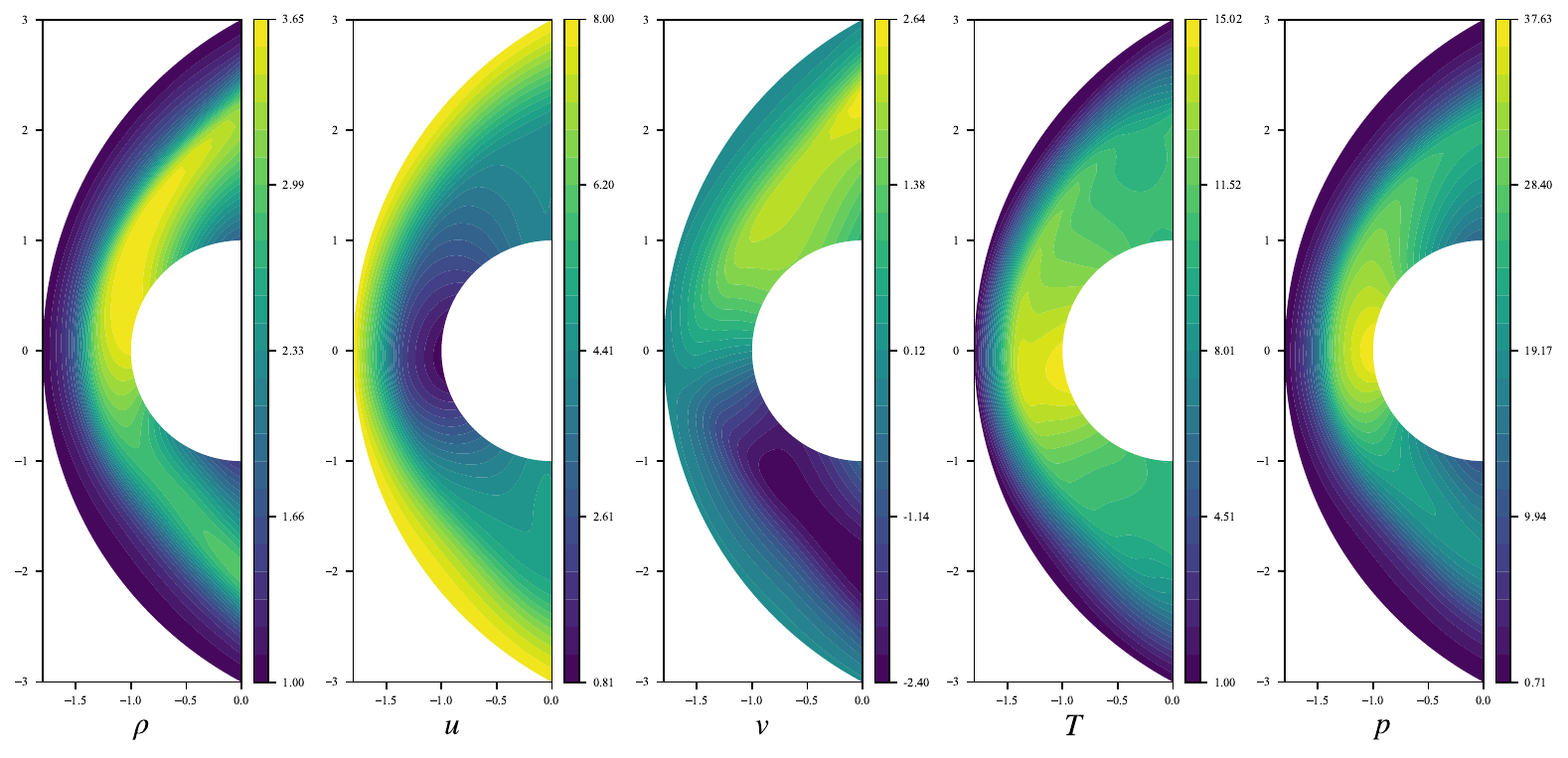}
        \caption*{(a) PINN\textsubscript{d}}
    \end{minipage}
    \hspace{0.04\textwidth}
    \begin{minipage}{0.47\textwidth}
        \centering
        \includegraphics[width=1\textwidth]{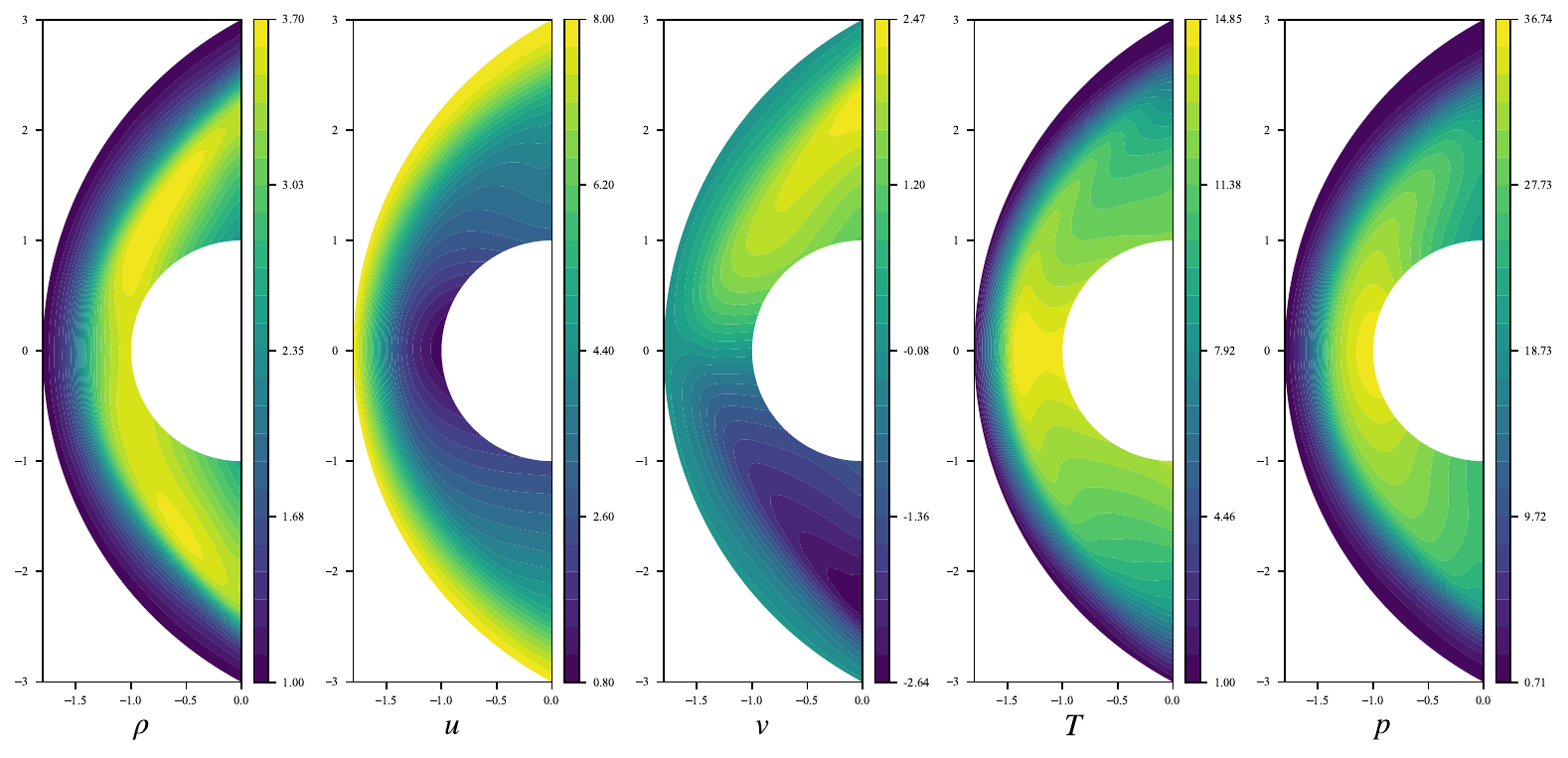}
        \caption*{(b) PIGNN}
    \end{minipage}
    \caption{Sub-optimal state induced by MUSCL reconstruction for discretization-constrained models: hypersonic inviscid cylinder flow at $\mathrm{Ma}=8$, trained for 30000 epochs.}
    \label{fig:bowshock_muscl}
\end{figure}

Similar to the observations in the Poisson cases, higher order schemes suffer from instability under strong nonlinearity. We introduce MUSCL reconstruction here, which brings oscillatory loss curves and impedes convergence. Even when trained for 30000 epochs, the models still get stuck in a sub-optimal state, as illustrated in \cref{fig:bowshock_muscl}. Other low-dissipation Riemann flux schemes produce similar behaviour, whereas the high-dissipation Rusanov flux yields more robust training. Such behaviour closely resembles the previously discussed characteristics of AD-based solvers, yet AD-based solvers exhibit far more severe instability. This demonstrates the inherent instability associated with high-order-like neural network training in stiff problems.

\begin{figure}[!ht]
    \centering
    \begin{minipage}{0.3225\textwidth}
        \centering
        \includegraphics[width=1\textwidth]{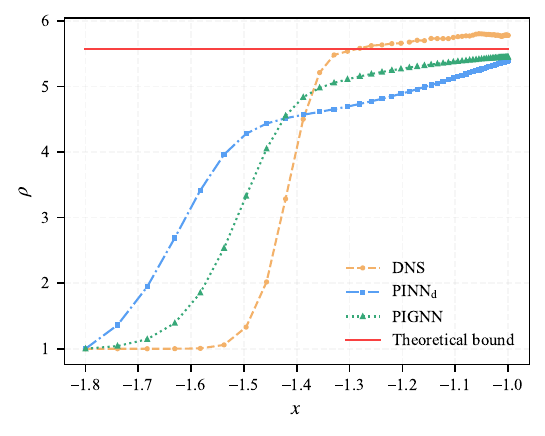}
        \caption*{(a)}
    \end{minipage}
    \begin{minipage}{0.33\textwidth}
        \centering
        \includegraphics[width=1\textwidth]{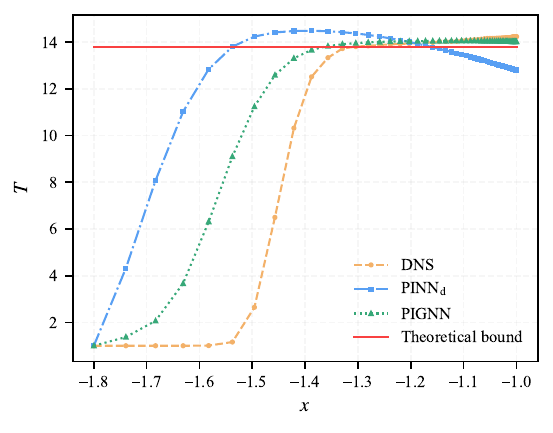}
        \caption*{(b)}
    \end{minipage}
    \begin{minipage}{0.33\textwidth}
        \centering
        \includegraphics[width=1\textwidth]{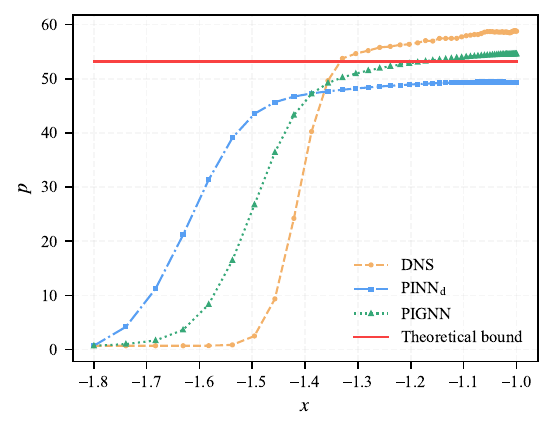}
        \caption*{(c)}
    \end{minipage}
    \caption{Profiles along the stagnation streamline ($y=0$): (a) density $\rho$; (b) temperature $T$; (c) pressure $p$.}
    \label{fig:stagnation_line}
\end{figure}

\begin{figure}[!ht]
    \centering
    \begin{minipage}{0.3225\textwidth}
        \centering
        \includegraphics[width=1\textwidth]{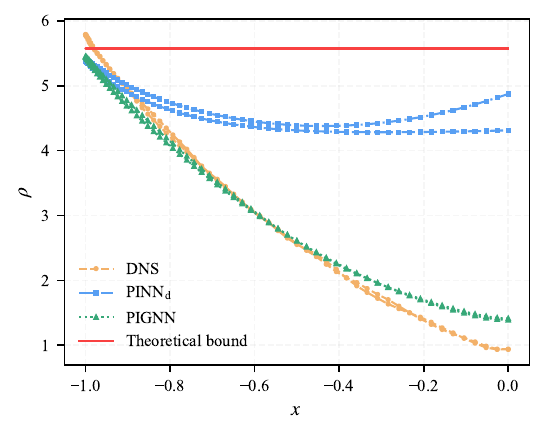}
        \caption*{(a)}
    \end{minipage}
    \begin{minipage}{0.33\textwidth}
        \centering
        \includegraphics[width=1\textwidth]{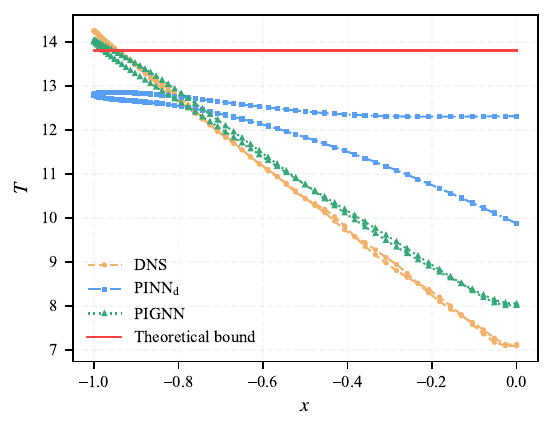}
        \caption*{(b)}
    \end{minipage}
    \begin{minipage}{0.33\textwidth}
        \centering
        \includegraphics[width=1\textwidth]{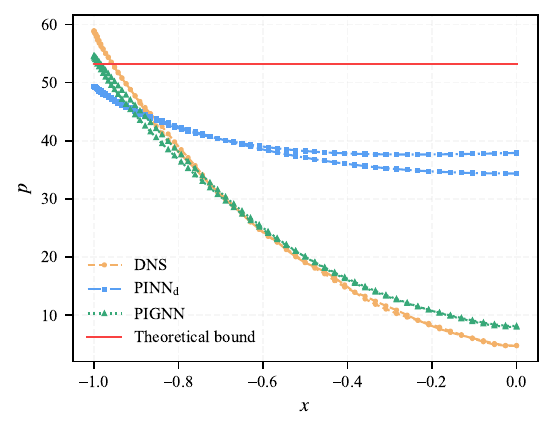}
        \caption*{(c)}
    \end{minipage}
    \caption{Wall-surface distributions for the cylinder: (a) density $\rho$; (b) temperature $T$; (c) pressure $p$.}
    \label{fig:surface_distribution}
\end{figure}

Considering shock-dominated problems, the achievable approximation bounds of neural networks are substantially inferior to those observed in the preceding test cases. Consequently, global metrics such as the $L_2$ error evaluated against high-fidelity reference solutions lose discriminative power and become less meaningful for assessment. We therefore focus on local flow details of practical engineering interest, namely the stagnation-line profile along $y=0$ and the wall-surface distribution. Since the results obtained by PINN\textsubscript{c} are of insufficient quality for fair comparison, it is excluded from further analysis. The stagnation-line and wall-surface distributions are presented in \cref{fig:stagnation_line} and \cref{fig:surface_distribution}, respectively.
The deficiencies of PINN\textsubscript{d} within the fan-shaped expansion regions shown in \cref{fig:bowshock_comparison} become more evident in these local profiles. PINN\textsubscript{d} exhibits an immediate steep rise of gradients starting from the inflow boundary, and the temperature peak is shifted upstream away from the wall. Although wall-side quantities lie within the theoretical bounds derived from normal-shock relations \parencite{2003_anderson_shockrelations}, the predictions deviate markedly from the physically expected trend of monotonic decrease when moving away from the stagnation point toward both sides. By contrast, despite observable shock smearing in PIGNN, its gradients grow gradually from the inflow. The wall-surface distributions agree reasonably well with DNS results and maintain favorable symmetry. While PIGNN cannot deliver high-fidelity engineering grade solutions, it still provides indicative results for preliminary engineering assessments.

\begin{figure}[!ht]
    \centering
    \begin{minipage}{0.3025\textwidth}
        \centering
        \includegraphics[width=1\textwidth]{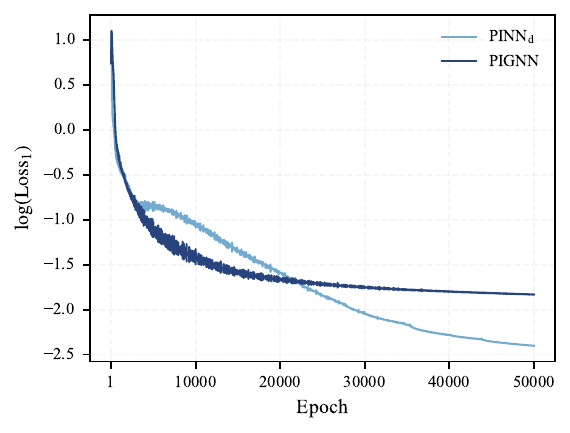}
        \caption*{(a)}
    \end{minipage}
    \begin{minipage}{0.3025\textwidth}
        \centering
        \includegraphics[width=1\textwidth]{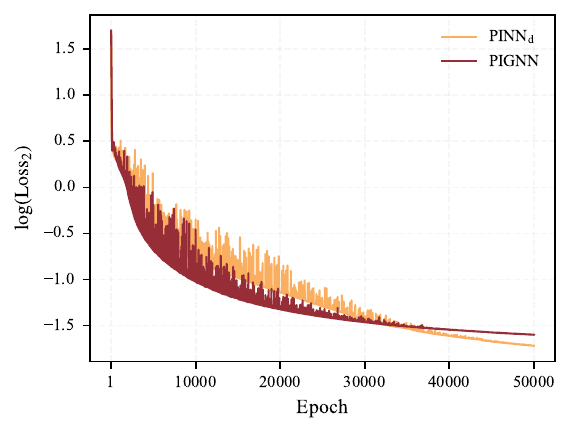}
        \caption*{(b)}
    \end{minipage}
    \begin{minipage}{0.3825\textwidth}
        \centering
        \includegraphics[width=1\textwidth]{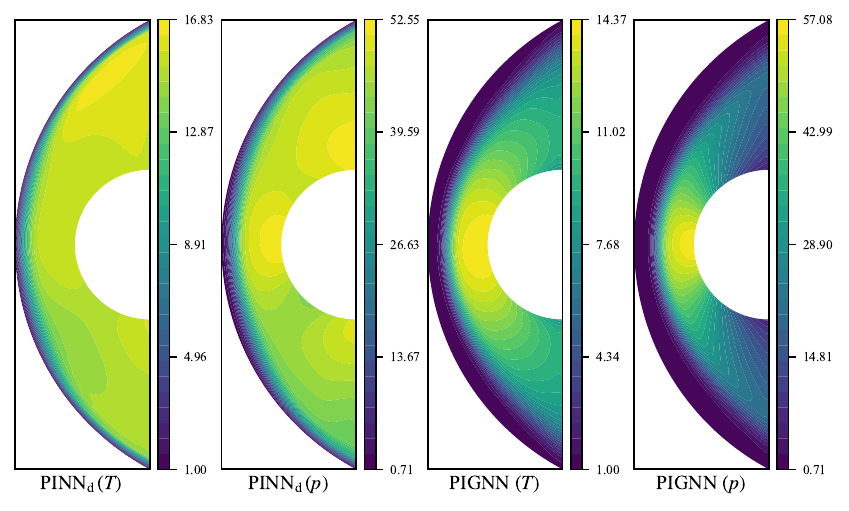}
        \caption*{(c)}
    \end{minipage}
    \caption{Illusory convergence of PINN\textsubscript{d}: (a) Soft-constraint BC residual loss history; (b) PDE residual loss history; (c) Flow‑field comparison of PINN\textsubscript{d} and PIGNN trained for $50000$ epochs.}
    \label{fig:bowshock_convergence}
\end{figure}

To investigate whether extended training can mitigate shock smearing, we increase the training epochs up to $50000$. The results reveal that shock smearing constitutes an inherent limitation rather than an artifact of insufficient training. Nevertheless, excessive training exerts adverse effects on PINN\textsubscript{d}, which exhibits unanticipated illusory convergence behaviour similar to PINN\textsubscript{c}, as illustrated in \cref{fig:bowshock_convergence}. Although PINN\textsubscript{d} achieves substantially lower boundary-condition and PDE residuals compared with PIGNN, it produces severely over-smoothed solutions within the shock layer. In contrast, PIGNN yields no further improvement with prolonged training yet stays close to its achievable performance bound, accompanied by relatively stable loss oscillations. Such deceptive convergence behaviour can also occur under discretization-constrained formulations, demonstrating that it is not a problem exclusive to AD-based paradigms and may be linked to intrinsic optimisation properties of MLP networks.
This contrasts with our earlier Poisson test-case observations. As discussed previously, GNNs tend to exhibit oscillatory or diverging loss under poor learning performance arising from inappropriate or ill-suited model configurations, which acts as an explicit diagnostic warning signal to users. By comparison, the deceptive illusory convergence observed here is far more perilous for practical engineering practice: the loss decreases favourably without any obvious warning, while the underlying physical predictions instead degrade toward worse solutions.

\begin{figure}[!ht]
    \centering
    \begin{minipage}{0.47\textwidth}
        \centering
        \includegraphics[width=1\textwidth]{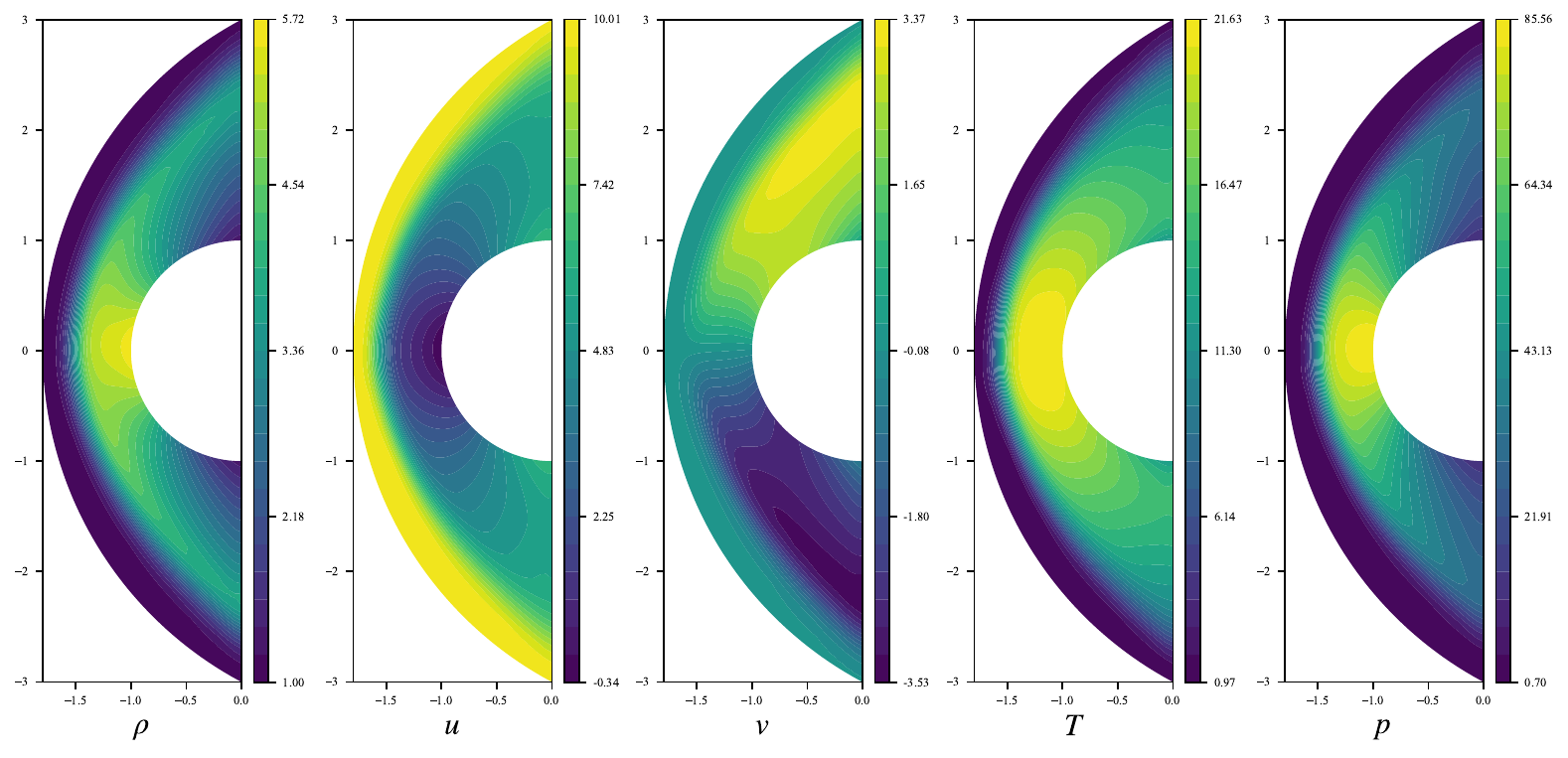}
        \caption*{(a) $\mathrm{Ma}=10$}
    \end{minipage}
    \hspace{0.04\textwidth}
    \begin{minipage}{0.47\textwidth}
        \centering
        \includegraphics[width=1\textwidth]{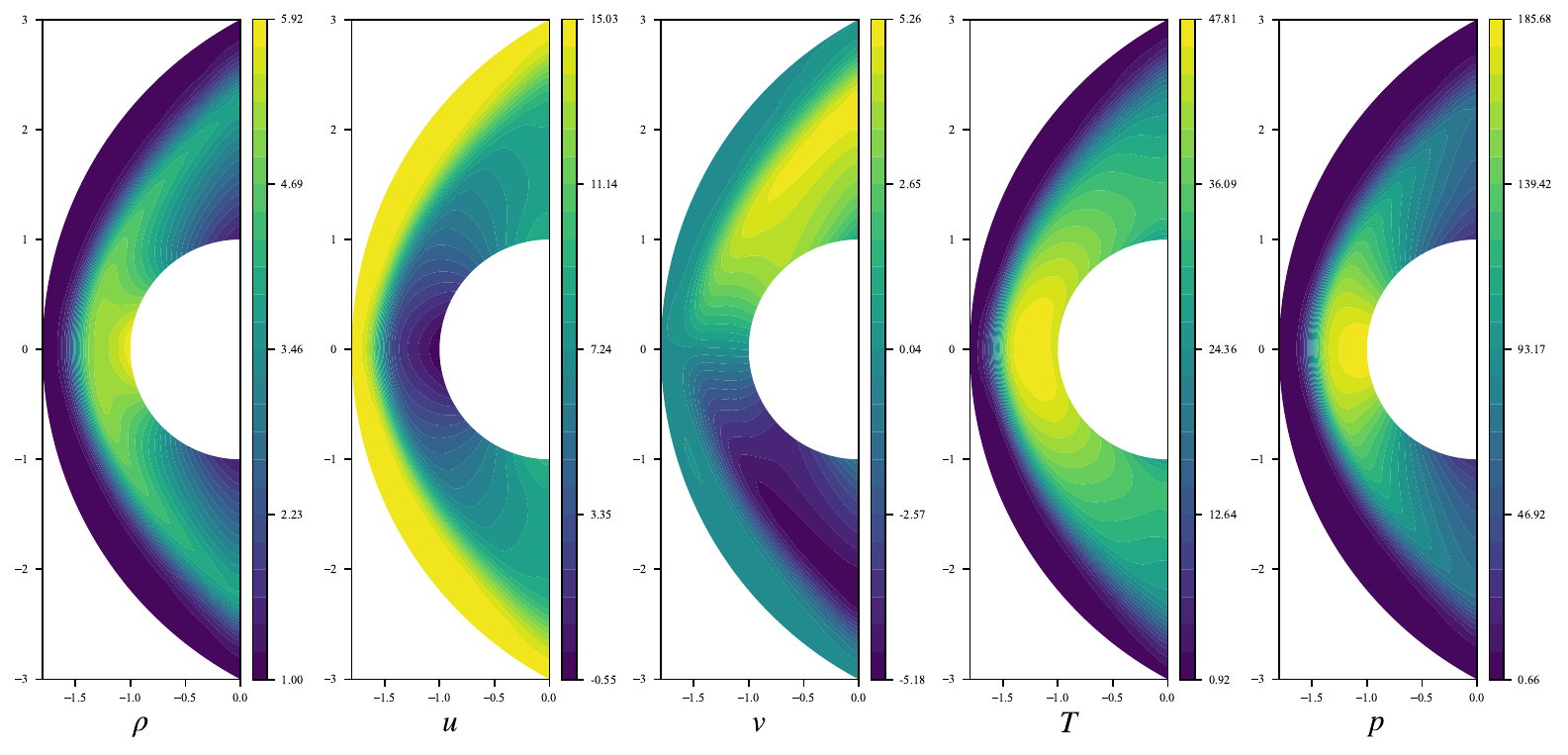}
        \caption*{(b) $\mathrm{Ma}=15$}
    \end{minipage}
    \begin{minipage}{0.47\textwidth}
        \centering
        \includegraphics[width=1\textwidth]{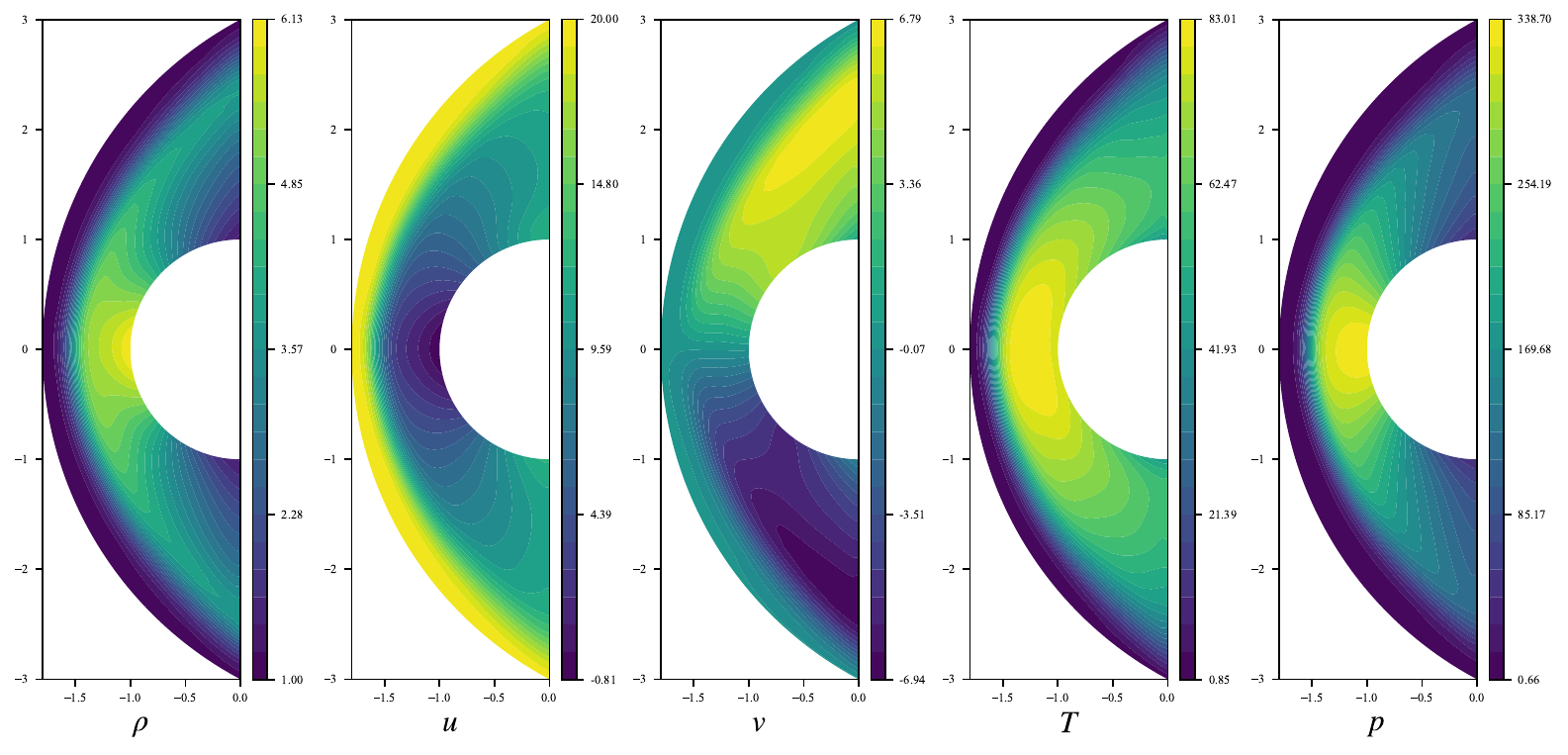}
        \caption*{(c) $\mathrm{Ma}=20$}
    \end{minipage}
    \hspace{0.04\textwidth}
    \begin{minipage}{0.47\textwidth}
        \centering
        \includegraphics[width=1\textwidth]{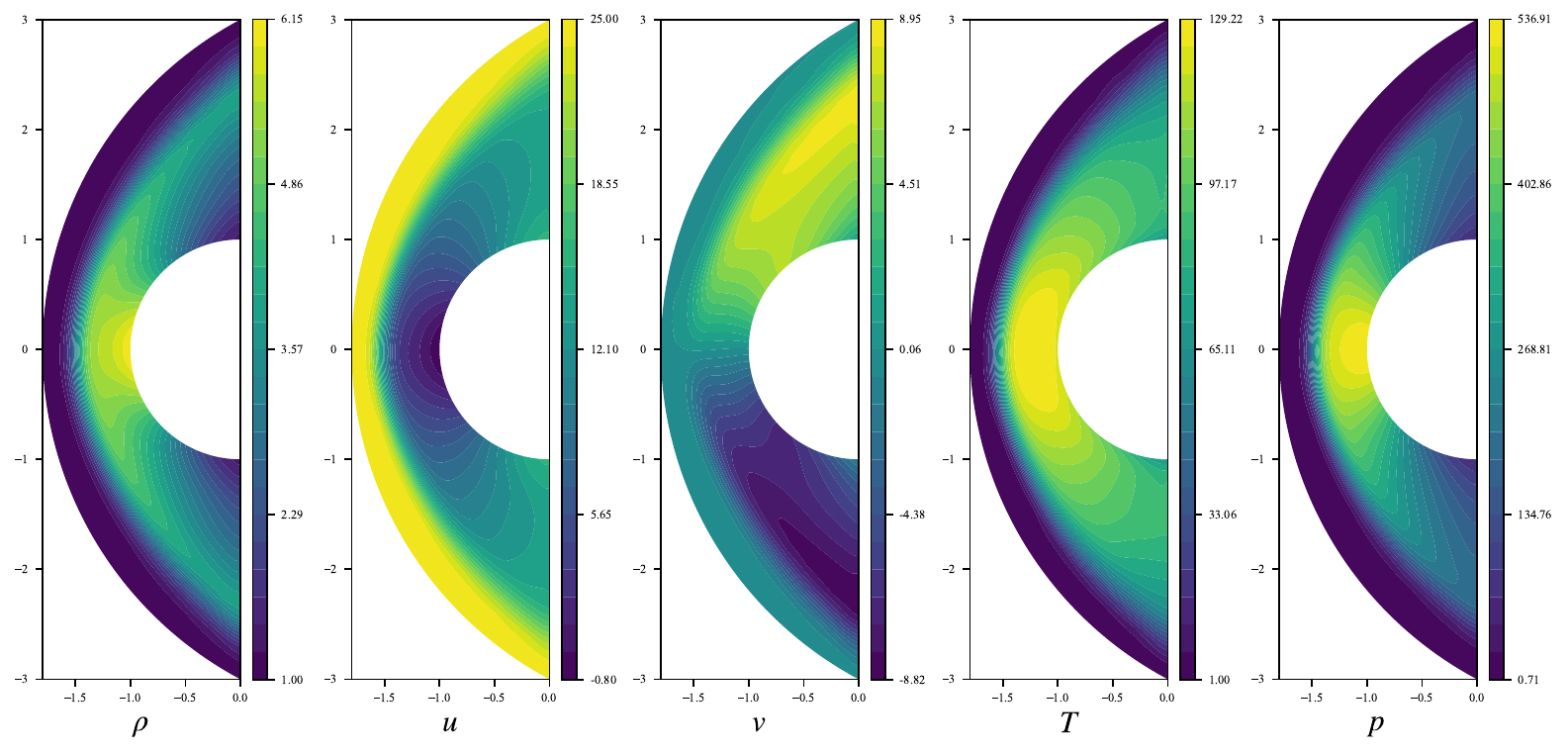}
        \caption*{(d) $\mathrm{Ma}=25$}
    \end{minipage}
    \caption{PIGNN performance for hypersonic inviscid cylinder flow at higher Mach numbers. From low to high Mach number, each case is trained with $20000$, $30000$, $30000$, and $40000$ epochs.}
    \label{fig:higher_mach}
\end{figure}

We further verify the robustness of PIGNN across increasing Mach number configurations. Higher Mach cases demand longer training, with 20000, 30000, 30000 and 40000 epochs used for Mach of 10, 15, 20 and 25, respectively, as illustrated in \cref{fig:higher_mach}. PIGNN remains stable even at a Mach number as high as 25. Although a slight upstream shift of the temperature peak emerges at elevated Mach numbers, the shock location does not move appreciably forward with increasing Mach number. Overall, the flow field predictions stay close to theoretical values, demonstrating PIGNN’s value as a coarse solver for rapid flow-field prediction.
We have also attempted simulations of viscous compressible NS flows but failed to achieve desirable outcomes. The no-slip boundary condition for viscous flows induces shock smearing-like artifacts for large velocity gradients near the wall. The combined smearing effects within both the near-wall region and the shock layer create substantial challenges for neural network training. By contrast, the stagnation region in Euler computations also features large incoming velocity gradients. Under the slip wall boundary condition, the neural network stabilizes the near wall gradient by stretching sharp velocity transitions via small amounts of permitted non-physical backflow during training.

\subsection{Efficiency analysis}\label{sec:efficiency}
Most existing studies evaluate efficiency within a limited scope, typically counting only training epochs or iterations while overlooking other additional overhead. From an engineering viewpoint, however, as highlighted in \parencite{2024_nmi_weakbaselines}, efficiency should be interpreted in a broader sense. Training steps and wall-clock time constitute only measurable parts of the overall cost. Additional implicit overhead, including preprocessing steps such as mesh generation and discretization, should also be accounted for in a complete efficiency assessment.
Meanwhile, the applicable boundary of efficiency improvement needs to be clarified. For instance, higher-order schemes can achieve results comparable to lower-order schemes at much lower resolutions. Nevertheless, such efficiency gains come at the cost of degraded stability for more complex problems. To some extent, robustness is even more critical than pure accuracy or efficiency in practical engineering applications, where numerical uncertainty poses major practical concerns. Therefore, computational efficiency cannot be evaluated in isolation, but should be comprehensively assessed from multiple dimensions involving accuracy, robustness, and practical computational overhead.

It is well-acknowledged that mesh generation and discretization can be fairly labour-intensive tasks. For simple problems such as the Poisson equation and when extreme accuracy is not required, PINN\textsubscript{c}'s mesh-free nature offers clear efficiency advantages even though DNS can produce solutions within seconds; thus, we do not elaborate further on this point. For fair comparison, we test the NS and Euler flow problems under identical hardware. In \cref{tab:training_efficiency}, we report wall-clock GPU computation time obtained on an NVIDIA GeForce RTX 4070 Laptop GPU with 8 GB of memory, under the training configurations summarized in \cref{tab:training_params}.
It can be seen that, for the lid-driven cavity and backward-facing step cases under identical training configurations, the growth in computation time of PIGNN with increasing mesh resolution even outpaces that of the AD processing in PINN\textsubscript{c}. This arises because the number of topological edges increases drastically as the mesh is refined. Consequently, PIGNN offers little practical merit when it provides no substantial accuracy improvement. For the hypersonic case, by contrast, MLP-based models suffer from severe numerical instability and deceptive convergence; hence their higher computational efficiency alone does not make them acceptable for rapid engineering assessment. Under the same challenging flow conditions, PIGNN retains stable training while achieving computation time comparable to that of explicit schemes.
Note that the hypersonic inviscid flow case consists of 6480 cells ($120 \times 54$) with four output variables. Despite having more state-variable dimensions, its per-epoch training time is shorter compared with the incompressible NS cases featuring only three output dimensions. This discrepancy arises from implementation differences in discrete constraint computation. The FD-based constraints for incompressible NS involve multiple matrix-based differential operators and numerous block-matrix arithmetic operations. By contrast, for FV constraints, face-fluxes are computed for each equation component, followed by only a single sparse scatter-add operation to assemble the residual, which reduces computational overhead. This indicates that, under comparable problem scales, FV-based constraints can yield better training efficiency. However, this does not imply that FD formulations are without merits. For regular uniform computational domains, FD discretization is comparatively straightforward: it does not require extra storage for topological information, face normals and related mesh metadata, and higher-order operator constraints can be implemented more conveniently.

\begin{table}[!ht]
    \centering
    \caption{Wall-clock GPU computation time for different solution methods.}
    \label{tab:training_efficiency}
    \small
    \begin{tabularx}{0.8\textwidth}{@{\hspace{0.5em}}lXXXXX@{\hspace{0.5em}}}
        \toprule
        \multirow{2}{*}{PDE} & \multirow{2}{*}{Resolution} & \multicolumn{4}{c}{Methods} \\
        \cmidrule{3-6}
        & & DNS & PINN\textsubscript{c} & PINN\textsubscript{d} & PIGNN \\
        \midrule
        \multirow{4}{*}{Lid-driven cavity flow}
        & $50 \times 50$ & $11.92 \, \mathrm{s}$ & $6.09  \, \mathrm{min}$ & $3.75 \, \mathrm{min}$ & $4.99 \, \mathrm{min}$ \\
        & $100 \times 100$ & $12.54 \, \mathrm{s}$ & $6.76 \, \mathrm{min}$ & $3.94 \, \mathrm{min}$ & $7.29 \, \mathrm{min}$ \\
        & $150 \times 150$ & $13.82 \, \mathrm{s}$ & $9.67 \, \mathrm{min}$ & $4.16 \, \mathrm{min}$ & $11.16 \, \mathrm{min}$ \\
        & $200 \times 200$ & $16.71 \, \mathrm{s}$ & $14.56 \, \mathrm{min}$ & $4.61 \, \mathrm{min}$ & $17.80 \, \mathrm{min}$ \\
        \midrule
        \multirow{3}{*}{Backward-facing step flow}
        & $100 \times 50$ & $10.14\ \mathrm{s}$ & $5.45 \, \mathrm{min}$ & $3.89 \, \mathrm{min}$ & $5.28 \, \mathrm{min}$ \\
        & $200 \times 100$ & $14.49 \, \mathrm{s}$ & $6.77 \, \mathrm{min}$ & $3.92 \, \mathrm{min}$ & $8.90 \, \mathrm{min}$ \\
        & $300 \times 150$ & $16.75 \, \mathrm{s}$ & $11.21 \, \mathrm{min}$ & $4.42 \, \mathrm{min}$ & $15.52 \, \mathrm{min}$ \\
        \midrule
        Hypersonic inviscid flow & / & $2.57 \, \mathrm{min}$ & / & $1.19 \, \mathrm{min}$ & $2.28 \, \mathrm{min}$ \\
        \bottomrule
    \end{tabularx}
\end{table}

It should be emphasized that the reported wall-clock times serve only as a fair reference under identical training configurations. Shorter runtime does not equate to absolute superiority, and practical performance is governed by multiple trade-offs.
For low Reynolds number incompressible NS problems, the mesh-free property of PINN\textsubscript{c} remains a distinct advantage. Nevertheless, such efficiency gains hold only under weak-nonlinearity conditions, and the failure boundaries of such efficiency gains should be carefully acknowledged. Additionally, mesh-based discretization-constrained neural solvers possess unique merits unavailable to DNS. They enable stable training on collocated grids without extra treatments to suppress numerical oscillations. Compared with full staggered-grid implementations, the discretization workflow is substantially simpler. Staggered-grid approaches require separate discrete coefficient matrices for $u$, $v$ and $p$, together with elaborate and error-prone boundary-condition treatments. All of these factors introduce considerable implicit overhead that cannot be captured by wall-clock measurements. The notable strength of DNS lies in its traceable and controllable error sources, rendering it irreplaceable for final high-fidelity engineering evaluation. Accordingly, efficiency comparisons for neural solvers should focus on practical benefit-to-cost ratios. Approaches that introduce extra computational cost while delivering only limited performance gains are undesirable for engineering practice.
\FloatBarrier

\section{Conclusion and discussion}\label{sec:conclusion}
Our systematic analysis work is grounded in the practical perspective of engineering applications. Rather than benchmarking configurations to crown the best performers, this work examines three key axes that fundamentally shape neural solver behaviour: network architecture, constraint formulation, and problem regime. The aim is not to produce a ranking, but to interrogate how these axes interact to determine accuracy, efficiency, and robustness, and to extract transferable insights into when and why specific design choices succeed or fail. In doing so, we seek to offer engineers a principled basis for method selection in settings where a posteriori verification is unavailable. Through a systematically constructed spectrum of test cases with progressively escalating nonlinearity and boundary complexity, and with a rigorous error decomposition framework, the principal findings are summarized as follows:
\begin{itemize}
    \item AD-based formulations demonstrate clear advantages in weakly nonlinear regimes where boundary conditions permit global hard constraint enforcement and optimization remains well-conditioned, achieving high accuracy at relatively low sampling resolutions with the additional benefit of being mesh-free. However, this advantage is bounded: as nonlinearity strengthens, AD-based models become increasingly susceptible to illusory convergence and abrupt degradation. Engineering practice therefore requires explicit delineation of the applicable boundary of such efficiency gains, for which nonlinearity strength offers a practical and transferable criterion, whether manifested through source-term stiffness, Reynolds number, or analogous measures of coupling intensity in other PDE systems.
    \item For discretization-based formulations, the dominant error source shifts from truncation error toward approximation and optimization errors as nonlinearity strengthens. Higher-order schemes exert a positive effect only in weakly nonlinear regimes at low mesh resolutions; their benefits diminish and may even reverse under stronger nonlinearity, where approximation and optimization errors grow at a faster rate. The averaging nature of MSE-based optimization imposes a practical ceiling on achievable accuracy, such that performance gains from mesh refinement stagnate and eventually degrade beyond a critical resolution, with this threshold decreasing as nonlinearity increases.
    \item GNN-based models, by exploiting topological neighbourhood information, maintain superior robustness under strong nonlinearity and complex boundary conditions, and their loss trajectories offer a more reliable diagnostic signal than those of MLPs, which are prone to deceptive convergence without visible precursors. In weakly nonlinear problems of low complexity, however, this added aggregation capacity brings no appreciable accuracy benefit and introduces unnecessary computational overhead.
    \item Extrapolation-type boundary conditions, introduced in classical solvers for completeness, act on neural training as static, deterministic perturbations akin to truncation error. Their role is regime-dependent: in weakly nonlinear problems they add marginal unnecessary errors, while in strongly nonlinear, highly ill-conditioned cases they materially accelerate convergence. They should therefore be treated as problem-dependent regularization rather than universally beneficial or detrimental.
\end{itemize}

Our findings carry broader implications for the deployment of PINN-type methods in engineering practice. The persistent pursuit of accuracy heavily reliant on a-posteriori verification is neither sufficient nor constructive. What is more urgently needed is clear delineation of the applicable regime for each method, alongside an honest assessment of the attainable performance ceiling within such regimes. Such delineation helps direct research efforts toward avenues that yield meaningful practical gains. For scenarios where neural solvers serve only for preliminary screening, the marginal improvements in accuracy must be carefully weighed against attendant computational costs; methodological choices should be guided by benefit-to-cost ratio instead of raw accuracy alone. Equally pressing is the development of generic and formulation-agnostic supplementary diagnostic evaluation metrics that remain informative when loss history loses its discriminating capability, especially under deceptive-convergence conditions with increasing nonlinearity. Lacking such metrics, the field remains tied to a-posteriori reference solutions. This dependency proves impractical for real-world engineering deployment. Evaluation of neural PDE solvers should therefore be reframed within a more comprehensive, engineering-oriented perspective that accounts for applicability boundaries, recognisable failure modes and practical cost-effectiveness.
\FloatBarrier


\section*{Disclosure statement}
No potential conflict of interest was reported by the authors.

\section*{Funding}
The authors acknowledge the Research Fund of National Key Laboratory of Aerospace Physics in Fluids (Grant No. KT-APF-2025-018).

\printbibliography

\appendix
\setcounter{figure}{0}
\renewcommand\thefigure{\Alph{section}\arabic{figure}}
\setcounter{table}{0}
\renewcommand\thetable{\Alph{section}\arabic{table}}
\setcounter{equation}{0}
\renewcommand\theequation{\Alph{section}\arabic{equation}}
\makeatletter
\@addtoreset{equation}{section}
\makeatother
\section*{Appendix}
\section{Details of numerical solvers} \label{appendix_a}
\begin{figure}[!ht]
    \centering
    \begin{minipage}{0.24\textwidth}
        \centering
        \includegraphics[width=1\textwidth]{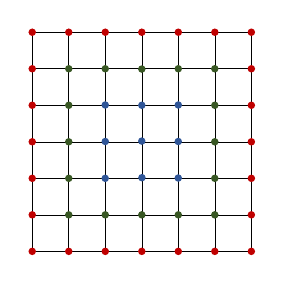}
        \caption*{(a)}
    \end{minipage}
    \begin{minipage}{0.24\textwidth}
        \centering
        \includegraphics[width=1\textwidth]{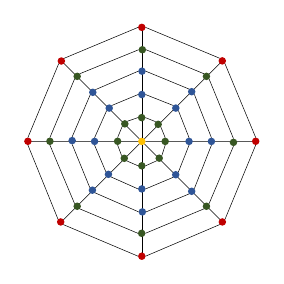}
        \caption*{(b)}
    \end{minipage}
    \begin{minipage}{0.24\textwidth}
        \centering
        \includegraphics[width=1\textwidth]{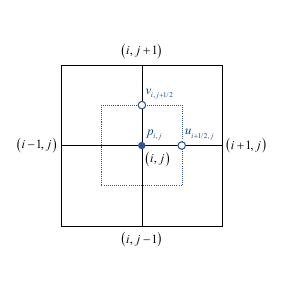}
        \caption*{(c)}
    \end{minipage}
    \begin{minipage}{0.24\textwidth}
        \centering
        \includegraphics[width=1\textwidth]{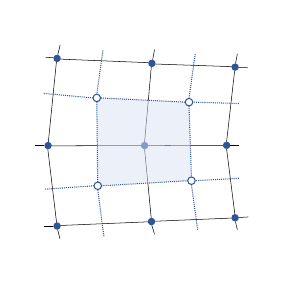}
        \caption*{(d)}
    \end{minipage}
    \caption{Schematic of discretization nodes and control volumes: (a) Cartesian nodal discretization; (b) Polar-coordinate nodal discretization; (c) Variable arrangement on staggered grid; (d) Dual control-volume for vertex-centered FV scheme.}
    \label{fig:mesh_discretization}
\end{figure}

To maintain consistency with the constraint formulations adopted in neural network solvers, no conventional ghost points or ghost cells are utilized throughout all numerical discretizations in the present work. 
\cref{fig:mesh_discretization} illustrates four distinct node layouts employed for Cartesian-grid standard FD, polar-coordinate FD, staggered-grid FD, and vertex-centered FV discretizations, respectively. Detailed descriptions for each discretization strategy are provided in the following subsections.

\subsection{FDM for Poisson-type equations}\label{appendix_a1}
For the Cartesian-grid Poisson test cases (Sine and Polynomial), we present the standard second-order and fourth-order central difference approximations for the Laplacian operator $\nabla^{2}u$. The second-order approximation reads
\begin{equation}\label{eq:laplacian_cartes}
    \nabla^{2}u_{i,j} = \frac{u_{i+1, j} - 2u_{i,j} + u_{i-1, j}}{\Delta x^{2}}  + \frac{u_{i, j+1} - 2u_{i,j} + u_{i, j-1}}{\Delta y^{2}} + \mathcal{O}(\Delta x^{2}) + \mathcal{O}(\Delta y^{2}).
\end{equation}
The fourth-order central difference stencil uses two neighbouring nodes on each side of the target node.
For the $x$-direction second-derivative approximation, the stencil coefficient vector along the $x$-axis reads $[-1, \, 16, \, -30, \, 16, \, -1] / (12\Delta x^2)$ corresponding to nodes $[u_{i+2,j}, \, u_{i+1,j}, \, u_{i,j}, \, u_{i-1,j}, \, u_{i-2,j}]$. The corresponding $y$-direction derivative follows the same coefficient pattern. This stencil yields a formal truncation error of $\mathcal{O}(\Delta x^{4})+\mathcal{O}(\Delta y^{4})$.

For the fourth-order scheme, the fourth-order stencil is applied at interior nodes (blue nodes in \cref{fig:mesh_discretization}(a)). Nodes adjacent to boundaries (green nodes) are treated with a second-order central difference stencil for order-reduction treatment, and Dirichlet boundary conditions are imposed directly at boundary nodes (red nodes). These difference operators cast the Poisson equation \cref{eq:poisson_formulation} into a sparse matrix system $\mathbf{A}u = f$, where $\mathbf{A}$ denotes the sparse coefficient matrix.
The linear Sine case is purely linear, we directly solve the resulting algebraic system using the direct sparse solver from {\tt scipy.sparse.linalg.spsolve}. The nonlinear Polynomial case is solved via standard Newton iteration.
At each iteration step, the residual function is defined as 
\begin{equation}\label{eq:residual_fun}
    g(u) = \mathbf{A}u - f(u),
\end{equation}
with the update given by
\begin{equation}\label{eq:newton_iteration}
    u_{n+1} = u_n - \alpha \big(g'(u_n)\big)^{-1} g(u_n),
\end{equation}
where $g'(u_n)$ denotes the Jacobian matrix $\mathbf{A}-f'(u)$ evaluated at $u_n$, and $\alpha \in [0, 1]$ is the relaxation parameter. Linearization at each iteration yields an algebraic system to be solved.

For the Liouville case, discretization is carried out on the uniform polar-coordinate mesh, with Laplacian operator given by
\begin{equation}\label{eq:laplacian_polar}
\nabla^2 u = \frac{\partial^2 u}{\partial r^2} + \frac{1}{r}\frac{\partial u}{\partial r} + \frac{1}{r^2}\frac{\partial^2 u}{\partial \theta^2},
\end{equation}
where the individual differential terms are discretized following analogous central difference procedures to the earlier presented Cartesian-coordinate case. Second-order and fourth-order central difference approximations are employed for the differential terms in polar coordinates. For the fourth-order scheme, order-reduction to second-order stencils is performed for nodes (marked green in \cref{fig:mesh_discretization}(b)) near both the outer boundary and the polar centre (orange). Boundary nodes (red) enforce Dirichlet boundary conditions identical to the previous setting. At the polar centre node, direct polar-coordinate FD stencils are ill-defined as $r\rightarrow 0$, where we adopt a special averaging discretization given by
\begin{equation}\label{eq:fdm_polar_center}
    \frac{4}{{\Delta r}^{2}}\left( u_{0} - \frac{1}{N_\theta} \sum_{k \in \mathcal{N}_{0}} u_{k}  \right) = 0,
\end{equation}
where $u_{0}$ denotes the numerical solution at the polar centre, $\mathcal{N}_{0}$ is the set of neighbouring nodes on the first ring surrounding the centre, and $N_\theta = |\mathcal{N}_{0}|$ is the cardinality of this set.

This discretization yields a residual definition $g(u)$ consistent with the foregoing. Conventional Newton iteration suffers instability for this stiff source-term problem; hence we adopt the higher-order Halley iteration. Specifically, we perform a second-order Taylor expansion of $g(u)$:
\begin{equation}\label{taylor_expansion}
    g(u) = g(u_{n}) + g'(u_{n}) \Delta u_{n} + \frac{1}{2}g''(u_{n}) \Delta u_{n}^{2} + \mathcal{O}(\Delta u_{n}^{3}).
\end{equation}
Directly solving for the update from this quadratic form is not feasible; the Halley iteration approximates the quadratic-term increment using the Newton step, which leads to the update formulation
\begin{equation}\label{eq:halley_iteration}
\Bigl(g'(u_{n})-\frac12\,g''(u_{n}) \Delta \hat{u}_{n}\Bigr) (u_{n+1}-u_{n}) = -g(u_{n}),
\end{equation}
where $\Delta \hat{u}_{n} = -\big(g'(u_{n})\big)^{-1}g(u_{n})$ is derived from the standard Newton update as in \cref{eq:newton_iteration}, and $g''(u) = -f''(u)$ denotes the Hessian matrix. Owing to this approximation, the nonlinear scheme achieves third-order convergence. Nevertheless, relaxation factors are still required in practical iterations to maintain stable convergence. For the test cases in \cref{sec:poisson}, relaxation factors of $0.2$, $0.2$, $0.3$, and $0.5$ are adopted for the four resolution levels from coarse to fine, respectively.

\subsection{Staggered-grid FDM for incompressible NS equations}\label{appendix_a2}
For the incompressible NS equations, conventional collocated-grid formulations are prone to pressure-velocity decoupling issues. Here we adopt the staggered-grid arrangement, as illustrated in \cref{fig:mesh_discretization}(c). The horizontal velocity $u$ is stored at half-grid nodes along vertical grid lines, the vertical velocity $v$ is stored at half-grid nodes along horizontal grid lines, and pressure $p$ is defined at the primary nodal positions. 
Then, \cref{eq:ns_formulation} is discretized into
\begin{equation}\label{eq:discrete_ns}
\begin{array}{lll}
    \big( uu'_{x} + vu'_{y} + p'_x - \frac{1}{\mathrm{Re}}(u''_{xx} + u''_{yy}) \big) \big|_{i+\frac{1}{2},\,j} = 0 \\[2pt]
    \big( uv'_{x} + vv'_{y} + p'_y - \frac{1}{\mathrm{Re}}(v''_{xx} + v''_{yy}) \big) \big|_{i,\,j+\frac{1}{2}} = 0 \\[2pt]
    \big(u'_{x} + v'_{y}\big)\big|_{i,\,j} = 0
\end{array}.
\end{equation}

Our test cases in \cref{sec:lid_driven} and \cref{sec:bward_step} adopt moderate Reynolds number values, so standard central difference schemes are directly employed for spatial discretization. Note that the pressure gradients appearing in momentum equations and the velocity gradients at the primary nodal position for the continuity equation are evaluated via half-grid difference operators (equivalent to standard second-order FD approximations), given by
\begin{equation}\label{eq:half_grid_diff}
\begin{aligned}
    & p'_{x}\big|_{i+\frac{1}{2}, j} = \frac{p_{i+1,j} - p_{i, j}}{\Delta x}, \quad p'_{y}\big|_{i, j+\frac{1}{2}} = \frac{p_{i,j+1} - p_{i, j}}{\Delta y}, \\[2pt]
    & u'_{x}\big|_{i, j} = \frac{u_{i+\frac{1}{2},j} - u_{i-\frac{1}{2}, j}}{\Delta x}, \quad v'_{y}\big|_{i, j} = \frac{v_{i,j+\frac{1}{2}} - v_{i, j-\frac{1}{2}}}{\Delta y}.
\end{aligned}
\end{equation}

In such a staggered scheme, we require separate gradient operators for each field variable.
Let $\mathbf{M}_{\mathrm{cx}}$, $\mathbf{M}_{\mathrm{cy}}$ denote the horizontal and vertical gradient operators for velocity in momentum equations; $\mathbf{M}_{\mathrm{gx}}$, $\mathbf{M}_{\mathrm{gy}}$ the pressure gradient operators; and $\mathbf{M}_{\mathrm{dx}}$, $\mathbf{M}_{\mathrm{dy}}$ the velocity gradient operators for the continuity equation. Let $\mathbf{M}_{\mathrm{vu}}$ and $\mathbf{M}_{\mathrm{vv}}$ be the discrete viscous term matrices for the horizontal and vertical momentum equations, respectively, where boundary condition treatments are embedded within these matrices, and let $\mathbf{M}_{\mathrm{p}}$ denote the boundary matrix associated with pressure constraints. Define the solution vector $\mathbf{U} = [u,\,v,\,p]^{\top}$. The partial differential NS equations are cast into a nonlinear algebraic system $G(\mathbf{U}) = \mathbf{0}$, where
\begin{equation}\label{eq:staggered_algebra_ns}
    \left[\begin{array}{lll}
        \mathrm{diag}(u) \mathbf{M}_{\mathrm{cx}} + \mathbf{M}_{\mathrm{vu}} & \mathrm{diag}(v) \mathbf{M}_{\mathrm{cy}} & \mathbf{M}_{\mathrm{gx}} \\
        \mathrm{diag}(u) \mathbf{M}_{\mathrm{cx}} & \mathrm{diag}(v) \mathbf{M}_{\mathrm{cy}} + \mathbf{M}_{\mathrm{vv}} & \mathbf{M}_{\mathrm{gy}} \\
        \mathbf{M}_{\mathrm{dx}} & \mathbf{M}_{\mathrm{dy}} & \mathbf{M}_{\mathrm{p}}
    \end{array}\right]
    \begin{bmatrix}
        u\\v\\p
    \end{bmatrix}
    =\mathbf{0}.
\end{equation}

Different from standard FD that imposes Dirichlet constraints directly at boundary nodes, $u$ and $v$ are not stored at physical boundary locations on the staggered mesh, with only pressure $p$ assignable directly as $p_b$. A three-point polynomial boundary constraint derived from quadratic fitting is adopted for velocity boundary conditions, taking $u$ as an example,
\begin{equation}\label{eq:staggered_bc}
\begin{array}{ll}
    \frac{15}{8} u(x_{\mathrm{bl}} + \frac{1}{2} \Delta x) - \frac{10}{8} u(x_{\mathrm{bl}} + \frac{3}{2} \Delta x) + \frac{3}{8} u(x_{\mathrm{bl}}+\frac{5}{2} \Delta x)  = u(x_{\mathrm{bl}}), \\[2pt]
    \frac{3}{8} u(x_{\mathrm{br}} - \frac{3}{2} \Delta x) + \frac{6}{8} u(x_{\mathrm{br}} - \frac{1}{2} \Delta x) -\frac{1}{8} u(x_{\mathrm{br}} + \frac{1}{2} \Delta x) = u(x_{\mathrm{br}}).
\end{array}
\end{equation}
These coefficients are assembled directly as rows into the sparse coefficient matrix, and the vertical velocity component for the bottom and top boundary is treated analogously. Additionally, the prescribed physical boundary conditions, illustrated in \cref{eq:lid_driven_bc} and \cref{eq:backward_step_bc}, do not yield a perfectly well-posed discrete system and lead to rank-deficiency for the implicit formulation. Therefore, supplementary Neumann-type conditions are required. Specifically, for no-slip walls, the condition $\frac{\partial p}{\partial \mathbf{n}} = 0$ is enforced. For the outlet boundary, $\frac{\partial u}{\partial \mathbf{n}} = 0$ and $\frac{\partial v}{\partial \mathbf{n}} = 0$ are imposed. Such Neumann-like boundary conditions are formulated as algebraic constraints of the form $\phi_{\mathrm{b}} - \phi_{\mathrm{in}} = 0$, where $\phi_{\mathrm{in}}$ denotes the field value at the interior node adjacent to the physical boundary along the normal direction.

The nonlinear algebraic system given in \cref{eq:staggered_algebra_ns} is solved via Newton iteration with analytical Jacobian, following the procedure described in the previous subsection. The Jacobian matrix $G'(\mathbf{U})$ takes the form
\begin{equation}\label{jacobian_ns}
\left[\begin{array}{lll}
    \mathbf{M}_{1}(u, v) & \mathrm{diag}(\mathbf{M}_{\mathrm{cy}}u) & \mathbf{M}_{\mathrm{gx}} \\
    \mathrm{diag}(\mathbf{M}_{\mathrm{cx}}v) & \mathbf{M}_{2}(u, v) & \mathbf{M}_{\mathrm{gy}} \\
    \mathbf{M}_{\mathrm{dx}} & \mathbf{M}_{\mathrm{dy}} & \mathbf{M}_{\mathrm{p}}
\end{array}\right],
\end{equation}
where $\mathbf{M}_{1}(u, v) = \mathrm{diag}(\mathbf{M}_{\mathrm{cx}}u) + \mathrm{diag}(u)\mathbf{M}_{\mathrm{cx}} + \mathrm{diag}(v)\mathbf{M}_{\mathrm{cy}} + \mathbf{M}_{\mathrm{vu}}$, and 
$\mathbf{M}_{2}(u, v)=\mathrm{diag}(u)\mathbf{M}_{\mathrm{cx}}+\mathrm{diag}(\mathbf{M}_{\mathrm{cy}}v) + \mathrm{diag}(v)\mathbf{M}_{\mathrm{cy}} + \mathbf{M}_{\mathrm{vv}}$. The linear sub-problem arising at each Newton step is solved by the generalized minimal residual method on CPU via {\tt scipy.sparse.linalg.gmres} or on GPU via {\tt cupyx.scipy.sparse.linalg.gmres}, for higher efficiency than the direct solver adopted previously. Finally, the obtained staggered-mesh solutions are mapped onto the vertex nodes via linear interpolation, as illustrated below:
\begin{equation}\label{eq:staggered_correct}
    u_{i, j} = \frac{u_{i-\frac{1}{2}, j} + u_{i+\frac{1}{2}, j}}{2}, \quad
    v_{i, j} = \frac{v_{i, j-\frac{1}{2}} + v_{i, j+\frac{1}{2}}}{2}.
\end{equation}

\subsection{Vertex-centered FVM for compressible Euler equations}\label{appendix_a3}
Conventional cell-centered FVM stores variables at cell centroids. To align with the neural network formulation, we adopt the vertex-centered FV formulation built upon dual mesh and dual-cell construction. Dual cells are formed by connecting centroids of cells sharing a common grid point, as illustrated in \cref{fig:mesh_discretization}(d). For each dual cell, we approximately calculate its area by
\begin{equation}\label{eq:dual_area}
    A_{i} = \frac{1}{4} \sum_{j \in \mathcal{N}_{i}} \tilde{A}_{j} \, ,
\end{equation}
where $\mathcal{N}_{i}$ is the index set of primal cells surrounding grid point $i$, and $\tilde{A}_{j}$ denotes the area of the $j$-th primal cell. This expression is exact for orthogonal meshes and reduces to an area-averaging approximation under non-orthogonal grid conditions.

We employ pseudo-time marching for iterative solution. The steady-state system given by \cref{eq:euler_formulation} is recast into its pseudo-time-dependent form as
\begin{equation}\label{eq:euler_formulation_pseudotime}
\begin{array}{lll}
    \frac{\partial \rho}{\partial \tau} + \nabla \cdot (\rho \mathbf{u}) = 0  \\[2pt]
    \frac{\partial \rho \mathbf{u}}{\partial \tau} + \nabla \cdot (\rho \mathbf{u}\mathbf{u}) + \nabla p = \mathbf{0} \\[2pt]
    \frac{\partial \rho E}{\partial \tau} + \nabla \cdot (\rho E \mathbf{u}) + \nabla \cdot (p \mathbf{u}) = 0
\end{array}.
\end{equation}
Based on this system, we work with conservative variables and define the solution vector $\mathbf{U} = [\rho, \, \rho u, \, \rho v, \, \rho E]^{\top}$. The equations are written in the compact vector form $\frac{\partial \mathbf{U}}{\partial\tau} + \mathcal{R}(\mathbf{U}) = \boldsymbol{0}$, where $\mathcal{R}(\mathbf{U})$ denotes the spatial operator collecting all inviscid flux terms. 

In the present hypersonic flow setting, dimensional quantities are non-dimensionalised using free-stream reference conditions. The free-stream speed of sound is $c_{\infty} = \sqrt{\gamma R T_{\infty}}$, where $T_{\infty}$ is the free-stream temperature. Taking $l$ as the geometric characteristic length, the key non-dimensional scaling relations read
\begin{equation}\label{eq:nondim_scaling}
\mathbf{x}^{*} = \frac{\mathbf{x}}{l}, \quad \nabla^{*} = l\,\nabla, \quad \tau^{*} = \frac{c_{\infty} \tau}{l}, \quad \rho^{*} = \frac{\rho}{\rho_{\infty}}, \quad \mathbf{u}^{*} = \frac{\mathbf{u}}{c_{\infty}}, \quad p^{*} = \frac{p}{\rho_{\infty} c_{\infty}^{2}}, \quad T^{*} = \frac{T}{T_{\infty}}, \quad E^{*} = \frac{E}{c_{\infty}^{2}}.
\end{equation}
Under this scaling, all free-stream reference quantities reduce to unity, and the non-dimensional speed of sound reads $c^{*}=\sqrt{T^{*}}$. The free-stream Mach number is defined as $\mathrm{Ma}_{\infty}=|\mathbf{u}_{\infty}^{*}|$. From the ideal-gas law one obtains $p^{*}=\frac{\rho^{*}T^{*}}{\gamma}$, leading to $p_{\infty}^{*} = \frac{1}{\gamma}$. The non-dimensional free-stream total energy is $E_{\infty}^{*}= \frac{1}{\gamma(\gamma-1)}+\frac{1}{2}\mathrm{Ma}_{\infty}^{2}$.
We take the cylinder radius as the characteristic length $l$, with $l=1$. The non-dimensional equations marked by $(\cdot)^{*}$ are equivalent to the original pseudo-time Euler system \cref{eq:euler_formulation_pseudotime}. For brevity, the asterisk superscript is omitted hereafter, and all variables are to be interpreted as non-dimensional.

For pseudo-time discretization, we employ the first-order forward Euler time-stepping scheme, given by $\mathbf{U}_{n+1}-\mathbf{U}_{n} = -\mathcal{R}(\mathbf{U}_{n}) \Delta \tau$, and $\mathbf{U}_{0}$ is initialized with the free-stream state. The pseudo-time step size of each dual cell is determined by
\begin{equation}\label{eq:pseudotime_step}
{\Delta \tau}_{i} = \mathrm{CFL} \frac{A_{i}}{\Lambda_i},
\end{equation}
where the CFL number ranges from $0$ to $1$ for explicit time-stepping. $A_{i}$ is the area of the dual control-volume for node $i$, and $\Lambda_{i}$ is the local spectral radius aggregated from edge-wise wave speeds, given by
\begin{equation}\label{eq:spectral_radius}
    \Lambda_{i} = \sum_{j \in \mathcal{N}_{i}} \big( |\mathbf{u}_{ij} \cdot \mathbf{n}_{ij}| + c_{ij}\big) d_{ij},
\end{equation}
where $\mathbf{u}_{ij}$, $c_{ij}$ and $\mathbf{n}_{ij}$ are respectively the velocity, sound speed and unit outward normal, all evaluated at the dual-edge midpoint, and $d_{ij}$ denotes the dual edge length.

The two-dimensional FV formulation converts the differential governing equations into integral form to evaluate the flux residual $\mathcal{R}(\mathbf{U})$ via the Gauss-Green theorem, given by
\begin{equation}\label{eq:gauss_green}
    \iint_\Omega \nabla \cdot \phi \, \mathrm{d}S = \oint_{\partial \Omega} \phi \mathbf{n} \, \mathrm{d}l.
\end{equation}
For each physical field $\phi$ on the left and right face states of the dual edge, we adopt the MUSCL reconstruction, given by
\begin{equation}\label{eq:muscl_reconstruction}
    \phi_{L} = \phi_{i} + \psi_{i} \nabla \phi_{i} \cdot \mathbf{r}_{i}, \quad \phi_{R} = \phi_{j} + \psi_{j} \nabla \phi_{j} \cdot \mathbf{r}_j,
\end{equation}
where $\mathbf{r}_{i}$ and $\mathbf{r}_{j}$ denote the position vectors from node $i$ and node $j$ to the dual-edge midpoint, respectively. The gradient $\nabla \phi$ is computed via an eight-point least-squares procedure, and $\psi$ is defined by the Venkatakrishnan limiter, given by
\begin{equation}\label{eq:venkata_limiter}
\psi(\theta) = \max \left( \frac{\theta^{2} + 2\theta}{\theta^{2} + \theta+2}, \, 0\right),
\end{equation}
where $\theta$ denotes the local slope ratio evaluated along each dual edge. Edge-wise limiting coefficients $\psi(\theta)$ are aggregated onto each vertex node by taking the minimum value over all connected edges to obtain the nodal limiter value $\psi$ used in the reconstruction above. The reconstructed left and right face states are employed to compute the numerical flux at the dual edge using the Rusanov scheme, which consists of a central difference term plus a maximum wave speed dissipation term,
\begin{equation}\label{eq:rusanov_face_flux}
\mathcal{F}_{ij}\big(\phi_{L},\phi_{R}\big) = \frac{1}{2}\big(\mathcal{F}(\phi_{L}) + \mathcal{F}(\phi_{R})\big)-\frac{1}{2} \max(\lambda_{i}, \lambda_{j}) \big(\phi_{R}-\phi_{L}\big)d_{ij},
\end{equation}
where $\lambda_{i}$ and $\lambda_{j}$ denote the local wave speeds $|\mathbf{u} \cdot \mathbf{n}| + c$ evaluated from the left and right reconstructed states on this dual edge, respectively. Subsequently, these edge-face fluxes are aggregated onto vertex nodes to obtain the volumetric residual contribution
\begin{equation}\label{eq:rusanov_volumetric_flux}
\mathcal{R}_{i} = \sum_{j \in \mathcal{N}_i} \mathrm{sgn}(i, j)\, \mathcal{F}_{ij}
\end{equation}
where $\mathrm{sgn}(i,j) = 1$ for $i<j$ and otherwise $-1$. $\mathcal{F}_{ij}$ denotes the numerical flux evaluated over the directed dual edge from index $i$ to $j$. We define the directed edge $e_{ij}$ to point from the smaller index toward the larger index, with positive flux corresponding to the outward-normal direction of this directed edge.

In this vertex-centered scheme, physical boundaries coincide with vertex locations. At the inlet boundary, free-stream conditions are given in primitive-variable form and Dirichlet-type constraints are applied to boundary-vertex conservative unknowns. Analogous to the auxiliary Neumann-type extrapolation boundaries widely adopted for completeness in incompressible flow computations, zero-order constant extrapolation is employed for the outlet boundary, which reads $\phi_{\mathrm{b}} = \phi_{\mathrm{in}}$, where $\phi_{\mathrm{in}}$ denotes the primitive variable value at the nearest interior vertex along the boundary nodal normal direction. For inviscid solid-wall boundaries, we impose the slip no-penetration condition $\mathbf{u} \cdot \mathbf{n} = 0$, together with auxiliary zero-normal-gradient extrapolation assumptions $\dfrac{\partial T}{\partial \mathbf{n}} = 0$, $\dfrac{\partial \mathbf{u}_{\tau}}{\partial \mathbf{n}} = 0$ and $\dfrac{\partial p}{\partial \mathbf{n}}=0$. Combining these constraints and the ideal-gas equation of state yields the primitive states at boundary nodes:
\begin{equation}\label{eq:euler_surface_constraints}
\rho_{\mathrm{b}} = \rho_{\mathrm{in}}, \quad
\mathbf{u}_{\mathrm{b}} = \mathbf{u}_{\mathrm{in}} - \big(\mathbf{u}_{\mathrm{in}} \cdot \mathbf{n}_{\mathrm{b}}\big) \mathbf{n}_{\mathrm{b}}, \quad
p_{\mathrm{b}} = p_{\mathrm{in}}.
\end{equation}
These boundary constraints expressed in primitive quantities are further transformed into corresponding conservative-variable boundary constraints, and the resulting conservative-variable constraints are enforced during the iterative solution procedure.
Due to the hypersonic flow characteristics and the MUSCL reconstruction procedure, the present compressible Euler case cannot achieve the same clean convergence behaviour observed for the incompressible test cases, and minor oscillations persist in the residual history. Accordingly, we adopt the wall-surface temperature distribution as the convergence check metric. The pseudo-time iteration is terminated when the wall-temperature field ceases to exhibit notable variations.

\section{BC enforcement and loss functions}\label{appendix_b}
All test cases in this work adopt unsupervised training without labeled observational data. The network optimization is constrained by a composite loss function consisting of PDE residual and boundary condition terms:
\begin{equation}\label{total_loss_term}
    \mathcal{L} = \omega_{1} \mathcal{L}_{\mathrm{PDE}} + \omega_{2} \mathcal{L}_{\mathrm{BC}},
\end{equation}
where both terms are formulated via MSE evaluation. For fully hard-constrained boundary implementations, $\mathcal{L}_{\mathrm{BC}} \equiv 0$. For partially hard-constrained BC cases targeting a priori-free engineering usage, equal weights $\omega_1 = \omega_2 = 1$ are assigned. Within each individual loss term, sub-components are treated with equal weighting, except for hypersonic test cases, where additional weighting factors are introduced to account for large magnitude discrepancies. Further details are discussed separately in the following subsections.

\subsection{Poisson-type equations}\label{appendix_b1}
These Poisson-type test cases feature a single unknown defined over convex domains, so global hard boundary constraints can be readily realized. We directly employ a distance function based ansatz analogous to \cref{eq:hard_bc_formulation}. Let $\tilde{u}$ be the network output, the prediction $\hat{u}$ for the Sine case reads
\begin{equation}\label{eq:sine_hard_constraints}
    \hat{u} = xy(1 - x)(1 - y) \tilde{u},
\end{equation}
for the Polynomial case
\begin{equation}\label{eq:polynomial_hard_constraints}
    \hat{u} = (x + \frac{\pi}{6})(x - \frac{\pi}{6})(y + \frac{\pi}{6})(y - \frac{\pi}{6}) \big( \tilde{u} - \tan (x + y) \big) + \tan (x + y),
\end{equation}
and for the Liouville case
\begin{equation}\label{eq:liouville_hard_constraints}
    \hat{u} = (x^{2} + y^{2} - 1) \tilde{u}.
\end{equation}

Network prediction and hard-constraint treatments are shared across all baseline models. Discretization-based neural models differ from AD-based PINN\textsubscript{c} in the formulation of PDE residual constraints, which are enforced via the discrete algebraic residual given by \cref{eq:residual_fun} at interior nodes.

\subsection{Lid-driven cavity flow}\label{appendix_b2}
Unlike the staggered-grid discretization scheme in Appendix~\ref{appendix_a2}, collocated grids are adopted for the neural discrete models, requiring only three discrete operator matrices $\mathbf{M}_{\mathrm{gx}}$, $\mathbf{M}_{\mathrm{gy}}$, and $\mathbf{M}_{\mathrm{v}}$, corresponding to the x-direction gradient operator, y-direction gradient operator, and viscous-term operator, respectively. Since boundary conditions are not embedded within these discrete matrices, a single viscous-term operator $\mathbf{M}_{\mathrm{v}}$ can be shared for both x and y directions. The resulting algebraic system reads
\begin{equation}\label{eq:collocated_algebra_ns}
    \left[\begin{array}{lll}
        \mathrm{diag}(u) \mathbf{M}_{\mathrm{gx}} + \mathbf{M}_{\mathrm{v}} & \mathrm{diag}(v) \mathbf{M}_{\mathrm{gy}} & \mathbf{M}_{\mathrm{gx}} \\
        \mathrm{diag}(u) \mathbf{M}_{\mathrm{gx}} & \mathrm{diag}(v) \mathbf{M}_{\mathrm{gy}} + \mathbf{M}_{\mathrm{v}} & \mathbf{M}_{\mathrm{gy}} \\
        \mathbf{M}_{\mathrm{gx}} & \mathbf{M}_{\mathrm{gy}} & \mathbf{O}
    \end{array}\right]
    \begin{bmatrix}
        u\\v\\p
    \end{bmatrix}
    =\mathbf{0}.
\end{equation}

For each PDE component of the NS system, equal weighting is adopted for residual constraints, i.e.,
\begin{equation}\label{eq:ns_pde_loss}
    \mathcal{L}_{\mathrm{PDE}} = \mathcal{L}_{\mathrm{mom}_{x}} + \mathcal{L}_{\mathrm{mom}_{y}} + \mathcal{L}_{\mathrm{div}},
\end{equation}
where $\mathcal{L}_{\mathrm{mom}_{x}}$, $\mathcal{L}_{\mathrm{mom}_{y}}$, and $\mathcal{L}_{\mathrm{div}}$ denote the MSE losses evaluated from the discrete residual of x-momentum, y-momentum, and divergence free constraint, respectively, to maintain consistent constraint intensity for each governing equation.
The coupled NS system involves distinct boundary conditions for different solution fields, thus separate constraint treatments are required for each variable. Let $\tilde{\mathbf{U}} = [\tilde{u}, \, \tilde{v}, \, \tilde{p}]$ be the network output, and the trial solution for each flow variable is constructed as
\begin{equation}\label{eq:lid_driven_hard_constraints}
    \hat{u} = \psi(\Omega) \tilde{u} + \big( 1 - \psi(\Omega) \big) 4x(1 - x)y, \quad \hat{v} = \psi(\Omega) \tilde{v}, \quad \hat{p} = (x^{2} + y^{2}),
\end{equation}
where $\psi(\Omega) = xy(1 - x)(1 - y)$ denotes the boundary-zeroing factor. This formulation enforces global hard boundary constraints, such that $\mathcal{L}_{\mathrm{BC}} \equiv 0$ holds naturally for \cref{eq:lid_driven_bc}.

\subsection{Backward-facing step flow}\label{appendix_b3}
The backward-facing step case adopts the same PDE constraint formulation as described in the lid-driven cavity section. Nevertheless, global hard boundary constraints cannot be realized for this configuration, so hard constraints are enforced only on partial boundaries. The boundary-condition residual reads
\begin{equation}\label{eq:bward_step_loss_bc}
    \mathcal{L}_{\mathrm{BC}} = \sum_{\Omega_{k} \in \partial \Omega} \mathcal{L}_{\Omega_{k}}(\hat{\mathbf{u}}) + \sum_{\Omega_{k} \in \partial \Omega} \mathcal{L}_{\Omega_{k}}(\hat{p}),
\end{equation}
with equal contributions from each boundary component. Note that terms over hard-constrained boundary segments vanish identically, $\mathcal{L}_{\Omega_{k}}(\cdot)\equiv 0$, which permits a unified loss formulation for both hard and soft boundary treatments. The hard-constrained ansatz is given by
\begin{equation}\label{eq:bward_step_hard_constraints}
    \hat{u} = y(2 - y) \tilde{u}, \quad \hat{v} = xy(2 - y) \tilde{v}, \quad \hat{p} = (4 - y) \tilde{p}.
\end{equation}
Under this construction, hard constraints for $u$ are satisfied at the top and bottom step walls; hard constraints for $v$ are imposed at the step walls as well as the inlet boundary. For pressure, only the outlet boundary condition needs to be satisfied, which is fully achieved.

For the discrete neural models PINN\textsubscript{d} and PIGNN, the DNS-compatible extrapolation constraints can either be explicitly enforced or omitted. When discrete operator stencils are applied, neighbouring boundary nodes are inherently coupled by the neural discrete constraints. We have tested the outlet boundary conditions $\frac{\partial u}{\partial \mathbf{n}} = 0$ and $\frac{\partial v}{\partial \mathbf{n}} = 0$, which are penalized in the form $\frac{1}{|\Omega_{\mathrm{b}}|}|\hat{\phi}_{\mathrm{b}} - \hat{\phi}_{\mathrm{in}}|^{2}$. Observable discrepancies arise between the two treatments, as shown in \cref{tab:bward_step_outlet_neumann}. Omitting these outlet constraints yields slightly lower velocity errors yet larger pressure errors, whereas enforcing them produces the opposite trend. Such discrepancies exert a milder influence on PIGNN, while PINN\textsubscript{d} is more sensitive, particularly for velocity quantities near the step corner. Considering that pressure gradient related quantities are of greater practical interest than absolute pressure values for incompressible engineering problems, we omit these extrapolation-type outlet constraints in this test case.

\begin{table}[!ht]
    \centering
    \caption{Global relative $L_{2}$ errors for the backward-facing step flow at $\mathrm{Re}=100$, for discrete models with different outlet boundary treatments and various mesh resolutions.}
    \label{tab:bward_step_outlet_neumann}
    \small
    \begin{tabularx}{0.95\textwidth}{@{\hspace{0.5em}}lllXXX@{\hspace{0.5em}}}
        \toprule
        \multirow{2}{*}{Outlet BC} & \multirow{2}{*}{Method} & \multirow{2}{*}{Metric} & \multicolumn{3}{c}{Resolution} \\
        \cmidrule{4-6}
        & & & $100 \times 50$ & $200 \times 100$ & $300 \times 150$ \\
        \midrule
        \multirow{6}{*}{$p = 0$}
        & \multirow{3}{*}{PINN\textsubscript{d}}
        & $u$ & $1.82\mathrm{e}{-2} \pm 6.19\mathrm{e}{-3}$ & $1.87\mathrm{e}{-2} \pm 5.18\mathrm{e}{-3}$ & $1.69\mathrm{e}{-2} \pm 5.89\mathrm{e}{-3}$ \\
        & & $v$ & $0.1473 \pm 4.60\mathrm{e}{-2}$ & $0.1497 \pm 3.26\mathrm{e}{-2}$ & $0.1445 \pm 3.33\mathrm{e}{-2}$ \\
        & & $p$ & $0.1245 \pm 2.26\mathrm{e}{-2}$ & $0.1545 \pm 9.05\mathrm{e}{-3}$ & $0.1518 \pm 1.30\mathrm{e}{-2}$ \\
        \cmidrule(lr){2-6}
        & \multirow{3}{*}{PIGNN}
        & $u$ & $6.13\mathrm{e}{-3} \pm 1.28\mathrm{e}{-3}$ & $5.19\mathrm{e}{-3} \pm 1.31\mathrm{e}{-3}$ & $7.74\mathrm{e}{-3} \pm 1.86\mathrm{e}{-3}$ \\
        & & $v$ & $8.47\mathrm{e}{-2} \pm 2.32\mathrm{e}{-2}$ & $7.29\mathrm{e}{-2} \pm 5.67\mathrm{e}{-3}$ & $0.1011 \pm 1.66\mathrm{e}{-2}$ \\
        & & $p$ & $0.1159 \pm 1.10\mathrm{e}{-2}$ & $0.1427 \pm 1.67\mathrm{e}{-2}$ & $0.1221 \pm 8.29\mathrm{e}{-3}$ \\
        \midrule
        \multirow{6}{*}{\shortstack{$p = 0$ \\ $\frac{\partial u}{\partial \mathbf{n}} = 0$ \\ $\frac{\partial v}{\partial \mathbf{n}} = 0$}}
        & \multirow{3}{*}{PINN\textsubscript{d}}
        & $u$ & $1.99\mathrm{e}{-2} \pm 4.61\mathrm{e}{-3}$ & $2.11\mathrm{e}{-2} \pm 5.92\mathrm{e}{-3}$ & $1.78\mathrm{e}{-2} \pm 6.20\mathrm{e}{-3}$ \\
        & & $v$ & $0.1710 \pm 2.81\mathrm{e}{-2}$ & $0.1651 \pm 3.68\mathrm{e}{-2}$ & $0.1458 \pm 3.09\mathrm{e}{-2}$ \\
        & & $p$ & $0.1075 \pm 1.21\mathrm{e}{-2}$ & $0.1340 \pm 8.83\mathrm{e}{-3}$ & $0.1500 \pm 1.53\mathrm{e}{-2}$ \\
        \cmidrule(lr){2-6}
        & \multirow{3}{*}{PIGNN}
        & $u$ & $6.48\mathrm{e}{-3} \pm 9.54\mathrm{e}{-4}$ & $4.61\mathrm{e}{-3} \pm 5.08\mathrm{e}{-4}$ & $9.61\mathrm{e}{-3} \pm 4.14\mathrm{e}{-3}$ \\
        & & $v$ & $9.09\mathrm{e}{-2} \pm 8.05\mathrm{e}{-3}$ & $7.64\mathrm{e}{-2} \pm 1.22\mathrm{e}{-2}$ & $0.1063 \pm 2.04\mathrm{e}{-2}$ \\
        & & $p$ & $7.72\mathrm{e}{-2} \pm 1.07\mathrm{e}{-2}$ & $0.1214 \pm 9.64\mathrm{e}{-3}$ & $0.1142 \pm 1.97\mathrm{e}{-2}$ \\
        \bottomrule
    \end{tabularx}
\end{table}

\subsection{Hypersonic inviscid cylinder flow}\label{appendix_b4}
In this case, neural network training directly imposes constraints for the steady Euler equations. Unlike pseudotime marching numerical solvers, conservative variables are therefore not required here, and the network predicts primitive variables directly, namely $\tilde{\mathbf{U}} = [\tilde{\rho}, \, \tilde{u}, \, \tilde{v}, \, \tilde{p}]$. In contrast to pressure in incompressible formulations, density $\rho$ and pressure $p$ must obey thermodynamic constraints, where negative values are non-physical. We therefore apply the softplus function $\vartheta$ to enforce positivity, given by
\begin{equation}\label{eq:softplus_fun}
    \vartheta(x) = \log \big(1 + \exp(x)\big).
\end{equation}
Prior to hard constraint enforcement, the network outputs for density and pressure are transformed as $\vartheta(\hat{\rho}) + \epsilon$ and $\vartheta(\hat{p}) + \epsilon$.

For the PDE residual, PINN\textsubscript{c} employs constraints in differential form, while the discrete models adopt the integral form Rusanov face flux, as illustrated in \cref{eq:rusanov_face_flux}, and aggregate to volumetric flux following \cref{eq:rusanov_volumetric_flux} consistent with the DNS implementation. Based on the nondimensional formulation in \cref{eq:nondim_scaling} and the unsteady system in \cref{eq:euler_formulation_pseudotime}, the momentum and energy flux terms scale as $\mathrm{Ma}$ and $\mathrm{Ma}^{2}$ relative to the mass flux, respectively. Magnitude discrepancies among different equation components grow with increasing Mach number. Accordingly, fixed weights are applied to each PDE residual component according to this scaling consideration, which reads
\begin{equation}\label{eq:bowshock_loss_pde}
    \mathcal{L}_{\mathrm{PDE}} = \mathcal{L}_{\mathrm{mass}} + \frac{1}{\mathrm{Ma}^{2}} \mathcal{L}_{\mathrm{mom}_{x}} + \frac{1}{\mathrm{Ma}^{2}} \mathcal{L}_{\mathrm{mom}_{y}} + \frac{1}{\mathrm{Ma}^{4}} \mathcal{L}_{\mathrm{enrg}}.
\end{equation}

The boundary condition residual is treated analogously to incompressible cases, which reads
\begin{equation}\label{eq:bowshock_loss_bc}
    \mathcal{L}_{\mathrm{BC}} = \sum_{\Omega_{k} \in \partial \Omega} \mathcal{L}_{\Omega_{k}}(\hat{\rho}) + \sum_{\Omega_{k} \in \partial \Omega} \mathcal{L}_{\Omega_{k}}(\hat{\mathbf{u}}) + \sum_{\Omega_{k} \in \partial \Omega} \mathcal{L}_{\Omega_{k}}(\hat{p}).
\end{equation}
Since the coordinate range leads to poorly scaled multipliers for inlet hard constraints, we normalize the location input as $\mathbf{x}' = \frac{\mathbf{x}}{\max(\mathbf{x})}$ for stable training. The hard Dirichlet enforcement on the inlet is given by
\begin{equation}\label{eq:bowshock_hard_constraints}
    \hat{\phi} = \psi(\Omega_{\mathrm{inlet}})\,\tilde{\phi} + \big(1 - \psi(\Omega_{\mathrm{inlet}})\big)\,\phi_{\infty},
\end{equation}
where $\psi(\Omega_{\mathrm{inlet}}) = \frac{3.4^{2}}{9} - \big(x - \frac{1.6}{3}\big)^2 - y^2$, the normalization factor satisfies $\max(\mathbf{x}) = 3$, and $\phi_{\infty}$ denotes the respective free-stream physical values $[1, \, \mathrm{Ma}, \, 0, \, \frac{1}{\gamma}]$ for $\rho$, $u$, $v$, $p$. At the outlet, Dirichlet type pressure constraints as employed for compressible inflow boundaries are unavailable, since all flow quantities remain unknown. The PINN\textsubscript{c} is defined over a continuous space domain, and nodal level coupling interactions present in discrete operators do not exist. We therefore incorporate outlet points into the PDE residual $\mathcal{L}_{\mathrm{PDE}}$. For solid wall boundaries, the slip condition $\hat{\mathbf{u}} \cdot \mathbf{n} = 0$ is enforced via soft boundary condition penalty terms. By contrast, discrete neural solvers achieve indirect outlet constraints through interior point PDE evaluations that reference boundary node values. Nevertheless, we observe faster convergence when additionally imposing the zero-normal-gradient outlet condition $\frac{\partial \phi}{\partial \mathbf{n}} = 0$. On solid wall surfaces, in addition to the slip velocity constraint, further zero-normal-gradient conditions $\frac{\partial T}{\partial \mathbf{n}} = 0$ and $\frac{\partial p}{\partial \mathbf{n}} = 0$ are also applied.

\section{Soft BC enforcement results} \label{appendix_c}
We conduct experiments with full soft boundary condition penalties for the fully-hard cases, Poisson-type cases and lid-driven cavity case to form a comparison, strictly following the same training configurations given in \cref{tab:training_params}. The experimental results for Sine, Polynomial, Liouville, and lid-driven cavity are reported in \cref{tab:sine_l2_softbc}, \cref{tab:polynomial_l2_softbc}, \cref{tab:liouville_l2_softbc}, and \cref{tab:lid_driven_l2_softbc}, respectively.
It can be observed that under identical training settings, the errors obtained with soft boundary condition penalties are evidently larger than those of the pure hard-constraint experiments presented in \cref{sec:poisson} and \cref{sec:lid_driven}. This effect is particularly prominent for PIGNN: it exhibits higher parameter sensitivity and tends to diverge under high-resolution, large-scale input conditions, which can be mitigated by adjusting network architectures, tuning training hyperparameters, or extending the training process. Across all test cases, PIGNN yields substantially larger errors compared with MLP-based PINN\textsubscript{c} and PINN\textsubscript{d}. Counter-intuitively, this performance gap diminishes as the degree of nonlinearity increases. The discrepancy can exceed one order of magnitude for the linear Sine case, whereas it becomes considerably smaller for the more strongly coupled NS equations. Moreover, the divergence frequently observed in simple benchmark cases no longer occurs. This observation further supports our earlier argument: the strength of PIGNN lies in its superior capability for capturing nonlinear and cross-coupling features, while such advantages are severely restricted by the over-squashing effect in simple, weakly nonlinear problems.
For MLP-based models, AD and high-order discrete operators bring no obvious performance gain on the simple Sine and Polynomial benchmarks. This indicates mutual interference among multiple optimization objectives. Under pure hard constraints, such interference is substantially suppressed and facilitates stable optimization.

\begin{table}[!ht]
    \centering
    \caption{Global relative $L_{2}$ errors for the linear Sine test case with soft BC enforcement, across different mesh resolutions.}
    \label{tab:sine_l2_softbc}
    \small
    \begin{threeparttable}
    \begin{tabularx}{0.85\textwidth}{@{\hspace{0.5em}}lXXXX@{\hspace{0.5em}}}
        \toprule
        \multirow{2}{*}{Method} & \multicolumn{4}{c}{Resolution} \\
        \cmidrule{2-5}
        & $20 \times 20$ & $50 \times 50$ & $100 \times 100$ & $200 \times 200$ \\
        \midrule
        PINN\textsubscript{c} & $7.46\mathrm{e}{-3} \pm 4.15\mathrm{e}{-3}$ & $5.83\mathrm{e}{-3} \pm 1.23\mathrm{e}{-3}$ & $5.37\mathrm{e}{-3} \pm 2.25\mathrm{e}{-3}$ & $4.79\mathrm{e}{-3} \pm 2.22\mathrm{e}{-3}$ \\
        PINN\textsubscript{d} ($2^{nd}$) & $5.17\mathrm{e}{-3} \pm 2.71\mathrm{e}{-3}$ & $4.55\mathrm{e}{-3} \pm 3.35\mathrm{e}{-3}$ & $4.53\mathrm{e}{-3} \pm 1.81\mathrm{e}{-3}$ & $6.51\mathrm{e}{-3} \pm 3.48\mathrm{e}{-3}$ \\
        PINN\textsubscript{d} ($4^{th}$) & $5.63\mathrm{e}{-3} \pm 2.75\mathrm{e}{-3}$ & $5.59\mathrm{e}{-3} \pm 1.58\mathrm{e}{-3}$ & $5.69\mathrm{e}{-3} \pm 3.22\mathrm{e}{-3}$ & $3.93\mathrm{e}{-3} \pm 2.47\mathrm{e}{-3}$ \\
        PIGNN ($2^{nd}$) & $8.98\mathrm{e}{-2} \pm 2.71\mathrm{e}{-2}$ & $6.72\mathrm{e}{-2} \pm 1.53\mathrm{e}{-2}$ & / & / \\
        PIGNN ($4^{th}$) & $9.91\mathrm{e}{-2} \pm 2.26\mathrm{e}{-2}$ & $8.70\mathrm{e}{-2} \pm 6.38\mathrm{e}{-2}$ & / & / \\
        \bottomrule
    \end{tabularx}
    \begin{tablenotes}
        \footnotesize 
        \item Note: "/" denotes configurations where optimization instability occurs frequently, such that no consistently converged results can be obtained across repeated trials.
    \end{tablenotes}
    \end{threeparttable}
\end{table}

\begin{table}[!ht]
    \centering
    \caption{Global relative $L_{2}$ errors for the nonlinear Polynomial test case with soft BC enforcement, across different mesh resolutions.}
    \label{tab:polynomial_l2_softbc}
    \small
    \begin{threeparttable}
    \begin{tabularx}{0.85\textwidth}{@{\hspace{0.5em}}lXXXX@{\hspace{0.5em}}}
        \toprule
        \multirow{2}{*}{Method} & \multicolumn{4}{c}{Resolution} \\
        \cmidrule{2-5}
        & $20 \times 20$ & $50 \times 50$ & $100 \times 100$ & $200 \times 200$ \\
        \midrule
        PINN\textsubscript{c} & $1.25\mathrm{e}{-2} \pm 1.33\mathrm{e}{-2}$ & $5.27\mathrm{e}{-2} \pm 4.75\mathrm{e}{-2}$ & $2.63\mathrm{e}{-2} \pm 2.03\mathrm{e}{-2}$ & $4.13\mathrm{e}{-2} \pm 2.54\mathrm{e}{-2}$ \\
        PINN\textsubscript{d} ($2^{nd}$) & $1.73\mathrm{e}{-2} \pm 1.34\mathrm{e}{-2}$ & $2.21\mathrm{e}{-2} \pm 2.30\mathrm{e}{-2}$ & $1.21\mathrm{e}{-2} \pm 1.29\mathrm{e}{-3}$ & $4.54\mathrm{e}{-2} \pm 2.60\mathrm{e}{-2}$ \\
        PINN\textsubscript{d} ($4^{th}$) & $4.07\mathrm{e}{-3} \pm 2.66\mathrm{e}{-4}$ & $2.10\mathrm{e}{-2} \pm 2.06\mathrm{e}{-2}$ & $1.80\mathrm{e}{-2} \pm 1.18\mathrm{e}{-2}$ & $2.61\mathrm{e}{-2} \pm 1.91\mathrm{e}{-2}$ \\
        PIGNN ($2^{nd}$) & $0.1254 \pm 1.42\mathrm{e}{-2}$ & $0.1723 \pm 5.14\mathrm{e}{-2}$ & / & / \\
        PIGNN ($4^{th}$) & $7.50\mathrm{e}{-2} \pm 1.46\mathrm{e}{-2}$ & $0.1204 \pm 5.16\mathrm{e}{-2}$ & / & / \\
        \bottomrule
    \end{tabularx}
    \begin{tablenotes}
        \footnotesize 
        \item Note: "/" denotes configurations where optimization instability occurs frequently, such that no consistently converged results can be obtained across repeated trials.
    \end{tablenotes}
    \end{threeparttable}
\end{table}

\begin{table}[!ht]
    \centering
    \caption{Global relative $L_{2}$ errors for the stiff Liouville test case with soft BC enforcement, across different mesh resolutions.}
    \label{tab:liouville_l2_softbc}
    \small
    \begin{tabularx}{0.85\textwidth}{@{\hspace{0.5em}}lXXXX@{\hspace{0.5em}}}
        \toprule
        \multirow{2}{*}{Method} & \multicolumn{4}{c}{Resolution} \\
        \cmidrule{2-5}
        & $N_{r}=20, \, N_{\theta}=12$ & $N_{r}=50, \, N_{\theta}=24$ & $N_{r}=100, \, N_{\theta}=36$ & $N_{r}=200, \, N_{\theta}=48$ \\
        \midrule
        PINN\textsubscript{c} & $5.71\mathrm{e}{-2} \pm 5.40\mathrm{e}{-3}$ & $0.1925 \pm 4.39\mathrm{e}{-2}$ & $0.1759 \pm 5.47\mathrm{e}{-2}$ & $0.1991 \pm 4.48\mathrm{e}{-2}$ \\
        PINN\textsubscript{d} ($2^{nd}$) & $5.28\mathrm{e}{-2} \pm 3.34\mathrm{e}{-3}$ & $0.1538 \pm 4.13\mathrm{e}{-2}$ & $0.1575 \pm 5.36\mathrm{e}{-2}$ & $0.2191 \pm 4.03\mathrm{e}{-3}$ \\
        PINN\textsubscript{d} ($4^{th}$) & $7.26\mathrm{e}{-2} \pm 9.44\mathrm{e}{-3}$ & $0.1425 \pm 5.13\mathrm{e}{-2}$ & $0.1490 \pm 5.54\mathrm{e}{-2}$ & $0.1333 \pm 4.41\mathrm{e}{-2}$ \\
        PIGNN ($2^{nd}$) & $0.2499 \pm 6.51\mathrm{e}{-2}$ & $0.2632 \pm 3.29\mathrm{e}{-3}$ & $0.2117 \pm 2.09\mathrm{e}{-2}$ & $0.2185 \pm 1.45\mathrm{e}{-2}$ \\
        PIGNN ($4^{th}$) & $0.2998 \pm 5.97\mathrm{e}{-2}$ & $0.2706 \pm 1.23\mathrm{e}{-2}$
        & $0.2225 \pm 1.67\mathrm{e}{-2}$ & $0.2024 \pm 2.92\mathrm{e}{-3}$ \\
        \bottomrule
    \end{tabularx}
\end{table}

\begin{table}[!ht]
    \centering
    \caption{Global relative $L_{2}$ errors for the lid-driven cavity flow at $\mathrm{Re} = 100$ with soft BC enforcement, across different mesh resolutions.}
    \label{tab:lid_driven_l2_softbc}
    \small
    \begin{tabularx}{0.85\textwidth}{@{\hspace{0.5em}}llXXXX@{\hspace{0.5em}}}
        \toprule
        \multirow{2}{*}{Method} & \multirow{2}{*}{Metric} & \multicolumn{4}{c}{Resolution} \\
        \cmidrule{3-6}
        & & $50 \times 50$ & $100 \times 100$ & $150 \times 150$ & $200 \times 200$ \\
        \midrule
        \multirow{3}{*}{PINN\textsubscript{c}}
        & $u$ & $2.22\mathrm{e}{-2} \pm 2.30\mathrm{e}{-3}$ & $2.30\mathrm{e}{-2} \pm 2.83\mathrm{e}{-3}$ & $2.39\mathrm{e}{-2} \pm 2.04\mathrm{e}{-3}$ & $2.62\mathrm{e}{-2} \pm 2.63\mathrm{e}{-3}$ \\
        & $v$ & $3.70\mathrm{e}{-2} \pm 2.33\mathrm{e}{-3}$ & $3.72\mathrm{e}{-2} \pm 3.71\mathrm{e}{-3}$ & $3.82\mathrm{e}{-2} \pm 2.60\mathrm{e}{-3}$ & $4.33\mathrm{e}{-2} \pm 5.21\mathrm{e}{-3}$ \\
        & $p$ & $0.1395 \pm 2.52\mathrm{e}{-3}$ & $9.96\mathrm{e}{-2} \pm 2.40\mathrm{e}{-3}$ & $9.26\mathrm{e}{-2} \pm 1.48\mathrm{e}{-3}$ & $8.99\mathrm{e}{-2} \pm 4.08\mathrm{e}{-3}$ \\
        \midrule
        \multirow{3}{*}{PINN\textsubscript{d}}
        & $u$ & $2.64\mathrm{e}{-2} \pm 1.82\mathrm{e}{-3}$ & $2.30\mathrm{e}{-2} \pm 1.97\mathrm{e}{-3}$ & $2.33\mathrm{e}{-2} \pm 2.10\mathrm{e}{-3}$ & $2.30\mathrm{e}{-2} \pm 2.10\mathrm{e}{-3}$ \\
        & $v$ & $5.05\mathrm{e}{-2} \pm 2.81\mathrm{e}{-3}$ & $3.93\mathrm{e}{-2} \pm 2.83\mathrm{e}{-3}$ & $3.75\mathrm{e}{-2} \pm 3.58\mathrm{e}{-3}$ & $3.61\mathrm{e}{-2} \pm 2.22\mathrm{e}{-3}$ \\
        & $p$ & $0.1505 \pm 7.74\mathrm{e}{-3}$ & $0.1000 \pm 4.10\mathrm{e}{-3}$ & $8.79\mathrm{e}{-2} \pm 3.41\mathrm{e}{-3}$ & $8.58\mathrm{e}{-2} \pm 1.30\mathrm{e}{-3}$ \\
        \midrule
        \multirow{3}{*}{PIGNN}
        & $u$ & $3.71\mathrm{e}{-2} \pm 7.86\mathrm{e}{-3}$ & $4.38\mathrm{e}{-2} \pm 1.32\mathrm{e}{-2}$ & $4.48\mathrm{e}{-2} \pm 1.68\mathrm{e}{-2}$ & $5.24\mathrm{e}{-2} \pm 1.37\mathrm{e}{-2}$ \\
        & $v$ & $6.46\mathrm{e}{-2} \pm 1.03\mathrm{e}{-2}$ & $6.63\mathrm{e}{-2} \pm 1.43\mathrm{e}{-2}$ & $6.32\mathrm{e}{-2} \pm 1.80\mathrm{e}{-2}$ & $7.44\mathrm{e}{-2} \pm 1.64\mathrm{e}{-2}$ \\
        & $p$ & $0.1542 \pm 1.01\mathrm{e}{-2}$ & $0.1213 \pm 9.19\mathrm{e}{-3}$ & $0.1106 \pm 1.70\mathrm{e}{-2}$ & $0.1194 \pm 1.36\mathrm{e}{-2}$ \\
        \bottomrule
    \end{tabularx}
\end{table}

It should be specially emphasized that these soft-penalty experiments do not represent the optimal achievable performance of soft boundary conditions. They are intended solely for fair side-by-side comparison under identical experimental settings. Although better results could be achieved by prolonged training, we focus on evaluating general model capacity in this work. Hence, no posterior parameter tuning guided by ground-truth solutions is performed, since such post-hoc tuning is of limited practical value for real-world engineering applications.

\section{Hyperparameter analysis} \label{appendix_d}
\begin{figure}[!ht]
    \centering
    \begin{minipage}{0.185\textwidth}
        \centering
        \includegraphics[width=1\textwidth]{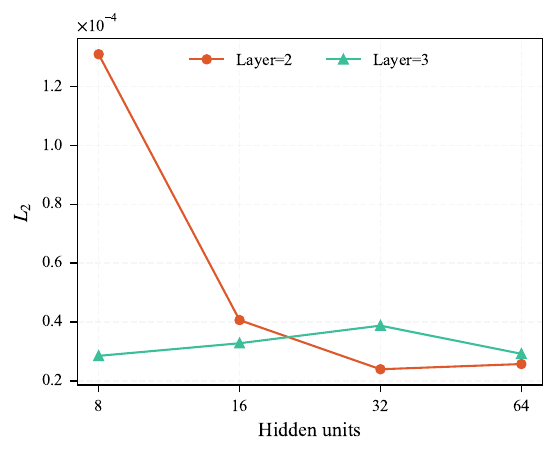}
        \caption*{(a)}
    \end{minipage}
    \begin{minipage}{0.185\textwidth}
        \centering
        \includegraphics[width=1\textwidth]{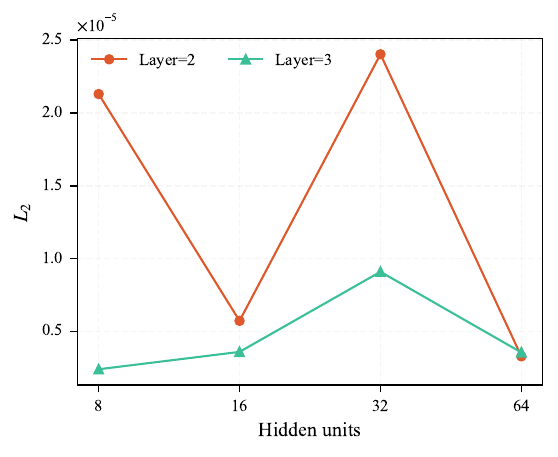}
        \caption*{(b)}
    \end{minipage}
    \begin{minipage}{0.195\textwidth}
        \centering
        \includegraphics[width=1\textwidth]{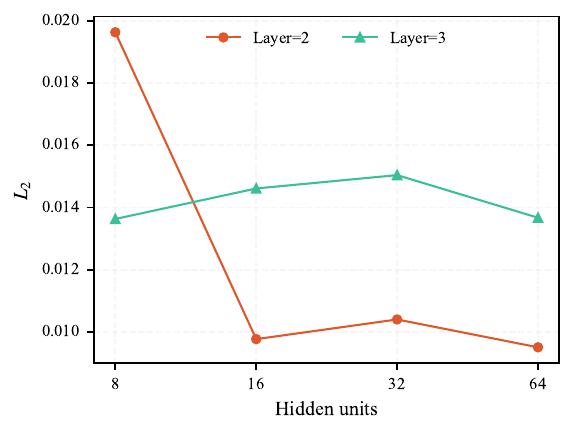}
        \caption*{(c)}
    \end{minipage}
    \begin{minipage}{0.195\textwidth}
        \centering
        \includegraphics[width=1\textwidth]{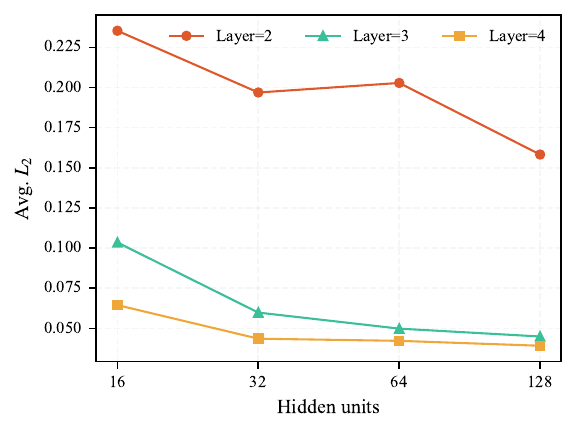}
        \caption*{(d)}
    \end{minipage}
    \begin{minipage}{0.195\textwidth}
        \centering
        \includegraphics[width=1\textwidth]{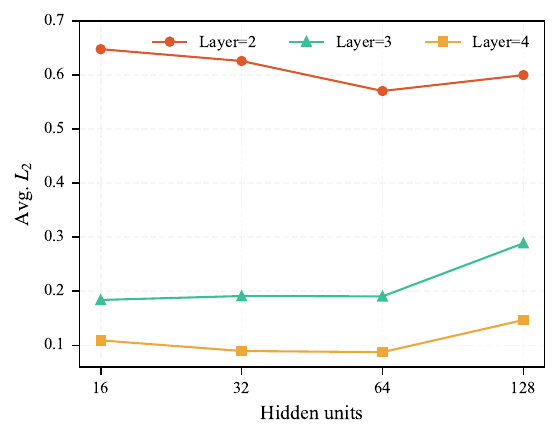}
        \caption*{(e)}
    \end{minipage}
    \caption{Hyperparameter analysis of hidden units and number of layers: (a)–(e) for Sine, Polynomial, Liouville, Lid-driven cavity, and backward-facing step.}
    \label{fig:hyperparams_hiddens}
\end{figure}

We carry out hyperparameter-tuning experiments for Poisson-type and incompressible NS cases, covering network width, network depth, and training duration. These experiments are restricted to the PINN\textsubscript{c} model, with each configuration executed only once. For Poisson-type cases, we test network depths of 2 and 3 layers, with hidden-layer units set to $[8, 16, 32, 64]$. For incompressible NS cases, depths of $[2, 3, 4]$ layers are examined with hidden-layer units $[16, 32, 64, 128]$. For this hyperparameter analysis, the training epoch is fixed to $20000$, and the corresponding results are shown in \cref{fig:hyperparams_hiddens}. Despite the randomness introduced by single-run evaluations, clear general trends can still be observed: model performance is substantially more sensitive to network depth than to network width. For comparable parameter budgets, deeper yet narrower architectures outperform shallower wider counterparts. Once the network size exceeds a certain threshold, further increases in depth or width yield marginal additional performance gains.

\begin{figure}[!ht]
    \centering
    \begin{minipage}{0.185\textwidth}
        \centering
        \includegraphics[width=1\textwidth]{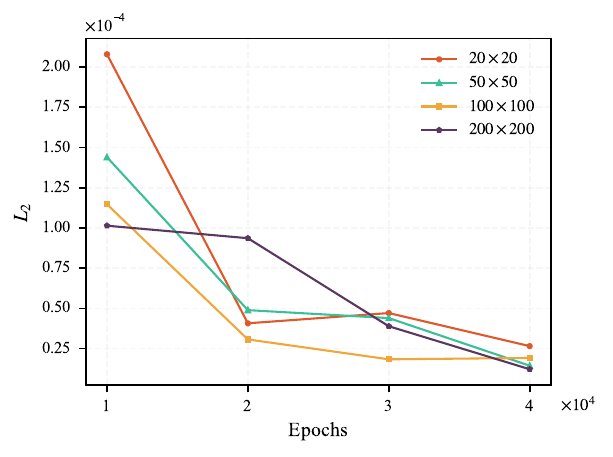}
        \caption*{(a)}
    \end{minipage}
    \begin{minipage}{0.185\textwidth}
        \centering
        \includegraphics[width=1\textwidth]{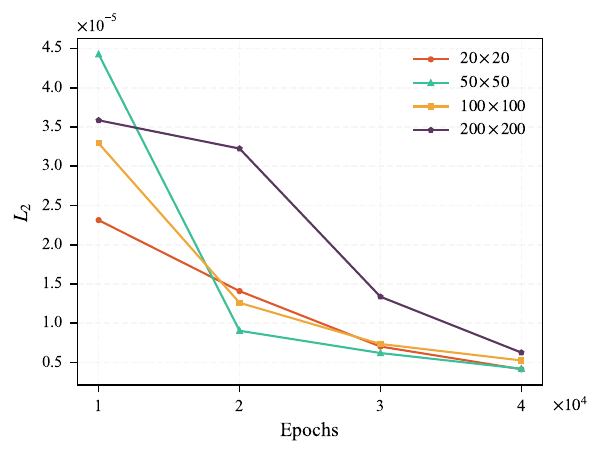}
        \caption*{(b)}
    \end{minipage}
    \begin{minipage}{0.195\textwidth}
        \centering
        \includegraphics[width=1\textwidth]{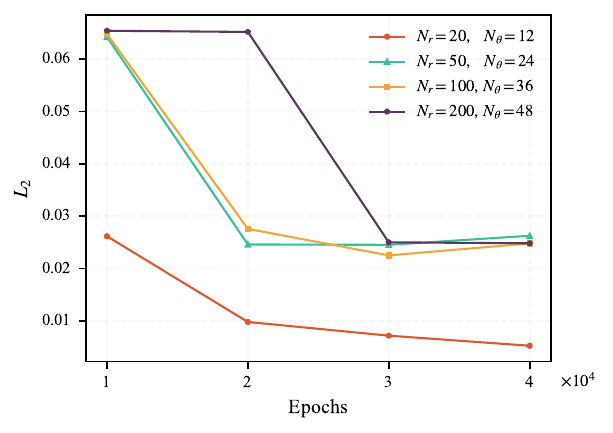}
        \caption*{(c)}
    \end{minipage}
    \begin{minipage}{0.195\textwidth}
        \centering
        \includegraphics[width=1\textwidth]{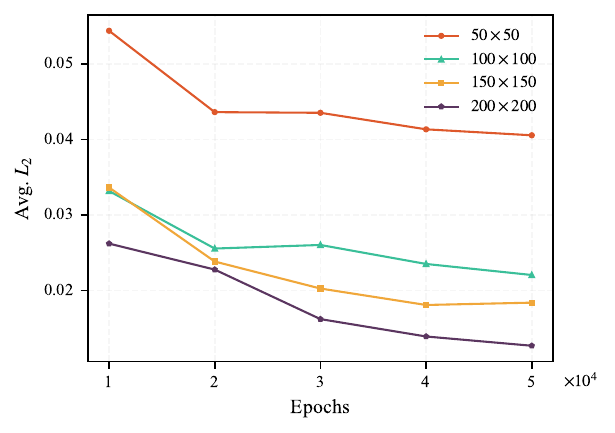}
        \caption*{(d)}
    \end{minipage}
    \begin{minipage}{0.195\textwidth}
        \centering
        \includegraphics[width=1\textwidth]{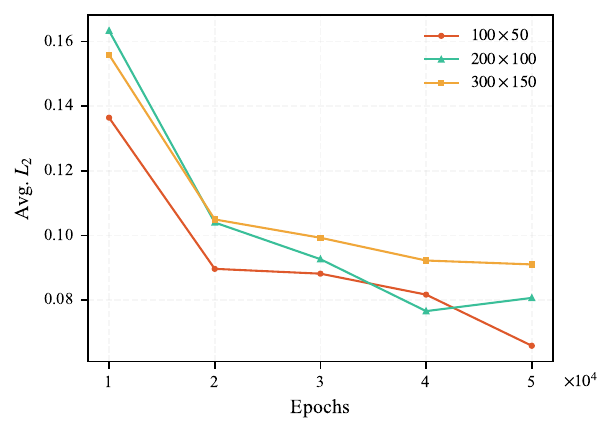}
        \caption*{(e)}
    \end{minipage}
    \caption{Hyperparameter analysis of training epochs under different resolutions: (a)–(e) for Sine, Polynomial, Liouville, Lid-driven cavity, and backward-facing step.}
    \label{fig:hyperparams_epochs}
\end{figure}

To investigate training duration, each test case at every resolution is trained up to $50000$ epochs, while the network architecture is kept identical to the setup in \cref{tab:training_params}; results are presented in \cref{fig:hyperparams_epochs}. In general, the largest performance improvement occurs within $10000$–$20000$ epochs. Further training can reduce the error slightly, yet the resulting gain becomes limited. Moreover, higher-resolution settings generally converge more slowly. Except for the Sine case, higher-resolution configurations cannot surpass the accuracy achieved at lower resolutions even with extended training. For incompressible NS cases, the decreasing averaged $L_2$ error originates mainly from ongoing pressure convergence, whereas velocity errors stagnate or even rise slightly. This indicates that the inferior accuracy at high resolution is not caused by insufficient training duration, but arises from inherent spectral bias induced by MSE-based optimization under high-resolution conditions. These observations provide practical guidance for engineering applications: priority should be given to network depth when allocating model parameters. Furthermore, for problems with increasing nonlinearity, we cannot rely on larger sampling sizes to improve performance, and lower-resolution configurations tend to be more robust.

\end{document}